\PassOptionsToPackage{unicode}{hyperref}
\PassOptionsToPackage{hyphens}{url}
\documentclass[12pt]{article}

\usepackage{amsmath,amssymb,amsbsy}
\usepackage{multirow}
\usepackage[T1]{fontenc}
\usepackage[utf8]{inputenc}
\usepackage{textcomp}
\usepackage{lmodern}
\usepackage{microtype}
\usepackage{xcolor}
\usepackage{booktabs}
\usepackage{longtable}
\usepackage{array}
\usepackage{calc}
\usepackage{graphicx}
\usepackage{subcaption}
\usepackage{float}
\usepackage{etoolbox}
\usepackage{bbm}

\usepackage[top=1in,bottom=1in,left=1.2in,right=1.2in]{geometry}

\usepackage{setspace}
\usepackage[authoryear,round]{natbib}
\usepackage{hyperref}
\hypersetup{
  colorlinks=true,
  linkcolor=blue,
  citecolor=blue,
  urlcolor=blue
}

\usepackage{amsthm}
\newcounter{thmcounter}

\newtheorem{theorem}[thmcounter]{Theorem}
\newtheorem{corollary}[thmcounter]{Corollary}
\newtheorem{remark}[thmcounter]{Remark}

\DeclareMathOperator*{\argmin}{arg\,min}

\usepackage{tikz}
\usetikzlibrary{calc,decorations.pathmorphing,shapes}
\newcounter{sarrow}
\newcommand\xrsquigarrow[2]{%
\stepcounter{sarrow}%
\mathrel{\begin{tikzpicture}[baseline={( $ (current bounding box.south) + (0, 1.5ex) $ )}]
  \node[inner sep=.5ex] (below\thesarrow) {$\scriptstyle #2$};
  \path[draw,<-,decorate,
    decoration={zigzag,amplitude=0.7pt,segment length=1.2mm,pre=lineto,pre length=4pt}]
      (below\thesarrow.north east) -- (below\thesarrow.north west);
  \node[above = 1.2ex, inner sep=.5ex] {$\scriptstyle #1$};
\end{tikzpicture}}%
}

\makeatletter
\def\maxwidth{\ifdim\Gin@nat@width>\linewidth\linewidth\else\Gin@nat@width\fi}
\def\maxheight{\ifdim\Gin@nat@height>\textheight\textheight\else\Gin@nat@height\fi}
\makeatother
\setkeys{Gin}{width=\maxwidth,height=\maxheight,keepaspectratio}
\makeatletter\def\fps@figure{htbp}\makeatother

\makeatletter
\ifx\paragraph\undefined\else
  \let\oldparagraph\paragraph
  \renewcommand{\paragraph}{\@ifstar\xxxParagraphStar\xxxParagraphNoStar}
  \newcommand{\xxxParagraphStar}[1]{\oldparagraph*{#1}\mbox{}}
  \newcommand{\xxxParagraphNoStar}[1]{\oldparagraph{#1}\mbox{}}
\fi
\ifx\subparagraph\undefined\else
  \let\oldsubparagraph\subparagraph
  \renewcommand{\subparagraph}{\@ifstar\xxxSubParagraphStar\xxxSubParagraphNoStar}
  \newcommand{\xxxSubParagraphStar}[1]{\oldsubparagraph*{#1}\mbox{}}
  \newcommand{\xxxSubParagraphNoStar}[1]{\oldsubparagraph{#1}\mbox{}}
\fi
\makeatother

\begin{document}

\begin{center}
  {\Large\bf Linear Regression Using Principal Components from\\[4pt]
    General Hilbert-Space-Valued Covariates}

  \bigskip
  Xinyi Li$^{a}$, Margaret Hoch$^{b}$, and Michael R.\ Kosorok$^{b}$

  \medskip
  {\small $^{a}$Clemson University, \;$^{b}$University of North Carolina at Chapel Hill}
\end{center}

\bigskip

\noindent\textbf{Abstract:}
We introduce Adaptive Subspace PCA (AS-PCA), a framework for principal component analysis of random elements in a general separable Hilbert space. AS-PCA projects the covariance operator onto a data-adaptive finite-dimensional subspace prior to eigendecomposition, requiring no kernel specification and accommodating multi-dimensional functional objects including images and surfaces. Under the second-moment condition, we prove a Donsker theorem for Hilbert-space-valued empirical processes and use it to establish uniform consistency and joint Gaussian limits for the leading eigenpairs. A data-driven diagnostic verifies projection accuracy, and a consistent proportion-of-variance-explained rule selects the number of components.
Building on AS-PCA, we construct Hilbert-Space Principal Component Regression (HS-PCR) for models combining Euclidean and Hilbert-space-valued covariates. The HS-PCR estimator is root-$n$ consistent and asymptotically normal, with an explicit influence function decomposition accounting for eigenfunction estimation uncertainty. Both nonparametric and wild bootstrap procedures are shown to be asymptotically valid.
Simulations with two- and three-dimensional imaging predictors confirm accurate eigenstructure recovery and nominal bootstrap coverage. HS-PCR is applied to Alzheimer's Disease Neuroimaging Initiative data in regression and precision-medicine settings.


\bigskip\noindent\textit{Address for correspondence:} Xinyi Li, School of Mathematical and Statistical
Sciences, Clemson University, Clemson, SC 29634, USA. Email:
\href{mailto:lixinyi@clemson.edu}{\texttt{lixinyi@clemson.edu}}

\newpage

\section{Introduction}\label{SEC:intro}

Multimodal learning, which integrates data from multiple sources to extract complementary information about complex systems, has become a dominant paradigm in biomedical research. In clinical settings, data routinely include medical imaging, genomic profiles, electronic health records, and biosignals; their joint analysis enables more comprehensive characterization of disease mechanisms, more accurate diagnosis, and more grounded guidance for personalized treatment. Alzheimer's disease (AD) exemplifies both the promise and the statistical challenges of this paradigm. As a leading cause of dementia, AD exhibits highly heterogeneous etiology shaped by neuroimaging biomarkers, genetic predisposition, and clinical risk factors. Pathophysiological changes detectable by neuroimaging can precede clinical dementia onset by more than a decade \citep{Bateman:etal:12}, making imaging biomarkers central to early intervention research. Understanding how these biomarkers, alongside genetic and environmental factors, jointly influence disease progression is a crucial yet statistically demanding task, owing to the high dimensionality and complex internal structure of imaging data.

Medical imaging data are naturally represented as functions defined over a spatial domain $\Omega\subset \mathbb{R}^d (d\geq 2)$, 
rendering functional data analysis (FDA) an efficient and natural modeling framework for neuroimaging. 
Recent advances have extended traditional FDA techniques to accommodate images and other complex functional data structures \citep{Zhu:Fan:Kong:14,Chen:Goldsmith:Ogden:19,Li:Wang:Wang:21,Zhang:etal:23}. 
In FDA, two principal perspectives have been employed to analyze functional data: (i) treating functional data as realizations of underlying stochastic processes and (ii) viewing them as random elements in Hilbert spaces. 
Under joint measurability assumptions, these perspectives are equivalent \citep{Hsing:Eubank:15}. 
The Hilbert-space perspective is particularly advantageous when the domain is multi-dimensional, irregular, or lacks a natural temporal ordering, as is typical of neuroimaging data. This motivates our focus on covariates that are random elements of a general separable Hilbert space.

We consider the model \vspace{-6pt}
\[
    Y = \alpha + \beta^{\top} X + \langle \gamma,Z \rangle + \varepsilon, \vspace{-6pt}
\] 
where $Y\in\mathbb{R}$ is the response, the $X\in\mathbb{R}^d$ is a Euclidean covariate vector, and $Z$ is a Hilbert space $\mathcal{H}$-valued random element, such as, for example, with images which can be viewed as a functional object defined over a domain $\Omega\subset \mathbb{R}^d (d\geq 2)$; $\varepsilon$ is conditionally mean zero given $X$ and $Z$; and $\alpha\in\mathbb{R}$, $\beta\in\mathbb{R}^d$, and $\gamma\in\mathcal{H}$ are unknown coefficients. 
The vector $\beta$ provides interpretable coefficients quantifying how the Euclidean covariates $X$ affect the response, while the contribution of $\mathcal{H}$-valued random element $Z$ is captured through $\gamma$.
Our objectives are twofold: to develop an interpretable, finite-dimensional functional principal component approximation to $\langle\gamma,Z\rangle$ that is computationally feasible in high-dimensional imaging settings, and to establish root-$n$ consistent estimation and inference procedures for all unknown parameters. Our framework extends traditional functional regression by incorporating an implicit, data-driven representation of the inner product structure within the Hilbert space. 
Existing principal component methods for functional data, however, are either restricted to one-dimensional temporal domains or lack the inferential theory needed for valid downstream regression analysis, motivating the development of a new, inference-ready PCA framework for general Hilbert-space-valued data.

Functional principal component analysis (FPCA) is the foundational dimension-reduction tool for functional data 
\citep{Morris:15,Wang:Chiou:Muller:16}. 
Classical asymptotic theory for empirical covariance operators in separable Hilbert spaces was established by \cite{Dauxois:etal:82} and substantially extended by \citet{Bosq:00} for function-space processes. For functional data on compact one-dimensional domains, \citet{Yuan:Cai:10} developed a reproducing kernel Hilbert space (RKHS) regularization framework that achieves minimax-optimal convergence rates under sparse and dense sampling regimes, while \citet{Cardot:etal:99,Cardot:etal:03} established asymptotic theory for functional linear regression based on estimated principal components. Further contributions include \citet{Hall:HosseiniNasab:06} on asymptotic properties of FPCA estimators, \citet{Dai:Muller:18} and \citet{Lin:Yao:19} on extensions to Riemannian manifold-valued data, \citet{Kim:etal:20} on autocovariance operators for Hilbert-space-valued processes, and \citet{Perry:etal:25} on inference procedures.

These foundational contributions, however, are predominantly confined to functional data indexed by one-dimensional temporal domains. Extensions to multi-dimensional functional data have received growing attention. \citet{Lila:Aston:Sangalli:16} developed PCA for functions on two-dimensional manifolds using finite element discretization, \citet{Happ:Greven:18} proposed multivariate FPCA for data on compact domains, and \citet{Shi:etal:22} constructed a spline-based two-dimensional FPCA framework for imaging feature extraction.
Despite these methodological advances, existing approaches for multi-dimensional functional data either lack comprehensive theoretical development or fail to provide a unified framework that ensures both efficient covariance operator estimation and optimal basis recovery.  
In particular, no existing framework simultaneously delivers consistent eigenfunction estimation, joint Gaussian limits for leading eigenpairs, and valid inference for regression parameters when the covariate is a multi-dimensional functional object.

A parallel and highly sophisticated body of work by Koltchinskii and collaborators \citep{Koltchinskii:Lounici:16, Koltchinskii:Lounici:17:Bernoulli, Koltchinskii:Lounici:17:Sankhya, Koltchinskii:Lounici:17:AOS, Koltchinskii:Loffler:Nickl:20} has established a comprehensive statistical theory for infinite-dimensional PCA in abstract Hilbert spaces. This line of work introduces the effective rank as a complexity measure and derives non-asymptotic concentration bounds, asymptotic distributions, and semiparametric efficiency results for spectral projectors and eigenvector functionals under minimal structural assumptions, with no requirement of Gaussianity or finite effective rank. These are elegant and far-reaching results. However, these theories are developed exclusively for the empirical covariance operator, which is the object obtained by forming sample outer products directly in the ambient Hilbert space, and does not address the practical problem of estimating the covariance operator when its reproducing kernel is unknown and must itself be estimated from data. Consequently, this framework does not extend naturally to settings where efficient finite-dimensional representation of the covariance operator is required for multi-resolution analysis, nor does it provide an integrated regression pipeline for predicting scalar outcomes from Hilbert-space-valued covariates.

Kernel-based approaches such as \citet{Yuan:Cai:10} address covariance estimation for functional regression, but they impose smoothness assumptions encoded through a predetermined reproducing kernel, require explicit specification of this kernel, and rely on predetermined basis truncation rather than data-driven selection. When the Hilbert space is abstract, as with multi-dimensional images, surfaces, or multi-index random fields, the covariance kernel is unknown, and RKHS methods do not directly apply.

What is missing from the existing literature is a PCA framework for general Hilbert-space-valued data that simultaneously: (i) requires no kernel specification; (ii) controls projection error with a verifiable, data-driven condition; (iii) establishes valid distributional theory under minimal moment assumptions; and (iv) supports a complete regression and bootstrap inference pipeline. We develop such a framework, Adaptive Subspace PCA (AS-PCA), and demonstrate its application in the regression setting of Hilbert-Space Principal Component Regression (HS-PCR).
This paper closes these gaps by developing a unified methodological and theoretical framework for principal component regression with general Hilbert-space-valued covariates. Our framework makes three distinct contributions.

First, we develop a framework for principal component analysis of random elements in a general separable Hilbert space. AS-PCA constructs a data-adaptive finite-dimensional subspace capturing nearly all variation of $Z$, estimates the covariance operator within that subspace, and recovers eigenpairs by standard matrix eigendecomposition, without requiring specification of a reproducing kernel or explicit covariance function. A data-driven diagnostic statistic allows practitioners to verify the projection accuracy condition from the observed data, and the number of retained components is selected by a consistent proportion-of-variance-explained rule.

Our second contribution is inferential theory under minimal conditions. We establish a Donsker theorem for empirical processes indexed by bounded bilinear forms of Hilbert-space-valued random elements, requiring only finite second moments: no Gaussianity, no sub-Gaussian tails, no effective rank constraints. This result drives uniform consistency and joint Gaussian limits for the leading estimated eigenvalues and eigenfunctions of AS-PCA (Theorems \ref{THM:op}--\ref{THM:asymp}), and validity of nonparametric and wild bootstrap procedures for uncertainty quantification (Theorems \ref{THM:bootstrap_v}--\ref{THM:bootstrap}). To our knowledge, this is the first complete inferential theory for PCA in general separable Hilbert spaces established under such minimal distributional conditions.

Third, we apply AS-PCA to construct a regression framework combining Euclidean and Hilbert-space-valued covariates. HS-PCR yields root-$n$ consistent and asymptotically normal estimators of all regression parameters (Theorem \ref{THM:LM_Consistency}), with an explicit influence function decomposition separating the contribution of eigenfunction estimation from the standard regression influence function. Bootstrap confidence regions are valid for all parameters simultaneously. HS-PCR extends classical functional principal component regression \citep[FPCR; ][]{Cardot:etal:99, Ramsay:Silverman:05} from one-dimensional functional covariates to abstract Hilbert-space-valued observations, with full inferential guarantees.

To handle the irregular anatomical boundaries and spatially heterogeneous structure of neuroimaging data, we employ multivariate splines over triangulations \citep{Lai:Schumaker:07} as the projection basis, which adapt naturally to complex spatial domains and have well-established approximation-theoretic properties ensuring the projection accuracy condition required by our theory. Simulation studies with two- and three-dimensional imaging predictors demonstrate accurate recovery of eigenstructures and regression coefficients, and bootstrap confidence intervals with coverage close to nominal levels across all settings studied. The methodology is applied to neuroimaging data from the Alzheimer's Disease Neuroimaging Initiative (ADNI), where we jointly analyze positron emission tomography (PET) imaging biomarkers, genetic risk factors, and demographic covariates to draw inferences about cognitive decline, and extend the framework to a precision-medicine formulation incorporating imaging-by-treatment interactions.

The remainder of the paper is organized as follows. Section \ref{SEC:HilbertSp} develops AS-PCA: the KL expansion and Donsker theorem (Section \ref{SUBSEC:SepHilbertSpace}), the adaptive subspace estimation procedure (Section \ref{SUBSEC:Eigen}), consistency and Gaussian limits for eigenpairs (Section \ref{SUBSEC:Asymp}), the projection accuracy diagnostic (Section \ref{SUBSEC:A2}), and the PVE component selection rule (Section \ref{SUBSEC:PVE}). Section \ref{SEC:LM} introduces HS-PCR: the regression estimator and its asymptotic linearity (Section \ref{SUBSEC:LMEST}), and bootstrap validity (Section \ref{SUBSEC:Bootstrap}). Section \ref{SEC:NeuroImple} details computational implementation for neuroimaging. Sections \ref{SEC:Simulation} and \ref{SEC:DA} present simulation studies and the ADNI data application. Section \ref{SEC:Conclusions} discusses extensions.

\section{Adaptive subspace PCA for Hilbert-space-valued random elements}
\label{SEC:HilbertSp}

In this section, we develop Adaptive Subspace PCA (AS-PCA), a framework for principal component analysis of random elements in a general separable Hilbert space. AS-PCA proceeds in three stages: adaptive subspace construction, empirical covariance estimation within the subspace, and eigendecomposition; and is supported by a complete inferential theory developed in Sections \ref{SUBSEC:SepHilbertSpace}--\ref{SUBSEC:PVE}.

\subsection{Hilbert-Space Framework: KL Expansion and Donsker Theorem}
\label{SUBSEC:SepHilbertSpace}

We begin with notations and basic assumptions that underpin our theoretical framework.
Let $\mathbb{P}_n$ denote the empirical measure based on $n$ observations, $P$ the expectation over a single generic observation,
and define the random measure $\mathbb{G}_n = \sqrt{n}(\mathbb{P}_n - P)$. Let $\mathbb{G}$ be a mean-zero Gaussian generalized Brownian bridge process indexed by a function class $\mathcal{F}$, with covariance $P(fg)-(Pf)(Pg)$, for all $f,g\in\mathcal{F}$.
We write $\mathbb{E}$ for the expectations of both Euclidean and Hilbert-space-valued random variables; for $\mathcal{H}$-valued objects, this is the Bochner expectation. 

Our first assumption concerns the basic structure of the random variable $Z$ taking values in the Hilbert space $\mathcal{H}$: 
\newline
\noindent (A1) (Square-integrable Hilbert-space setting) Let $\mathcal{H}$ be a Hilbert space with inner product $\langle\cdot, \cdot\rangle$ and norm $\|\cdot\|$, and assume the random variable $Z \in \mathcal{H}$ is separable with $\mathbb{E}\|Z\|^2<\infty.$
\newline
\noindent This assumption ensures that the covariance operator is well-defined and enables the Karhunen-Lo\`{e}ve (KL) expansion presented in Theorem \ref{THM:expansion} below.
We denote the mean of $Z$ by $\mu = \mathbb{E} Z$.

The following theorem is the classical Karhunen-Lo\`{e}ve expansion, extended with a summability condition that follows from Assumption (A1). 
While this result is well-established in the literature, we include it here, without proof, for completeness and to fix notation for subsequent development.

\begin{theorem}
\label{THM:expansion} (Karhunen-Lo\`{e}ve expansion)
Under Assumption (A1),
\[
	Z = \mu + \sum_{j=1}^{\infty} \lambda_j^{1/2} U_j\phi_j,
\]
where $\{\lambda_j\}_j$ are scalars with $\infty > \lambda_1 \geq \lambda_2 \geq \cdots$, $\{\phi_j\}_j$ form an orthonormal basis in $\mathcal{H}$, $\{U_j\}_j$ are mean zero and variance one random variables with $\mathbb{E}(U_jU_{j'}) = 0$ for $1 \leq j \ne j' \leq \infty$, and $\mu$ is finite. 
Moreover, $\lambda_j^{1/2} U_j = \langle \phi_j, Z - \mu \rangle$ and $\sum_{j=1}^\infty \lambda_j < \infty$. The span of the bases $\{\phi_j\}$ is a separable Hilbert space $\mathcal{H}_0 \subset \mathcal{H}$, with $\Pr(Z \in \mathcal{H}_0) = 1.$
\end{theorem}

Theorem \ref{THM:expansion} characterizes the structure of $Z$ and provides the foundation for developing our asymptotic framework. Building on this characterization, we now establish a Donsker-type result for Hilbert-valued random elements under mild moment conditions. This result constitutes our primary theoretical contribution in this section and has independent theoretical interest beyond its application to the estimators developed subsequently. Proofs are given in the supplemental material.

For notational simplicity, we state our main result under a zero-mean assumption, leaving the more general non-centered version to Corollary \ref{COR2:product}. This foundational result is essential for establishing the asymptotic normality of our proposed estimators later.

\begin{theorem}
\label{THM:product} (Donsker theorem for Hilbert-valued random elements)
Let Hilbert spaces $\mathcal{H}_k$ equipped with inner products $\langle\cdot, \cdot\rangle_k$ and norms $\|\cdot\|_k$ and random variables $Z_k \in \mathcal{H}_k$ satisfy Assumption (A1), for $1 \leq k \leq K$, where $K < \infty$. Assume also that $\mathbb{E} Z_k = 0$ and $\mathbb{E} (\prod_{k=1}^K \| Z_k \|_k^2) < \infty$, and let $\mathcal{B}_k = \{h \in \mathcal{H}_k : \|h\|_k \leq 1\}$, $k = 1, \cdots, K$. Then $\mathcal{F} = \{f(Z)=\prod_{k=1}^K \langle h_k, Z_k\rangle_k : h_k \in \mathcal{B}_k, 1 \leq k \leq K\}$ is Donsker.
\end{theorem}

\begin{remark}
Theorem \ref{THM:product} established a Donsker-type result for a rich class of functionals of Hilbert-space-valued random elements. Specifically, for any functional $f$ of the form $f(Z)=\prod_{k=1}^K \langle h_k, Z_k\rangle_k$ with $h_k$ in the unit ball, the corresponding empirical process $\sqrt{n}(\mathbb{P}_n f - P f)$ converges in distribution to a mean-zero Gaussian limit. This is a functional central limit theorem that guarantees all such empirical averages behave asymptotically as a tight Gaussian process. 
Theorem \ref{THM:product} is the theoretical engine of AS-PCA: because the empirical covariance operator is built from functionals of the form $\langle a, Z \rangle$ with $K = 2$, the Donsker property established here implies that all subsequent consistency and distributional results for AS-PCA hold under the same second-moment condition (A1), without any additional distributional assumptions.
This result guarantees that the empirical covariance operator constructed in Section \ref{COR1:product} behaves asymptotically as a tight Gaussian perturbation of the population operator, enabling the uniform convergence and distributional results of Sections \ref{SUBSEC:Asymp}--\ref{SUBSEC:PVE} and the regression inference of Section \ref{SEC:LM}.
\end{remark}

\begin{corollary}\label{COR1:product} (Single Hilbert space)
	Let $\mathcal{H}$ equip with inner product $\langle \cdot, \cdot \rangle$ and norm $\|\cdot \|$ and $Z \in \mathcal{H}$ satisfy Assumption (A1). Then $\mathcal{F} = \{f(Z)=\langle h, Z \rangle : \|h \| \leq 1, h \in \mathcal{H}\}$ is Donsker. 
\end{corollary}
Corollary \ref{COR1:product} is a special case of Theorem \ref{THM:product} by setting $K=1$. And the following Corollary \ref{COR2:product} extends Theorem \ref{THM:product} by removing the centering assumption. The moment condition ensures the necessary uniform integrability after centering.
\begin{corollary}\label{COR2:product} (Nonzero mean)
	Let $\mathcal{H}_k$ equip with inner product $\langle\cdot, \cdot\rangle_k$ and norm $\|\cdot\|_k$, and random variable $Z_k$ satisfy Assumption (A1), for each $1 \leq k \leq K$ for some $K < \infty$. Let $\mu_k = \mathbb{E} Z_k$, $1 \leq k \leq K$, and assume that $\sum_{k=1}^K \sum_{r \in N_k} \mathbb{E} (\prod_{k^{\prime} = 1}^K \|Z_{k^\prime}-\mu_{k^\prime} \|_{k'}^{2r_{k^{\prime}}}) < \infty$, where $\mathcal{N}_k =\{r = (r_1, \cdots, r_K) \in \{0,1\}^K : \sum_{k' = 1}^K r_{k'} = k\}$. 
	Let $\mathcal{B}_k = \{h \in \mathcal{H}_k : \|h\|_k \leq 1\}$, $k = 1, \cdots, K$.
	Then $\mathcal{F} = \{f(Z)=\prod_{k=1}^K \langle h_k, Z_k\rangle_k : h_k \in \mathcal{B}_k, 1 \leq k \leq K\}$ is Donsker.
\end{corollary}

\subsection{The AS-PCA Estimation Procedure}
\label{SUBSEC:Eigen}

AS-PCA estimates the principal components of $Z$ through three stages. We describe each stage in turn, followed by a remark on the statistical properties of the construction.

\noindent {\bf AS-PCA Stage 1: Adaptive Subspace Construction}

\noindent Let $\Psi_N^\ast = (\psi_1^\ast, \cdots, \psi_N^\ast)^\top$ 
be a finite basis (not necessarily orthonormal) of the projection subspace $\mathrm{span}\{\psi_1^\ast, \cdots, \psi_N^\ast\}\subset \mathcal{H}$. We allow the dimension $N\equiv N(n)$ to grow with the sample size and even to be larger than $n$. 
Define the $N \times N$ Gram matrix $L=(\langle \psi_\ell^\ast, \psi_{\ell'}^\ast \rangle)_{1\leq \ell,\ell'\leq N}$.
Let $\Lambda$ be the diagonal matrix of the non-zero eigenvalues of $L$, and let $\Gamma$ be the matrix of the corresponding eigenvectors. Set $L_{\ast}^{-1/2} = \Lambda^{-1/2}\Gamma^{\top}$, so that $\psi_{\ell} = \sum_{\ell'=1}^N [L_{\ast}^{-1/2}]_{\ell, \ell'}\psi_{\ell'}^{\ast}$, $l=1,\ldots,N$, from an orthonormal basis spanning the the same subspace as $\Psi_N^\ast$.
We denote by $\mathrm{G}_N x=\sum_{\ell=1}^N\langle x,\psi_l \rangle \psi_l$ the orthogonal projector onto this subspace for any $x\in\mathcal{H}$.
While the orthonormal basis $\{\psi_{\ell}\}$ is introduced for theoretical convenience, in applications one may work directly with the original basis $\Psi_N^\ast$ and the Gram matrix $L$; all quantities can be computed from $\Psi_N^\ast$ without explicitly forming $\{\psi_{\ell}\}$, which is a practical advantage of the proposed construction.

With a slight abuse of notation, for any element $a \in \mathcal{H}$, we write $a^\top \equiv \langle a, \cdot\rangle$, so that for any $a_1, a_2 \in \mathcal{H}$, we have $a_1^\top a_2 = \langle a_1, a_2 \rangle$ for their inner product.
We also write $a_1 a_2^\top$ for the rank-one operator $x \mapsto \langle x,a_2\rangle a_1$, $x\in\mathcal{H}$.
In particular, $\mathrm{G}_N=\sum_{\ell=1}^N \psi_{\ell}\psi_{\ell}^{\top}$.
For finite matrices, $(\cdot)^{\top}$ denotes the usual transpose.

To control the error introduced by projection, we impose the following condition:
\newline
(A2) (Projection accuracy) 
	$\delta_n \equiv \mathbb{E}\|(Z-\mu) - \mathrm{G}_N (Z-\mu)\|^2=o(n^{-1})$.
\newline
Assumption (A2) requires that the subspace spanned by $\{\psi_{\ell}\}_{\ell=1}^N$ captures almost all variation of $Z$, with the part it misses being negligible compared to the $n^{-1/2}$ stochastic error. If needed, we can enlarge the basis $\Psi_N^\ast$ until (A2) holds; we also describe a way to test this assumption later in this section.
Specific examples of the candidate bases are given in Section \ref{SEC:NeuroImple}.

The population covariance operator $V_0 = \mathbb{E}\{(Z-\mu)(Z-\mu)^\top\}$ acts on the entire (possibly infinite-dimensional) space $\mathcal{H}$, but from $n$ observations we can estimate only a finite-rank operator. By first projecting the data onto the finite-dimensional subspace spanned by $\mathrm{G}_N$, we ensure that the population covariance operator we want to estimate and the sample object we construct are on the same finite-dimensional space. Assumption (A2) guarantees that the projection error is negligible compared to the statistical error, so working in $\mathrm{span}\{\psi_{\ell}\}$ does not change the asymptotic conclusions but makes the estimation problem finite-dimensional.

\noindent {\bf AS-PCA Stage 2: Subspace Covariance Estimation}

\noindent Let $Z_1,\ldots, Z_n$ be independent and identically distributed copies of $Z$, and let $\bar{Z}_n=\mathbb{P}_nZ$ be the sample mean.
Form the $N\times N$ empirical moment matrix 
$M_n^\ast = n^{-1} \sum_{i=1}^n (\langle \psi_\ell^\ast, Z_i - \bar{Z}_n\rangle \langle Z_i - \bar{Z}_n, \psi_{\ell'}^\ast\rangle)_{1 \leq \ell, \ell' \leq N}$.
Note that $\langle Z_i - \bar{Z}_n, \psi_{\ell}^\ast\rangle$ is the coefficient of $Z_i - \bar{Z}_n$ in terms of the (possibly non-orthonormal) basis $\psi_{\ell}^\ast$, and thus $M_n^\ast$ is the sample covariance matrix of the projected data in the original (possibly non-orthonormal) $\Psi_N^\ast$. 
To obtain the covariance in the orthonormalized space spanned by $\{\psi_{\ell}\}$, we transform $M_n^\ast$ by $L_{\ast}^{-1/2}$ to obtain $\widehat{M}_n = L_{\ast}^{-1/2}M_n^\ast (L_{\ast}^{-1/2})^\top$.
By construction, $\widehat{M}_n$ is the empirical covariance matrix of the coefficient of $Z_i - \bar{Z}_n$ with respect to the orthonormal projector $\mathrm{G}_N$, and is therefore symmetric and positive semi-definite.
In particular,
\[
	\widehat{M}_n 
	= n^{-1} \sum_{i=1}^n \left(\langle \psi_\ell, Z_i - \bar{Z}_n\rangle \langle Z_i - \bar{Z}_n, \psi_{\ell'}\rangle\right)_{1 \leq \ell, \ell' \leq N}.
\]

\noindent {\bf AS-PCA Stage 3: Eigendecomposition and PC Recovery}

\noindent 
Let $\widehat{M}_n = \sum_{j = 1}^{J} \widehat{\lambda}_j \omega_j \omega_j^\top$ be the eigen-decomposition of $\widehat{M}_n$, where the rank $J\leq \min\{n,N\}$,
$\infty > \widehat{\lambda}_1 \geq \widehat{\lambda}_2 \geq \cdots \geq \widehat{\lambda}_{J} \geq 0$, 
and $\omega_{1}, \cdots, \omega_{J}\in\mathbb{R}^N$ are orthonormal eigenvectors. 
Each $\omega_j=(\omega_{j1}, \ldots,\omega_{jN})^{\top}$ defines an element of $\mathcal{H}$ by $\widehat{\phi}_{j} = \sum_{\ell =1}^N \omega_{j \ell}\psi_\ell$, $j = 1, \cdots, J$.
Because $\{\phi_{\ell}\}$ is orthonormal and the $\omega_j$ are orthonormal in $\mathbb{R}^N$, the functions $\{\widehat{\phi}_j\}_{j=1}^J$ are orthonormal in $\mathcal{H}$, which are our estimated functional principal bases, and $\{\widehat{\lambda}_j\}_{j=1}^J$ are the associated estimated eigenvalues.

It is convenient to rewrite the construction in operator form. Define the empirical covariance operator $\widehat{V}_n: \mathcal{H} \rightarrow \mathcal{H}$ by
\[
	\widehat{V}_n 
	= \frac{1}{n} \sum_{i=1}^n \mathrm{G}_N^\top (Z_i - \bar{Z}_n)(Z_i - \bar{Z}_n)^\top \mathrm{G}_N
	= \frac{1}{n}\sum_{\ell=1}^N \sum_{\ell' =1}^N \sum_{i =1}^n \psi_\ell \langle \psi_\ell, Z_i - \bar{Z}_n\rangle \langle Z_i - \bar{Z}_n, \psi_{\ell'}\rangle \psi_{\ell'}^\top.
\]
That is, $\widehat{V}_n$ first projects a vector onto the subspace $\mathrm{span}\{\psi_{\ell}\}$ and then applies the empirical covariance in that subspace. With this notation, we can represent $\widehat{V}_n$ as $\widehat{V}_n = \sum_{j=1}^{J} \widehat{\lambda}_j \widehat{\phi}_j\widehat{\phi}_j^\top$, and equivalently, $\widehat{V}_n \widehat{\phi}_j = \widehat{\lambda}_j \widehat{\phi}_j$, $j=1,\ldots,J$.
Hence $\widehat{V}_n$ is a valid empirical covariance operator that we could compare with $V_0$.

\begin{remark}(Adaptive Subspace PCA)
The three-stage procedure above constitutes AS-PCA. The term ``adaptive subspace'' reflects three distinct forms of adaptivity: 
(i) the projection subspace $\mathrm{span}\{\psi_{\ell}\}_{\ell=1}^N$ is selected to match the geometry of the domain $\Omega$ and the structure of $Z$, rather than being predetermined by a kernel; 
(ii) the dimension $N \equiv N(n)$ grows with the sample size to ensure projection accuracy (A2), verified empirically via $\widehat{T}_n$ in Section \ref{SUBSEC:A2}; and 
(iii) the number of retained components $m$ is chosen from the data by the PVE rule of Section \ref{SUBSEC:PVE}. 
Unlike kernel-based FPCA \citep{Yuan:Cai:10}, AS-PCA requires no specification of a reproducing kernel; unlike spectral methods applied directly to the empirical covariance operator \citep{Koltchinskii:Lounici:17:AOS}, it produces a finite-rank estimator directly comparable to the population operator and amenable to the regression framework of Section \ref{SEC:LM} . When $\mathcal{H} = L^2[0,1]$ and $\Psi_N^{\ast}$ is a B-spline basis, AS-PCA reduces to the basis-expansion approach of \citet{Ramsay:Silverman:05}, confirming that the latter is a one-dimensional special case of the general framework.
\end{remark} \vspace{-6pt}

We have thus obtained the empirical covariance operator $\widehat{V}_n$ together with its eigenpairs $(\widehat{\lambda}_j, \widehat{\phi}_j)$; in the next Section \ref{SUBSEC:Asymp}, we study their large-sample behavior, using the Donsker theorem of Section \ref{SUBSEC:SepHilbertSpace}.

\subsection{Asymptotic Theory for AS-PCA: Consistency and Gaussian Limits}
\label{SUBSEC:Asymp}

This section establishes the large-sample theory for AS-PCA, showing that the adaptive subspace construction of Section \ref{SUBSEC:Eigen} yields consistent and asymptotically Gaussian estimates of the leading eigenvalues and eigenfunctions of the population covariance operator. 
Proofs are provided in the supplemental material.
Throughout, we maintain Assumptions (A1)--(A2).

We compare operators as bounded bilinear forms on the unit ball $\mathcal{B}\subset\mathcal{H}$. For a bounded linear operator $V:\mathcal{H}\to\mathcal{H}$, define $\| V \|_{\mathcal{B} \times \mathcal{B}} = \sup_{a_1, a_2 \in \mathcal{B}} |a_1^\top V a_2|$. In this sense, $V_0$, $\widehat{V}_n$, and similar style operators are elements of $\ell^{\infty}(\mathcal{B} \times \mathcal{B})$.
When measurability of the indexing class cannot be guaranteed, we use outer probability and expectation \citep{vanderVaart:Wellner:23, Kosorok:08} and write $X_n \xrightarrow{as\ast} X$ for outer almost sure convergence, meaning there exists a sequence $\{\Delta_n\}$ of measurable random variables satisfying $\| X_n - X \| \leq \Delta_n$ for all $n$ and $\Pr(\limsup_{n \rightarrow \infty} \Delta_n = 0) = 1$.

\begin{theorem}\label{THM:op} (Uniform convergence and Gaussian limit of $\widehat{V}_n$)
Assume (A1)--(A2). Then
\begin{itemize}
	\item[(a)] $\| \widehat{V}_n -V_0\|_{\mathcal{B} \times \mathcal{B}} \xrightarrow{as\ast} 0$.
	\item[(b)] Moreover, provided that $\mathbb{E}\|Z\|^4 < \infty$, we also have the following:
	\begin{itemize}
		\item[(i)] $\|\sqrt{n}(\widehat{V}_n -V_0) - \mathbb{G}_n\left((Z-\mu)(Z-\mu)^\top\right)\|_{\mathcal{B} \times \mathcal{B}} = o_P(1)$.
		\item[(ii)] $\mathcal{F}_2 = \{a_1^\top (Z-\mu)(Z-\mu)^\top a_2 : a_1, a_2 \in \mathcal{B}\}$ is Donsker. 
		\item[(iii)] $\|\sqrt{n}(\widehat{V}_n - V_0)\|_{\mathcal{B} \times \mathcal{B}} = O_P(1)$.
	\end{itemize}
\end{itemize}
\end{theorem}

\begin{remark}
    Part (a) establishes that the empirical covariance operator on the projected space converges uniformly to the population covariance as a bilinear form. Part (b) provides a functional central limit theorem: 
\[
	\sqrt{n}(\widehat{V}_n -V_0) 
	= \mathbb{G}_n\left((Z-\mu)(Z-\mu)^\top\right) + o_P(1)~~~~\text{in } \ell^{\infty}(\mathcal{B} \times \mathcal{B}).
\]
The Donsker property established in Theorem \ref{THM:product} ensures tightness and weak convergence to a Gaussian limit for this class.
\end{remark} \vspace{-6pt}

Since $V_0$ is self-adjoint, positive, and compact, its spectrum is a (possibly infinite) sequence of nonnegative eigenvalues $\{\lambda_j\}_j$ with finite multiplicities and only possible accumulation at zero. Ties (equal eigenvalues) are allowed.
To ensure identifiability of individual principal components, we impose a strict spectral gap assumption.
\newline
	(A3) (Spectral gap) For some $1 \leq m < \infty$, $\infty > \lambda_1 > \lambda_2 > \cdots > \lambda_m > \lambda_{m+1} \geq 0$. 
\newline
This condition ensures the first $m$ eigenfunctions $\{\phi_1,\ldots,\phi_m\}$ are uniquely defined up to sign. Accordingly, we can assume without loss of generality that differences between eigenfunctions are chosen to minimize norm distance, with random sign selection in case of a tie. If $\lambda_j=\cdots=\lambda_{j+r-1}$ has multiplicity $r>1$, then $\widehat{\phi}_j,\ldots,\widehat{\phi}_{j+r-1}$ consistently estimate the eigenspace $\text{span}\{\phi_j,\ldots,\phi_{j+r-1}\}$. Results for individual eigenfunctions require Assumption (A3); without it, subspace consistency and $\sqrt{n}$ limits only hold for the orthogonal projector onto the eigenspace.
Relaxing (A3) (e.g., handling eigenvalue multiplicities) is possible with additional structure but lies beyond our scope.

\begin{theorem}
\label{THM:maxes}
(Consistency of PCs)
Under (A1)--(A3),
$$
	\max_{1 \leq j \leq m+1} | \widehat{\lambda}_j - \lambda_j | \xrightarrow{as\ast} 0, ~~~~
	\max_{1 \leq j \leq m} \|\widehat{\phi}_j - \phi_j\| 
	\xrightarrow{as\ast} 0.
$$
\end{theorem}

Our final result establishes joint asymptotic linearity and asymptotic normality of the estimated eigenvalues and eigenfunctions.
\begin{theorem}\label{THM:asymp} (Joint Gaussian $\sqrt{n}$-limit for leading eigenpairs)
Assume (A1)--(A3) and $\mathbb{E}\|Z\|^4 < \infty$. Then, for $j = 1, \cdots, m$,  
\begin{align*}
	\sqrt{n}(\widehat{\lambda}_j - \lambda_j) 
	&= \phi_j^\top \mathbb{G}_n\left((Z-\mu) (Z-\mu)^\top\right)\phi_j + o_P(1)\\ 
		\sqrt{n}(\widehat{\phi}_j - \phi_j) 
	&= \!\!\!\!\sum_{j' \neq j : j' \geq 1}\!\! \frac{\phi_{j'}}{\lambda_j - \lambda_{j'}}\phi_{j'}^\top \mathbb{G}_n\left((Z-\mu) (Z-\mu)^\top\right)\phi_j +o_P(1).
\end{align*}
Moreover, all of these quantities on the left are jointly, asymptotically a tight, mean-zero Gaussian process with covariance defined by the covariance generated by the influence functions on the right, which are all Donsker classes. 
\end{theorem}

\begin{remark}
Let $\xi_j=\langle Z-\mu,\phi_j\rangle$. The expansions in Theorem \ref{THM:asymp} show
\[
	\sqrt{n}(\widehat{\lambda}_j - \lambda_j)
	=\mathbb{G}_n(\xi_j^2-\lambda_j) + o_P(1), ~~~~
		\sqrt{n}(\widehat{\phi}_j - \phi_j)
	=\sum_{j' \neq j : j' \geq 1} \frac{\mathbb{G}_n(\xi_{j'} \xi_j)}{\lambda_j-\lambda_{j'}} \phi_{j'} + o_P(1).
\]
Thus $\widehat{\lambda}_j$ depends only on the variance of the $j$-th score, while $\widehat{\phi}_j$ aggregates cross-products with all other components, weighted by the inverse spectral gaps $(\lambda_j-\lambda_{j'})^{-1}$.  Smaller gaps lead to larger sampling variability for $\widehat{\phi}_j$. Because these influence functions form a Donsker class, the eigenvalue and eigenfunction estimators admit a joint tight, mean-zero Gaussian limit.
\end{remark}

Theorems \ref{THM:op}--\ref{THM:asymp} together constitute the inferential theory of AS-PCA: they establish that the adaptive subspace construction introduces no asymptotic bias, and that the estimated eigenpairs inherit the distributional structure of the population eigenpairs at the $\sqrt{n}$ rate.

\subsection{Assessing Projection Accuracy: The AS-PCA Diagnostic}
\label{SUBSEC:A2}

A distinctive feature of AS-PCA is that the projection accuracy condition (A2), which ensures the adaptive subspace captures nearly all variation of $Z$, can be assessed directly from the data before proceeding to eigendecomposition and regression.
We now develop a data-driven estimator of $\delta_n$ and construct a test statistic for assessing whether (A2) holds for a chosen projection subspace.

Define the empirical projection error as 
\begin{equation}
\label{DEF:deltan_hat}
	\widehat{\delta}_n = \frac{1}{n} \sum_{i=1}^n \| Z_i - \bar{Z}_n \|^2 - \sum_{\ell = 1}^N \frac{1}{n}\sum_{i=1}^n \langle \psi_\ell, Z_i - \bar{Z}_n\rangle^2,
\end{equation}
which measures the sample variance not captured by the projection subspace $\mathrm{span}\{\psi_1, \cdots,$ $\psi_N\}$. Denote the projection residual for each subject $i$ as $r_i=(Z_i - \bar{Z}_n) - G_N(Z_i - \bar{Z}_n)$. By orthogonality of $G_N$, it follows (with details provided in the supplement) that $\widehat{\delta}_n = n^{-1}\sum_{i=1}^n \|r_i\|^2$, demonstrating that the proposed estimator equals the average squared residual after projection across all subjects.
We then define the test statistics  
\begin{equation}
\label{DEF:Sn_Tn}
	\widehat{S}_n^2 = \frac{1}{n}\sum_{i=1}^n \left(\|r_i\|^2 - \widehat{\delta}_n\right)^2, ~~~~~~~~
	\widehat{T}_n = \frac{\sqrt{n}\widehat{\delta}_n}{\sqrt{\widehat{S}_n^2+\frac{1}{n}}}.
\end{equation}
The term $1/n$ stabilizes the denominator when $\widehat{S}_n^2$ is small, as may occur with homoscedastic samples or small sample sizes. 
Heuristically, $\widehat{T}_n$ functions as a $t$-statistic for testing whether the mean residual projection error $\delta_n$ equals zero.
The following theorem formalizes the theoretical properties of $\widehat{T}_n$.

\begin{theorem}\label{THM:delta} (Assessing (A2))
Assume (A1), the Gram matrix $L$ is positive definite, and that $\sup_{j \geq 1} \mathbb{E}(\langle \phi_j, Z-\mu \rangle^4/\mathbb{E}^2\langle \phi_j, Z-\mu \rangle^2) < \infty$. Then 
\begin{itemize}
	\item[(i)] if $\limsup_{n \rightarrow \infty} n\delta_n = 0$, then $\widehat{T}_n \xrightarrow{P} 0$;
	\item[(ii)] for any $k < \infty$ and $\epsilon >0 $, there exists a $k_1 < \infty$ such that 
	\[
		\liminf_{n \rightarrow \infty}n \delta_n \geq k_1~~\text{implies}~~
		\liminf_{n \rightarrow \infty} \Pr(\widehat{T}_n > k) \geq 1-\epsilon.
	\]
\end{itemize}
\end{theorem}

\begin{remark}
The moment condition in Theorem \ref{THM:delta} implies that $E\|Z-\mu\|^4<\infty$. 
Part (i) guarantees that when (A2) holds, $\widehat{T}_n$ converges to zero in probability. Part (ii) establishes that when (A2) is violated, $\widehat{T}_n$ eventually exceeds any fixed threshold with high probability.
Therefore, $\widehat{T}_n$ provides a valid statistic for assessing Assumption (A2).
Proofs are deferred to the appendix.
In practice, the AS-PCA diagnostic $\widehat{T}_n$ is computed as a byproduct of Stage 1, adding negligible computational cost. If $\widehat{T}_n$ exceeds the critical value $z_{1-\alpha}$, the subspace dimension $N$ should be increased until the diagnostic is satisfied, at which point AS-PCA Stages 2 and 3 proceed.
 \end{remark}

In practice, we implement the one-sided test for (A2) through the following procedure:\vspace{-8pt}
\begin{enumerate}
	\item Compute $\widehat{\delta}_n$, $\widehat{S}_n^2$, and $\widehat{T}_n$ from (\ref{DEF:deltan_hat})--(\ref{DEF:Sn_Tn}). \vspace{-8pt}
	\item For a nominal level $\alpha \in (0, 0.05]$, reject (A2) if $\widehat{T}_n>z_{1-\alpha}$, where $z_{1-\alpha}$ is the upper standard normal quantile.
\end{enumerate}

\subsection{Adaptive Component Selection in AS-PCA: The PVE Rule}
\label{SUBSEC:PVE}

The third and final form of adaptivity in AS-PCA is the selection of the number of retained components $m$. Following the proportion of variance explained (PVE) approach
is to require that they explain a specified proportion/percentage of the total variation \citep{Kong:Staicu:Maity:16}. 
Let $Z \in \mathcal{H}$ satisfy (A1). Since $\mathbb{E}\|Z-\mu\|^2 = \sum_{j=1}^\infty \lambda_j$, for any $\tau \in (0, 1)$, define the population target
$m_0(\tau) = \inf \{m \geq 1 : \sum_{j = 1}^m \lambda_j > (1-\tau)\mathbb{E}\|Z-\mu\|^2\}$.
Thus, $m_0(\tau)$ is the smallest number of eigenvalues required to explain at least $(1-\tau)$ of the total variance.

In practice, we adopt the sample analogue and propose the PVE rule 
$\widehat{m}_n(\tau) = \inf \{m \geq 1: \sum_{j=1}^m \widehat{\lambda}_j > (1-\tau)\mathbb{P}_n\|Z-\bar{Z}_n\|^2\}$, which selects the smallest $m$ for which the estimated eigenvalues account for at least $(1-\tau)$ of the empirical variance. For example, $\widehat{m}_n(0.95)$ is the minimal number of components needed to explain 95\% of the sample variance. 
The threshold $\tau$ is commonly chosen within $[0.95, 0.99]$. For our implementation, we set $\tau=0.95$.

\begin{corollary}
\label{COR:PVE_Consistency} (Consistency of the PVE rule)
Assume (A1).
Let $\mathrm{G}_N$ satisfy (A2), and suppose the spectral gap assumption (A3) holds at $m = m_0(\tau)$. Then $\widehat{m}_n(\tau) \xrightarrow{as\ast} m_0(\tau)$.
\end{corollary}

\begin{remark}
Theorem \ref{THM:maxes} also yields $\widehat{\lambda}_{m+1} \rightarrow \lambda_{m+1}$; however, we will also need consistency of the accompanying eigenfunctions to be able to estimate all principal components for downstream analyses.
\end{remark}

\section{Hilbert-Space Principal Component Regression (HS-PCR)}
\label{SEC:LM} 

We now apply AS-PCA to construct Hilbert-Space Principal Component Regression (HS-PCR), a regression framework that combines Euclidean covariates $X \in \mathbb{R}^d$ with Hilbert-space-valued covariates $Z \in \mathcal{H}$. HS-PCR uses the estimated principal components from AS-PCA as finite-dimensional summaries of $Z$, and inherits the inferential guarantees of AS-PCA, including root-$n$ consistency, asymptotic normality, and bootstrap validity, through the asymptotic linearity results established below.

\subsection{The HS-PCR Estimator}
\label{SUBSEC:LMEST}

We consider the following linear model with both Euclidean and Hilbert-space covariates.
Let $(Y_1,X_1,Z_1),\ldots, (Y_n,X_n,Z_n)$ be i.i.d., where $Y\in\mathbb{R}$ is the outcome of interest, $X\in\mathbb{R}^d$ is a finite-dimensional covariate, and $Z\in\mathcal{H}$ is the random element taking values in Hilbert space $\mathcal{H}$ with inner product $\langle\cdot,\cdot\rangle$ and norm $\|\cdot\|$.
We assume the model 
\[
	Y=\alpha+\beta^{\top}X+\langle \gamma, Z\rangle+\varepsilon,
\]
where $\alpha\in\mathbb{R}$, $\beta\in\mathbb{R}^d$, $\gamma\in\mathcal{H}$, and $\varepsilon$ has mean zero conditional on $(X,Z)$.
Our goal is to obtain $\sqrt{n}$-consistent estimation and inference for a finite-dimensional approximation of this model based on estimated functional principal components of $Z$. 

Recall that $\{\phi_j\}_{j\geq1}$ and $\{\lambda_j\}_{j\geq1}$ denote the eigenfunctions and eigenvalues of the covariance operator of $Z$, and $\{\widehat{\phi}_j\}_{j\geq1}$ and $\{\widehat{\lambda}_j\}_{j\geq1}$ denote their estimators constructed in Section \ref{SEC:HilbertSp}.
Fix $m$ and define the ``oracle'' principal component vector 
$Z_\ast = (\langle \phi_1, Z \rangle, \cdots,$ $\langle \phi_m, Z \rangle)^\top\in\mathbb{R}^m$, 
and its empirical analogue based on estimated eigenfunctions
$\widehat{Z}_\ast = (\langle \widehat{\phi}_1, Z \rangle, \cdots, \langle \widehat{\phi}_m, Z \rangle)^\top\in\mathbb{R}^m$.
In practice, we choose $m$ by the PVE rule $\widehat{m}_n(\tau)$ in Section \ref{SUBSEC:PVE} for a prespecified $\tau\in(0,1)$.

We form the combined ``oracle'' and empirical regressor vectors
$U = (1, X^\top, Z_\ast^\top)^\top\in \mathbb{R}^{p+1}$, 
$\widehat{U}=(1, X^\top,\widehat{Z}_{\ast}^\top)^\top \in \mathbb{R}^{p+1}$, respectively,
where $p=d+m$. 
Note that $(\widehat{U}-U) = (0, \cdots, 0, \langle \widehat{\phi}_1 -\phi_1, Z \rangle, \cdots, \langle \widehat{\phi}_m - \phi_m, Z \rangle)^\top$.
Assumptions (A1)--(A3) are imposed on $Z$, on the projection dimension used in Section \ref{SUBSEC:PVE}, and on the eigenstructure of the covariance operator, and we additionally assume $\mathbb{E}\|Z\|^4<\infty$.
For clarity, we assume there is a single $\mathcal{H}$-valued $Z$; extension to finitely many such random elements is straightforward but omitted.

Let $\mu_U = \mathbb{E}U$ and $\mu_Z = \mathbb{E}Z$. 
For the random vector $\widehat{U}$, it is convenient to use the empirical-process notation $P$ to denote the expectation operator with respect to $(X,Z)$, holding the estimated eigenfunctions $\widehat{\phi}_1, \cdots, \widehat{\phi}_m$ fixed.
In this notation, 
$\mu_{\widehat{U}} = P\widehat{U} = (1, PX^\top, P\widehat{Z}_\ast^\top)^\top$,
where $PX =\mathbb{E}X \equiv \mu_X$ and $P\widehat{Z}_\ast = (\langle \widehat{\phi}_1, \mathbb{E}Z \rangle ^\top, \cdots, \langle \widehat{\phi}_m, \mathbb{E}Z \rangle^\top)^\top$.
Denote the second-moment matrix of the regressor $\Sigma_0=P(UU^{\top})$.

Define the population target parameter $\theta_0=(\alpha_0, \beta_0^{\top}, \gamma_0^{\top})^{\top}$ as the minimizer of the mean squared error in the model using the true principal components: 
$\theta_0=\argmin_{\theta} P(Y-\theta^\top U)^2
	= \Sigma_0^{-1}P(UY)$.
In practice, we only observe $\widehat{Z}_{\ast}$.
The HS-PCR estimator replaces the oracle principal component scores $Z_{\ast} = (\langle \phi_1, Z\rangle, \ldots, \langle \phi_m, Z\rangle)^\top$ with the AS-PCA estimated scores $\widehat{Z}_{\ast} = (\langle \widehat{\phi}_1, Z\rangle, \ldots, \langle \widehat{\phi}_m, Z\rangle)^\top$, and fits ordinary least squares on the combined regressor $\widehat{U} = (1, X^\top, \widehat{Z}_*^\top)^\top$.
We therefore define the empirical estimator $\widehat{\theta}_n$ be 
\begin{equation} \label{EQN:ThetaHat}
	\widehat{\theta}_n
	=(\widehat{\alpha}_n^{\top}, \widehat{\beta}_n^{\top}, \widehat{\gamma}_n^{\top})^{\top}
	=\argmin_{\theta} \mathbb{P}_n(Y-\theta^{\top}\widehat{U})^2
	= \left\{\mathbb{P}_n \left(\widehat{U}\widehat{U}^\top\right)\right\}^{-1}\mathbb{P}_n\left(\widehat{U}  Y\right).
\end{equation}
We impose the following additional moment and nondegeneracy condition:
\newline
(A4) (Moment and nondegeneracy) The outcome satisfies $\mathbb{E}|Y|^2 < \infty$ and $\mathbb{E}(\|U\|^2 |Y|^2) < \infty$, and we also assume $\mathbb{E}\|U\|^4<\infty$. Further, the matrix $\Sigma_0$is full rank. 
\newline
Assumption (A4) combines standard finite-moment conditions with the requirements that $\Sigma_0$ be full rank, the usual full-rank design condition ensuring identifiability of $\theta_0$.
For later use, we rewrite the eigenfunction expansion in Theorem \ref{THM:asymp} in operator form. 
For each $j=1,\ldots,m$, define a linear map $L_{1j}:\mathcal{L}(\mathcal{H})\to\mathcal{H}$ by $L_{1j}(K)=\sum_{j' \neq j} (\lambda_j - \lambda_{j'})^{-1}\phi_{j'}\phi_{j'}^{\top} K\phi_j$, $K\in\mathcal{L}(\mathcal{H})$, where $\mathcal{L}(\mathcal{H})$ denotes the space of bounded linear operators on $\mathcal{H}$. In particular, with $K=(Z-\mu_Z)(Z-\mu_Z)^{\top}$, Theorem \ref{THM:asymp} gives 
$\sqrt{n}(\widehat{\phi}_j - \phi_j)=\mathbb{G}_n L_{1j}((Z-\mu_Z)(Z-\mu_Z)^{\top})+o_P(1)$.
Equivalently, $L_{1j}$ is the Fr\'{e}chet derivative at $V_0$ of the map that sends the covariance operator to its $j$-th eigenfunction.
The next theorem gives an asymptotic linear representation of $\widehat{\theta}_n$, making explicit the contribution of estimating eigenfunctions and forming the basis for both jackknife and bootstrap inference.
\begin{theorem}\label{THM:LM_Consistency} (Asymptotic linearity of PCR estimator)
Suppose Assumptions (A1)--(A4) hold and $\mathbb{E}\|Z\|^4<\infty$. Then
\[
	\sqrt{n}(\widehat{\theta}_n - \theta_0)
	=\mathbb{G}_n\left(\Sigma_0^{-1}
		\left[U\left(Y -\theta_0^\top U\right) + L_0\left\{(Z-\mu_Z)(Z -\mu_Z)^\top\right\}\right]\right) +o_P(1),
\]
where $L_0\{(Z-\mu_Z)(Z -\mu_Z)^{\top}\}$ is a mean-zero vector-valued functional of $Z$ that collects the additional variability arising from estimating the principal components. 
For $K=(Z-\mu_Z)(Z-\mu_Z)^{\top}$, it can be written as $L_0(K)=(Q_1(K), Q_2(K)^{\top}, Q_3(K)^{\top})^{\top}$, with the scalar term $Q_1(K)=-\sum_{j=1}^m\gamma_{0j}\langle\mu_Z,L_{1j}(K)\rangle \in\mathbb{R}$;
the $d$-vector $Q_2(K)=(Q_{21}(K),\ldots,Q_{2d}(K))^\top$, where
$Q_{2\ell}(K)=-\sum_{j=1}^m\gamma_{0j}\langle \mathbb{E}(X_{\ell}Z),L_{1j}(K)\rangle\in\mathbb{R}$, $\ell=1,\dots,d$;
the $m$-vector $Q_3(K)=\{Q_{31}(K),\ldots,Q_{3p}(K)\}^\top$, where for $\ell=1,\ldots,p$,
\[
	Q_{3\ell}(K)
	=\left\langle \mathbb{E}\left\{Z\left(Y -\theta_0^{\top}U \right)\right\}, L_{1\ell}(K)\right\rangle
		-\sum_{j=1}^m\gamma_{0j}\left\langle \mathbb{E}({Z_{\ast}}_{\ell}Z), L_{1j}(K) \right\rangle \in \mathbb{R}.
\]
The influence function on the right-hand side is square integrable under the stated moment conditions, and hence $\sqrt{n}(\widehat{\theta}_n - \theta_0)$ is asymptotically normal.
\end{theorem}

\begin{remark} (HS-PCR influence function)
Theorem \ref{THM:LM_Consistency} decomposes the asymptotic influence function of the HS-PCR estimator into two interpretable components: the standard regression influence function $\Sigma_0^{-1} U(Y - \theta_0^\top U)$, which would obtain if the true principal components were known, and a correction term $\Sigma_0^{-1} L_0((Z-\mu_Z)(Z-\mu_Z)^\top)$ that captures the additional variability from estimating eigenfunctions via AS-PCA. The correction term vanishes when the eigenfunction estimates are exact (oracle case), confirming that HS-PCR reduces to standard PCR when the principal components are known. The explicit form of the correction, derived in the proof, enables the bootstrap procedures of Section \ref{SUBSEC:Bootstrap} to account for eigenfunction estimation uncertainty without requiring separate resampling of the PCA step.
\end{remark}

\begin{remark} (Block jackknife) 
By Lemma~19.8 of \cite{Kosorok:08}, the asymptotic linearity in Theorem \ref{THM:LM_Consistency} implies that a block jackknife can be used for inference on $\theta_0$. Briefly, choose an integer $r > p+1 = d+m+1$. Let $k_{r, n}$ be the largest integer with $rk_{r,n} \leq n$, and set $N_{r, n} \equiv rk_{r,n}$. After computing the $\check{\theta}_n$ on the full sample, we repeatedly omit blocks of observations with indices $\ell, r+\ell, 2r+\ell,\ldots,(k_{r,n}-1)r + \ell$ for $\ell=1,\ldots,r$, and recompute the estimator on the remaining $N_{r,n}-k_{r,n}$ observations. The resulting $r$ jackknife replicates can be used to estimate the variance of $\widehat{\theta}_n$; see \cite{Kosorok:08} for details.
    The bootstrap analysis in the next subsection provides an alternative inference procedure and verifies that both the nonparametric and wild bootstraps are valid.
\end{remark}

\begin{remark} 
    Suppose the usual linear model $Y = \check{\alpha}_0 + \check{\beta}_0^\top X + \check{\gamma}_0^\top Z_{\ast} + \varepsilon$ holds with $\mathbb{E}[\varepsilon | X, Z_{\ast}] = 0$, $\mathbb{E}\varepsilon^4 < \infty$, $\mathbb{E}\|X\|^4 < \infty$, and $\mathbb{E}\|Z\|^4 < \infty$. Then $\theta_0 = (\check{\alpha}_0, \check{\beta}_0, \check{\gamma}_0)$ and $\mathbb{E}Y^4 < \infty$.
    Thus, $\theta_0$ in our principal component regression coincides with the usual linear-model parameter when the slope $\gamma$ lies in the principal-component subspace.
\end{remark}

\subsection{Bootstrap Inference for HS-PCR}
\label{SUBSEC:Bootstrap}

We now establish that both the nonparametric and wild bootstrap are valid for inference about all HS-PCR parameters $\theta_0 = (\alpha_0, \beta_0^\top, \gamma_0^\top)^\top$, including the functional regression coefficients $\gamma_0$ whose estimation depends on the AS-PCA eigenfunctions.

Let $\widetilde{\mathbb{P}}_n f \equiv n^{-1} \sum_{i=1}^n W_{ni} f(Y_i, X_i, Z_i)$, where $W_n = (W_{n1}, \cdots, W_{nn})$ are random weights independent of data, and define the corresponding empirical process $\widetilde{\mathbb{G}}_n = \sqrt{n}(\widetilde{\mathbb{P}}_n - \mathbb{P}_n)$. 
We impose the following assumption on the weights:
\newline 
(A5) (Bootstrap weights) Each $W_{ni}$ is either drawn from a multinomial distribution with $n$ categories of probabilities $1/n$ (the nonparametric bootstrap); or of the form $W_{ni} = \xi_i/\bar{\xi}_n$, where $\xi_1, \cdots, \xi_n$ are i.i.d. positive random variables with mean and variance one, and $\bar{\xi}_n  = n^{-1}\sum_{i=1}^n \xi_i$, with $\|\xi\|_{2,1} = \int_{0}^\infty \sqrt{\Pr(|\xi| > x)}dx < \infty$ (the wild bootstrap).
\newline
Assumption (A5) specifies standard choices of bootstrap weights for multiplier empirical processes. This condition rules out excessively heavy-tailed weights and ensures that the bootstrap empirical process converges to the same Gaussian limit as the original process, which is crucial for the validity of Theorems \ref{THM:bootstrap_v}--\ref{THM:bootstrap}.

To describe conditional weak convergence, we use the following notation.
Let $\{\widetilde{X}_n\}$ be a sequence of bootstrapped processes, with weights denoted by $W$, and let $X$ be a tight process.
We denote the conditional convergence of bootstrap laws as $\widetilde{X}_n \xrsquigarrow{P}{W} X$ if, for all $h \in BL_1$, $\sup_{h \in BL_1}\lvert \mathbb{E}_Wh(\widetilde{X}_n) - \mathbb{E}h(X)\rvert \xrightarrow{P} 0$ and $\mathbb{E}_Wh(\widetilde{X}_n)^{\ast} - \mathbb{E}_Wh(\widetilde{X}_n)_{\ast} \xrightarrow{P} 0$, where $BL_1$ is the space of functions with Lipschitz norm bounded by one, $\mathbb{E}_W$ denotes expectation with respect to the weights $W$ conditional on the data, and $h(\widetilde{X}_n)^{\ast}$ and $h(\widetilde{X}_n)_{\ast}$ denote measurable majorants and minorants with respect to the joint data (including the weights $W$).
An analogous definition holds for $\widetilde{X}_n \xrsquigarrow{as*}{W} X$.

We apply this to the empirical processes indexed by the function class $\mathcal{F}_1 \equiv \{\langle a, Z\rangle: a\in \mathcal{B}\}$ defined in Section \ref{SEC:HilbertSp}. Since $\mathcal{F}_1 \cdot \mathcal{F}_1$ is Donsker by Theorem \ref{THM:product}, and $P^{\ast}[\sup_{f \in \mathcal{F}_1 \cdot \mathcal{F}_1} \{f(X) - Pf\}^2] < \infty$, Theorem 2.7 of \cite{Kosorok:08} yields $\widetilde{\mathbb{G}}_n \xrsquigarrow{as*}{W} \mathbb{G}$ in $\ell^\infty(\mathcal{F})$, where $\mathbb{G}$ is the same tight mean-zero Gaussian process as in the non-bootstrap theory.

Define the bootstrapped covariance operator $\widetilde{V}_n^\ast = \mathrm{G}_N^\top \widetilde{\mathbb{P}}_n\{(Z - \widetilde{Z}_n^{\ast})(Z- \widetilde{Z}_n^{\ast})^\top\}\mathrm{G}_N$, with $\widetilde{Z}_n^{\ast} = \widetilde{\mathbb{P}}_n Z$, where $\mathrm{G}_N$ is the projection operator introduced in Section \ref{SUBSEC:Eigen}.
The next two results provide bootstrap analogues of the covariance-operator asymptotics.

\begin{theorem}\label{THM:bootstrap_v} (Bootstrap approximation of the covariance operator)
Under Assumptions (A1), (A2), and (A5),
$\sqrt{n}(\widetilde{V}_n^{\ast} - \widehat{V}_n) - \widetilde{\mathbb{G}}_n\{(Z-\mu)(Z-\mu)^\top\} \xrsquigarrow{P}{W} 0 $ in $\ell^\infty(\mathcal{B} \times \mathcal{B})$,
where $\widehat{V}_n$ is the empirical covariance operator defined in Section \ref{SUBSEC:Eigen} and $\mathcal{B}$ is the unit ball of $\mathcal{H}$.
\end{theorem}

\begin{theorem}\label{THM:bootstrap_v2} (Bootstrap functional delta method)
Assume (A1)--(A3) and (A5). Let $f: \ell^\infty (\mathcal{B}\times\mathcal{B}) \rightarrow \mathbb{D}$ be Hadamard differentiable tangentially to $\ell^{\infty}(\mathcal{B}\times\mathcal{B})$, for some Banach space $\mathbb{D}$, and $\dot{f}_{V_0}$ be the Hadamard derivative of $f$ at the population covariance operator $V_0$.
Then, for $1 \leq m \leq \infty$,  
\[
	\sqrt{n}\{f(\widetilde{V}_n^{\ast}) - f(\widehat{V}_n)\}-\dot{f}_{V_0}\mathbb{G}_n\{(Z-\mu)(Z-\mu)^\top\} \xrsquigarrow{P}{W} 0.
\]
\end{theorem}

We now turn to the regression estimator. Let $(\widetilde{X}, \widetilde{Z}_\ast)$ denote the bootstrapped version of $(X, Z_\ast)$, and define the bootstrap analogue of $\widehat{\theta}_n$ by
\[
	\widetilde{\theta}_n 
	= \left[\widetilde{\mathbb{P}}_n\left\{
	\begin{pmatrix}
		1 & \widetilde{X}^\top & \widetilde{Z}_\ast^\top
	\end{pmatrix}^\top\right\}^{\otimes 2}\right]^{-1}\widetilde{\mathbb{P}}_n\left\{
	\begin{pmatrix}
	1 & \widetilde{X}^\top & \widetilde{Z}_\ast^\top
	\end{pmatrix}^\top Y\right\},
\] 
where $\otimes 2$ denotes the outer product. 
Recall that $\widehat{\theta}_n = 
	\mathbb{P}_n(\widehat{U}\widehat{U}^\top)^{-1}
	\mathbb{P}_n(\widehat{U}Y) 
$. 
The following theorem shows that the bootstrap estimator has the same limiting distribution, conditionally on the data, as $\widehat{\theta}_n$.
\begin{theorem}\label{THM:bootstrap} (Bootstrap validity for principal component regression)
    Suppose Assumptions (A1)--(A5) hold. Then
\[
    \sqrt{n}(\widetilde{\theta}_n-\widehat{\theta}_n) -\widetilde{\mathbb{G}}_n\left[\Sigma_0^{-1}\left\{U(Y-\theta_0^\top U)
        + L_0\left((Z-\mu)(Z-\mu)^\top\right)\right\}\right]
        \xrsquigarrow{P}{W} 0.
\]
        Since the influence function on the right side matches the influence function given in Theorem~3.1, any bootstrap satisfying (A5) is asymptotically valid for inference on $\hat{\theta}_n$.
\end{theorem}
Theorem \ref{THM:bootstrap} confirms that any bootstrap satisfying (A5) correctly captures the full sampling variability of the HS-PCR estimator, including the variability introduced by the AS-PCA eigenfunction estimation step. In practice, we recommend the wild bootstrap with Rademacher weights for its computational simplicity and robustness to heteroscedasticity.
\section{Implementing AS-PCA and HS-PCR for Neuroimaging Data}
\label{SEC:NeuroImple}

Implementing AS-PCA requires specification of the projection basis $\Psi_N^{\ast}$ in Stage 1. For neuroimaging data, this choice requires careful consideration of the high dimensionality, complex spatial structure, and irregular anatomical boundaries that characterize brain imaging. We describe our recommended implementation of AS-PCA for this setting, which then feeds directly into HS-PCR for regression analysis.

Selecting an efficient initial basis $\Psi_N^{\ast}$ is crucial for the entire estimation procedure. 
Various nonparametric methods have been developed for neuroimage analysis; see, for example, tensor-product kernel smoothing \citep{Zhu:Fan:Kong:14}, tensor product B-spline smoothing \citep{Shi:etal:22}, and multivariate spline over triangulation \citep[MST;][]{Lai:Schumaker:07}.
Among these approaches, MST has demonstrated strong performance for multi-dimensional imaging data, as evidenced by applications to bivariate penalized spline analysis of 2D images \citep{Lai:Wang:13} and trivariate penalized spline analysis of 3D images \citep{Li:etal:24}.

While univariate splines are well understood and widely applied, their multivariate counterparts introduce additional theoretical and computational challenges. 
One effective approach to constructing multivariate splines is through triangulations, which partition a domain into non-overlapping simplices (triangles in two dimensions, tetrahedra in three dimensions).
The construction of MST involves two main steps: first, building a triangulation to approximate the domain of interest; and second, constructing multivariate splines on this triangulation. 
A multivariate spline over a triangulation is a piecewise polynomial function on the triangulated domain that ensures smooth connections across simplex boundaries, typically by requiring continuity of derivatives up to a specified order. 
These splines are especially valuable for statistical applications on irregularly shaped domains, such as neuroimaging data. By using triangulations, multivariate splines provide flexible function approximations while maintaining computational efficiency. In contrast to tensor-product kernel or spline methods, which suffer from grid alignment constraints, triangulation-based splines adapt naturally to complex data distributions and heterogeneous spatial structures.
Theoretical properties and statistical efficiency of MST have been established in \cite{Lai:Wang:13} and \cite{Li:etal:24}, which ensures that Assumption (A2) is satisfied for our choice of $\Psi_N^{\ast}$.
For further technical details on the construction and properties of MST, we refer to  \cite{Lai:Schumaker:07}.

\section{Simulation study}
\label{SEC:Simulation}

We assess the finite-sample performance of AS-PCA and HS-PCR through simulation studies designed to evaluate: 
(i) eigenvalue and eigenfunction recovery by AS-PCA (Tables \ref{TAB:Eigens2D}, \ref{TAB:Eg3_MSE}); 
(ii) regression coefficient estimation accuracy of HS-PCR (Tables \ref{TAB:Eg2_MSE_beta}, \ref{TAB:Eg2_MSE_gamma}, \ref{TAB:Eg4_MSE}); and 
(iii) bootstrap coverage validity for HS-PCR inference (Tables \ref{TAB:Eg2_MSE_beta}, \ref{TAB:Eg2_MSE_gamma}, \ref{TAB:Eg4_MSE}).
The data are generated from:
\[
    Y_i = \alpha + \beta^{\top} X_i + \langle \gamma, Z_i\rangle + \varepsilon_i,
    \qquad i=1,\ldots,n.
\]
The Euclidean covariates $X_i\in\mathbb{R}^d$ are drawn independently from $\text{MVN}(0_d, \Omega_d(r))$, where the covariance matrix $\Omega_d(r)$ has entries $\{\Omega_d(r)\}_{\ell \ell^{\prime}}=r^{|\ell-\ell^{\prime}|}$,
which allows us to control the dependence among covariates via the parameter $r$. The random errors $\varepsilon_i$ are i.i.d. $N(0,1)$. For the Hilbert-valued component $Z_i$, we consider several different settings described below.

Throughout, we set $d = 4$, $\alpha_0=1$, and $\beta_0=(1,1,1,1)^{\top}$. We examine two correlation scenarios: $r = 0$ (independent covariates) and $r = 0.5$ (moderately correlated covariates), and three sample sizes: $n = 100, 500$, and $2000$. We conduct $100$ Monte Carlo replications for each setting, compare the estimated basis and eigenvalues to their underlying truth, measure the estimation accuracy of the Euclidean covariates by mean squared errors (MSEs), and quantify the inference performance of the bootstrap confidence interval.

\subsection{Two-dimensional imaging setting}
\label{SUBSEC:2D_ADNI}

\begin{figure}[htbp]
\begin{center}
\begin{tabular}{ccccccc}
	\includegraphics[height=.52in]{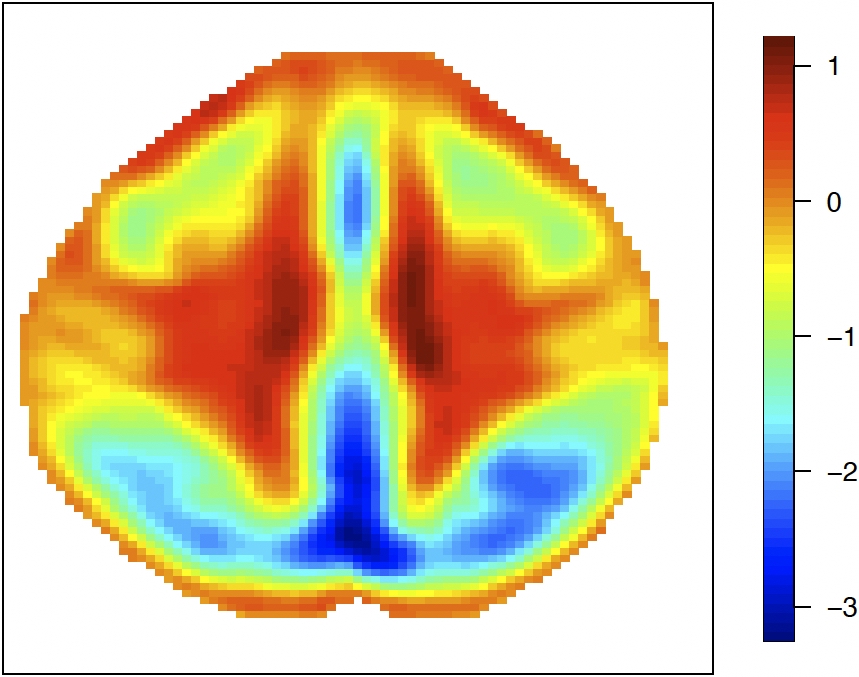} 
		& \includegraphics[height=.52in]{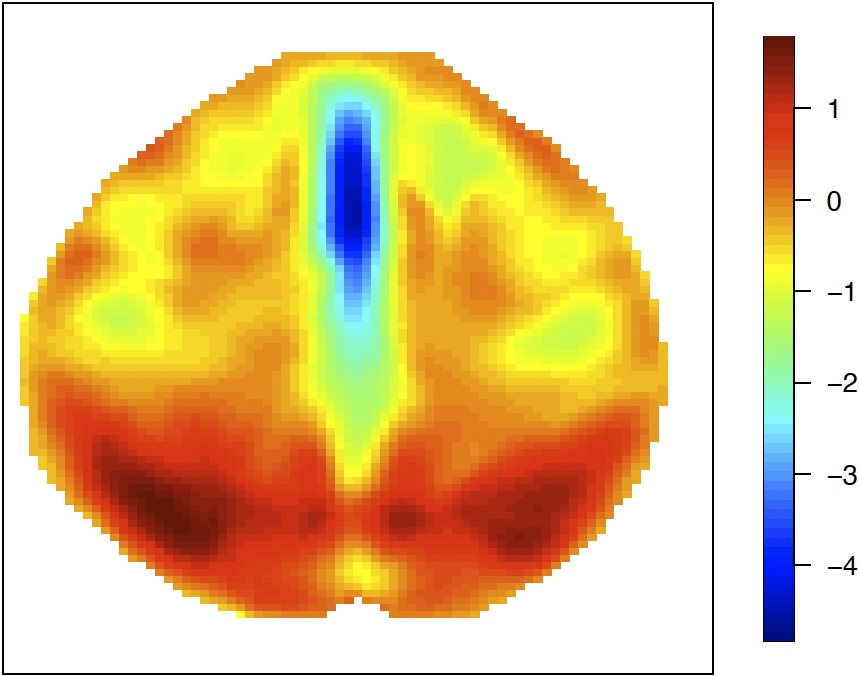} 
		& \includegraphics[height=.52in]{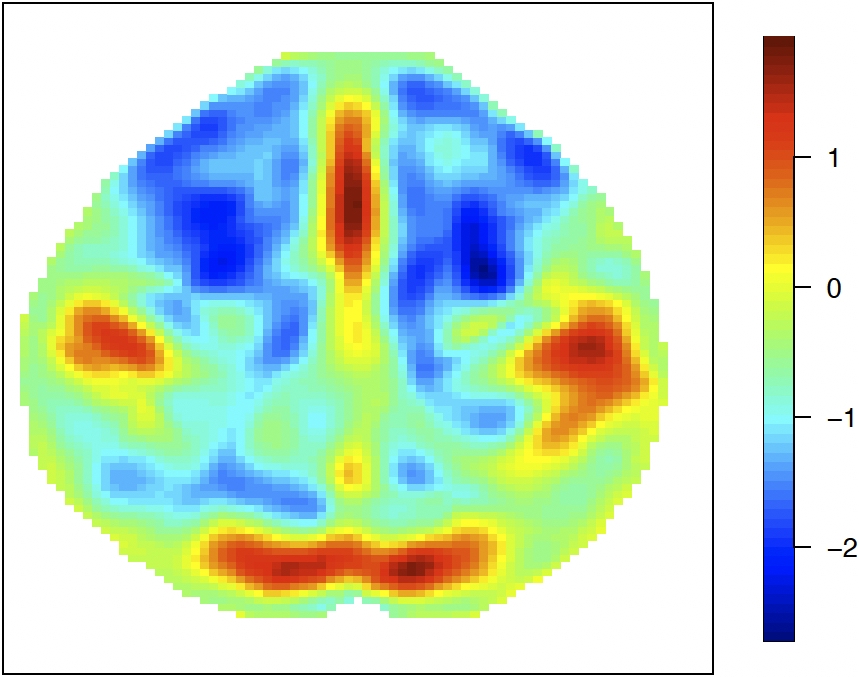}
		& \includegraphics[height=.52in]{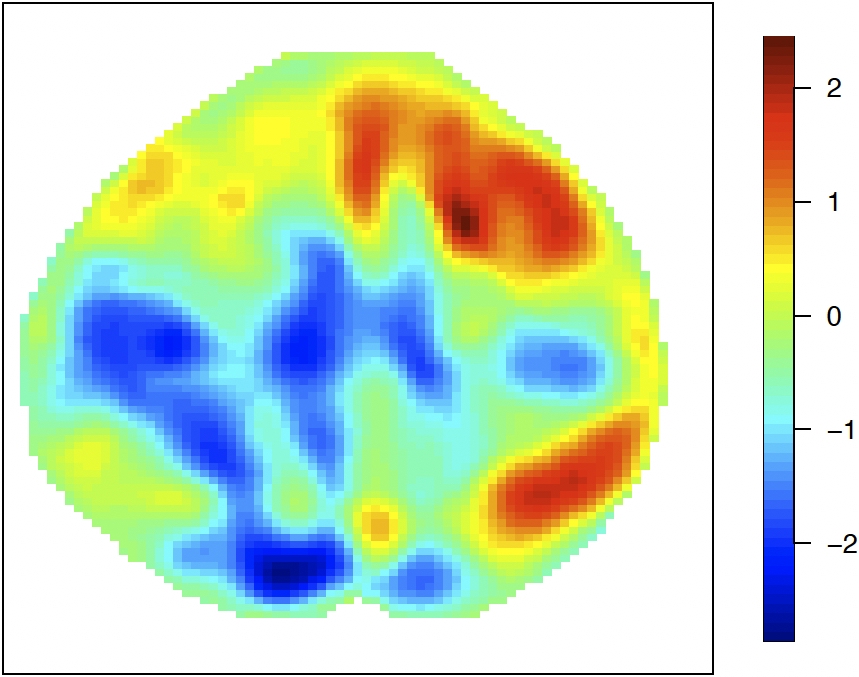} 
		& \includegraphics[height=.52in]{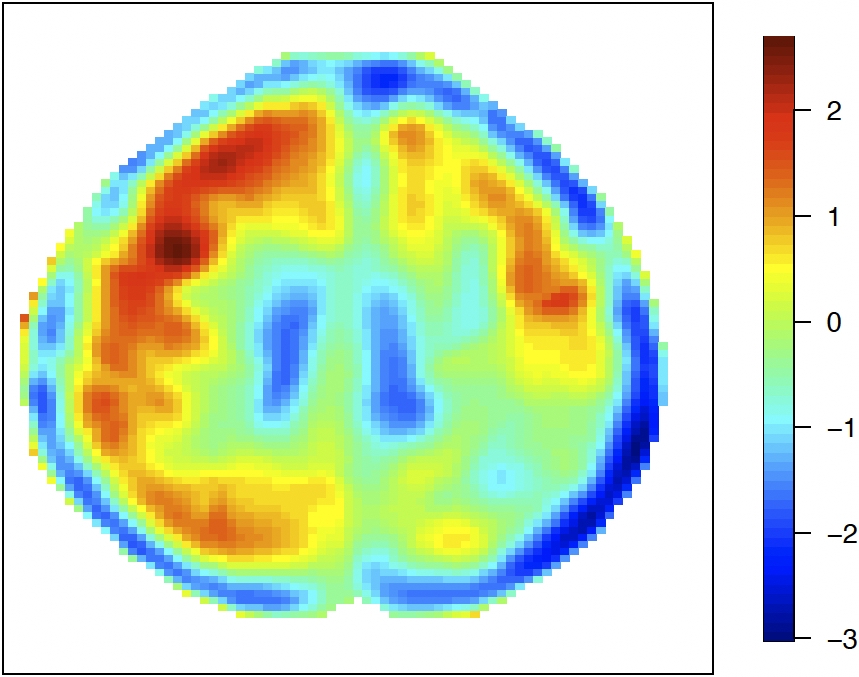} 
		& \includegraphics[height=.52in]{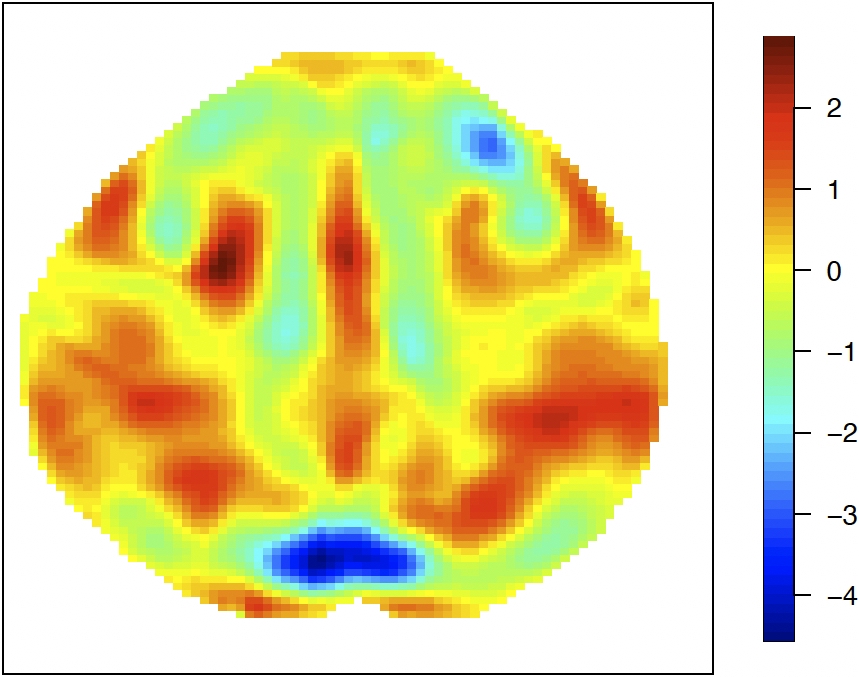}
		& \includegraphics[height=.52in]{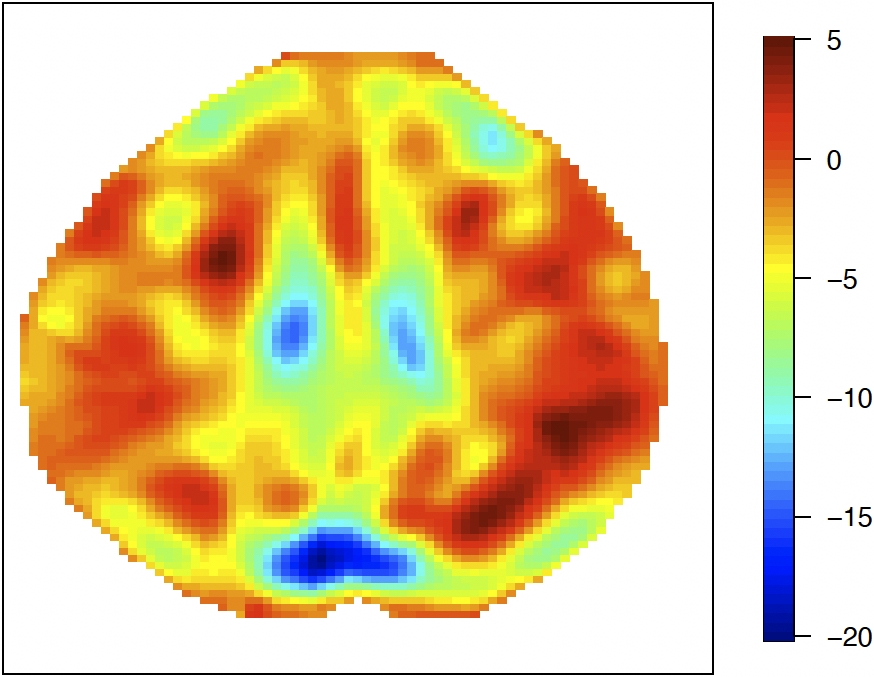}\vspace{-2pt} \\ 
	$\phi_1~~$ & $\phi_2$ ~~ & $\phi_3$~~ & $\phi_4~~$ & $\phi_5$ ~~ & $\phi_6$ & $\gamma_0$~~ \\
	\includegraphics[height=.52in]{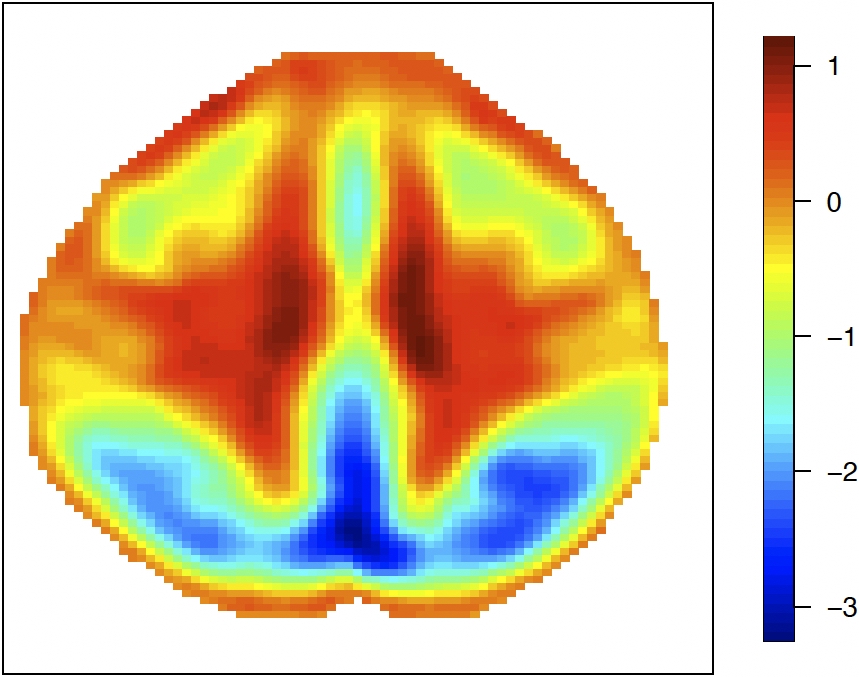} 
		& \includegraphics[height=.52in]{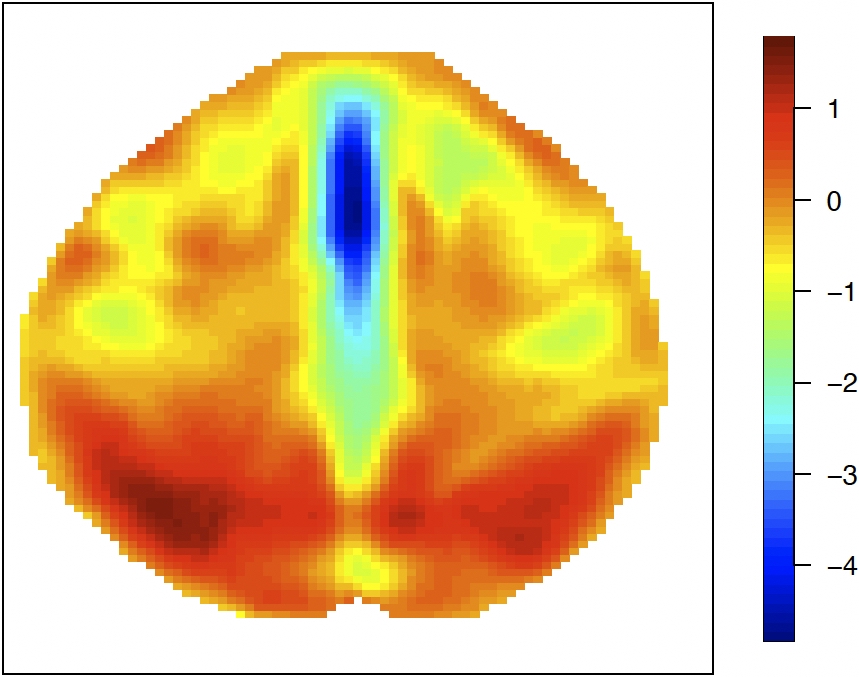} 
		& \includegraphics[height=.52in]{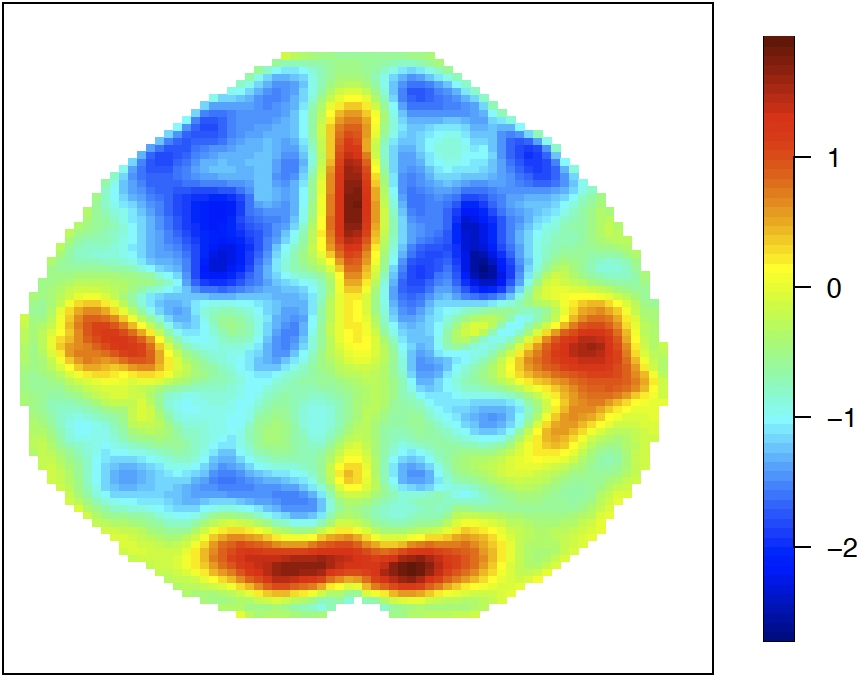}
		& \includegraphics[height=.52in]{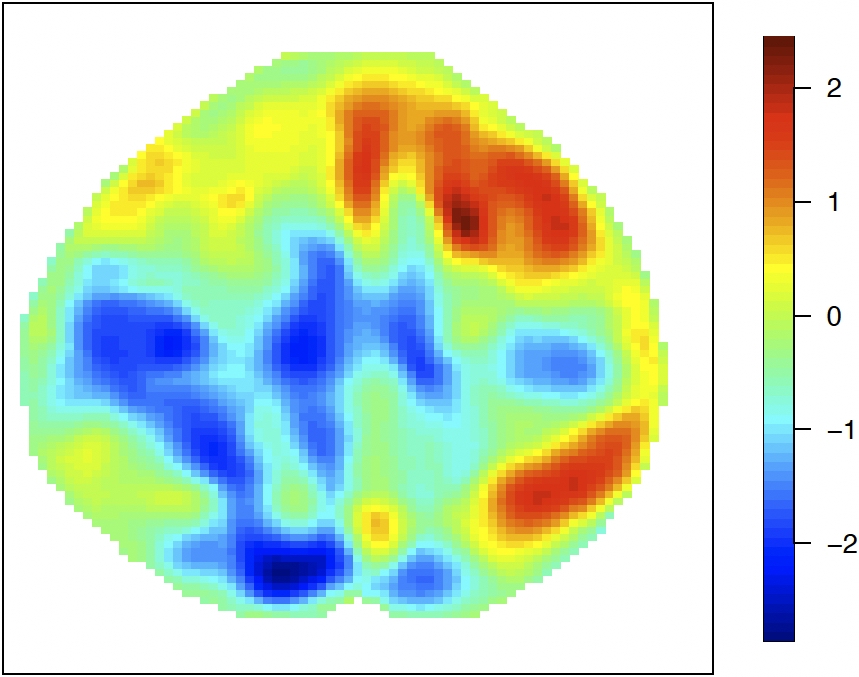} 
		& \includegraphics[height=.52in]{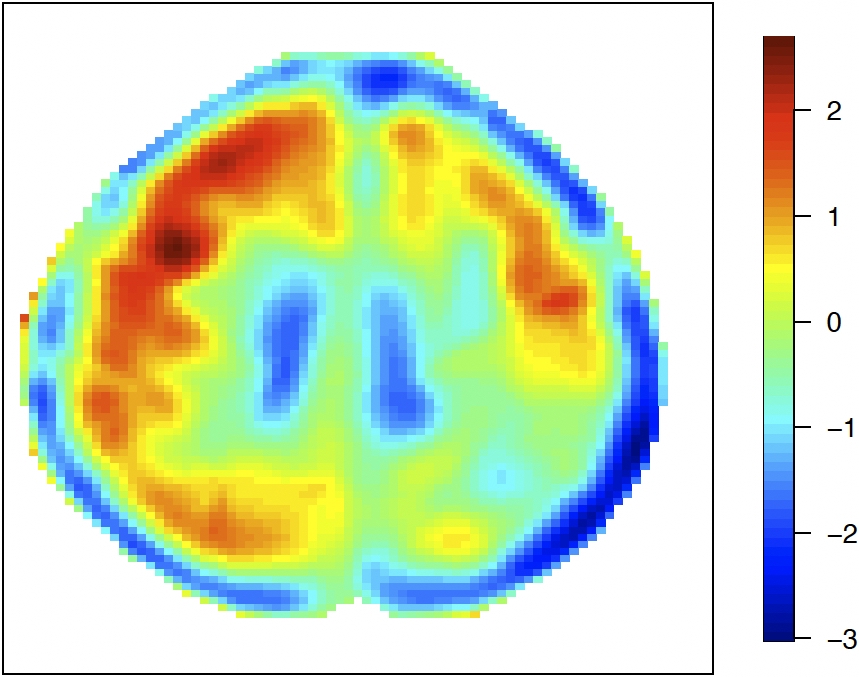} 
		& \includegraphics[height=.52in]{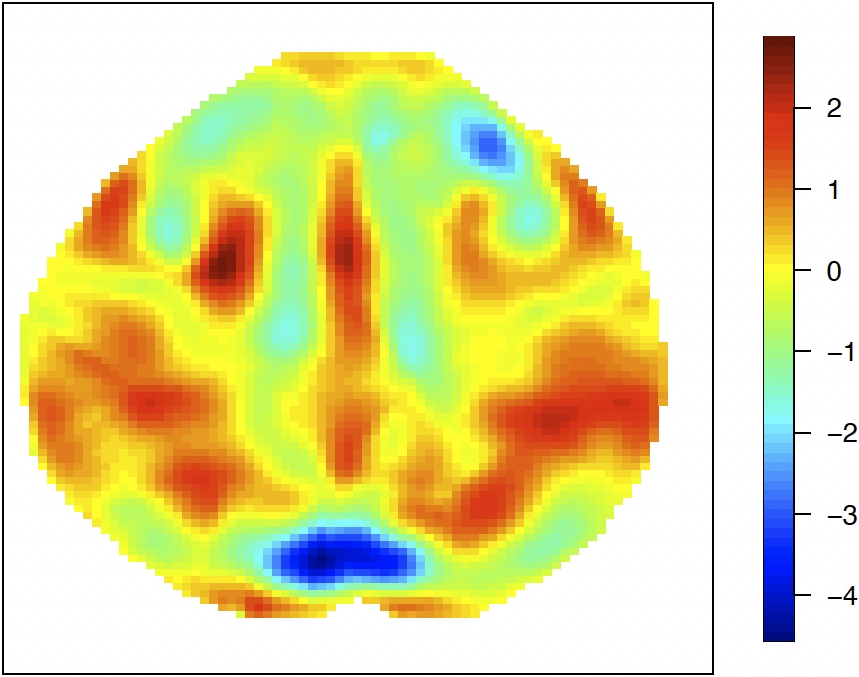}
		& \includegraphics[height=.52in]{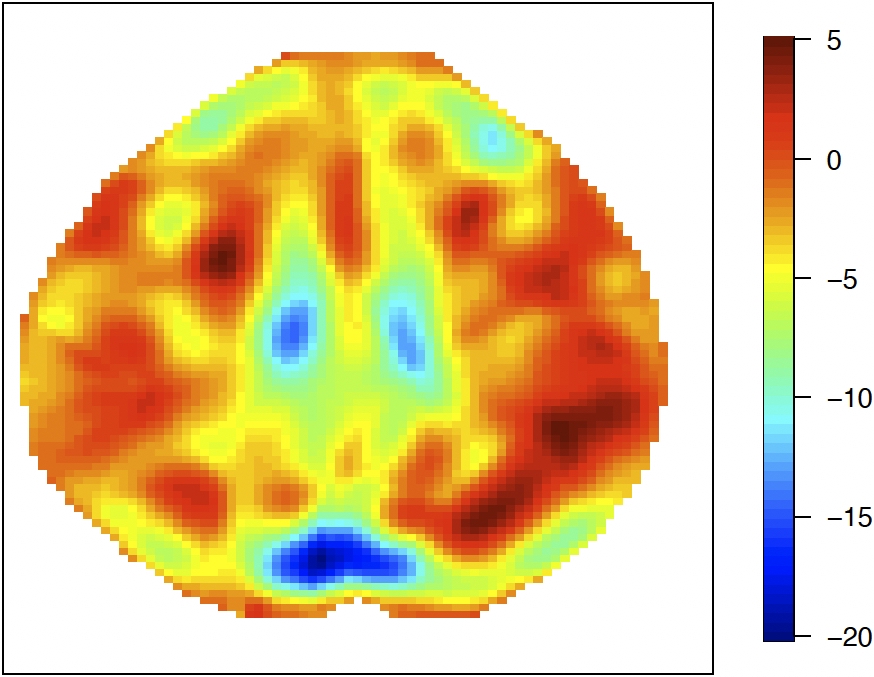} \vspace{-2pt} \\ 
	$\widehat{\phi}_1~~$ & $\widehat{\phi}_2$ ~~ & $\widehat{\phi}_3$~~ & $\widehat{\phi}_4~~$ & $\widehat{\phi}_5$ ~~ & $\widehat{\phi}_6$ & $\widehat{\gamma}$~~ \\
\end{tabular}
\end{center} \vspace*{-.6cm}
\caption{Comparison of true basis functions $\{\phi_j\}_{j=1}^6$ and coefficient function $\gamma$ (top) and their estimates $\{\widehat{\phi}_j\}_{j=1}^6$ and $\widehat{\gamma}$ (bottom) for a randomly selected iteration with $n=2000$ in the 2D imaging setting.}
\label{FIG:eg1-1}
\end{figure}

For the Hilbert-valued component, we first consider a 2D imaging setting. To mimic the full complexity of the real data and illustrate the performance of the proposed method in a challenging scenario with complicated patterns, we use the leading empirical basis functions obtained from the real data analysis in Section \ref{SEC:DA}. To be specific, we construct the imaging data as $Z_i(s_1,s_2) = \sum_{j=1}^{6} \lambda_j^{1/2}U_{ij}\phi_j(s_1,s_2)$ on a $79 \times 95$ pixel grid, where the eigenvalues $(\lambda_1,\lambda_2,\lambda_3,\lambda_4,\lambda_5,\lambda_6)=(3.5,3,2.5,2,1.5,1)$ and  $(U_{i1},\ldots,U_{i6})^\top$ follows $\text{MVN}(0_6, I_6)$. The basis functions $\{\phi_j\}_{j=1}^6$ are the leading six principal component bases that are selected by PVE from the real data analysis; see Section \ref{SUBSEC:LM} for more details. Figure \ref{FIG:eg1-1} (top row) and Figure \ref{FIG:beta2hat} (left six panels) display the plots of these basis functions. We set $\gamma_0=1.5\phi_1+\phi_2+2\phi_3+2.5\phi_4+1.5\phi_5+3\phi_6$.

The estimates of the basis functions are shown at the bottom left six panels of Figure \ref{FIG:eg1-1}, which is selected randomly from a single iteration. 
The figure shows that with $n = 2000$ observations, the estimated basis functions (bottom row) closely reproduce both the pattern and magnitude of the true functions (top row), providing empirical validation of Theorem \ref{THM:maxes}.
In addition, Table \ref{TAB:Eigens2D} presents the estimation accuracy of the six eigenvalues ($\lambda_1$ through $\lambda_6$) and the estimated number of eigenvalues under the PVE ($\widehat{m}_n$) criterion.
We observe that MSEs for all eigenvalues decrease as sample size increases across both correlation scenarios, demonstrating improved estimation precision with larger samples.
The results in Table \ref{TAB:Eigens2D} further support the conclusions in Theorem \ref{THM:maxes} and Corollary \ref{COR:PVE_Consistency}.

\begin{table}
\caption{Mean Squared Errors (MSEs$\times 10^{-2}$) for AS-PCA eigenvalue estimates $\lambda_j$ $(j=1,\ldots,6)$, and estimated number of principal components $\widehat{m}_n$ in the 2D imaging setting.
Results are based on $100$ Monte Carlo simulations with varying sample sizes ($n = 100, 500, 2000$) and correlation scenarios ($r = 0, 0.5$). \label{TAB:Eigens2D}}
\setlength{\tabcolsep}{3pt}
\renewcommand{\arraystretch}{0.85}
\begin{center}
\resizebox{0.8\textwidth}{!}{
\begin{tabular}{rrrrrrrrrrrrrrr} \hline
	\multirow{2}[0]{*}{$n$} & \multicolumn{7}{c}{$r=0$} & \multicolumn{7}{c}{$r=0.5$} \\ 
	\cmidrule(lr){2-8} \cmidrule(lr){9-15}
	& $\lambda_1$ & $\lambda_2$ & $\lambda_3$ & $\lambda_4$ & $\lambda_5$ & $\lambda_6$ & $\widehat{m}_n$ & $\lambda_1$ & $\lambda_2$ & $\lambda_3$ & $\lambda_4$ & $\lambda_5$ & $\lambda_6$ & $\widehat{m}_n$ \\ \hline
	100 & 39.66 & 9.04 & 6.82 & 7.58 & 4.57 & 4.24 & 5.98& 39.66 & 9.04 & 6.82 & 7.58 & 4.57 & 4.24 & 5.98  \\
	500 & 4.77 & 3.67 & 2.41 & 1.42 & 0.83 & 0.45 & 6.00 & 4.77 & 3.67 & 2.41 & 1.42 & 0.83 & 0.45 & 6.00\\
	2000 & 1.25 & 0.97 & 0.64 & 0.45 & 0.18 & 0.09 & 6.00 & 1.25 & 0.97 & 0.64 & 0.45 & 0.18 & 0.09 & 6.00 \\
	\hline
\end{tabular}}
\end{center}\vspace{-8pt}
\end{table}

To assess the estimation accuracy and uncertainty quantification of the proposed methods, we report MSEs and empirical coverage rates of the $95\%$ bootstrap confidence interval for the coefficients $\alpha$, $\beta$, and $\gamma$. 
Tables \ref{TAB:Eg2_MSE_beta} and \ref{TAB:Eg2_MSE_gamma} summarize the results under varying sample sizes ($n = 100, 500, 2000$) and correlation levels ($r = 0, 0.5$).
Note that with $100$ Monte Carlo replications, the standard error for the $95\%$ confidence level is around $2.18\%$.
As expected, the MSEs decrease significantly with increased sample size, reflecting the consistency of the estimators. The presence of correlation ($r = 0.5$) generally leads to larger estimation errors of linear coefficients $\alpha$ and $\beta$ compared to the independent setting ($r = 0$), indicating increased complexity in estimation under dependence. 
PVE criteria yield coverage rates close to the nominal level across all settings.
These findings support the robustness and practical applicability of the proposed estimation method, and provide numerical validation of the theoretical results established in Theorem \ref{THM:LM_Consistency}.
Figure \ref{FIG:eg1-1} also illustrates the effectiveness of our approach in estimating the functional-coefficient $\gamma$. 
The estimates accurately capture the true patterns of the function. The visual similarity across both plots confirms that the estimate performs well in reconstructing the true signal.

\begin{table}
\caption{Mean Squared Errors (MSEs$\times 10^{-3}$) and empirical coverage rates (in parentheses) of the 95\% HS-PCR bootstrap confidence intervals for $\alpha$ and $\beta$ in the 2D imaging setting. Results are reported for different sample sizes ($n = 100, 500, 2000$) and correlation scenarios ($r = 0, 0.5$), based on $100$ Monte Carlo replications.\label{TAB:Eg2_MSE_beta}}
\setlength{\tabcolsep}{3pt}
\renewcommand{\arraystretch}{.85}
\begin{center}
\resizebox{0.7\textwidth}{!}{
\begin{tabular}{rrrrrrrrrrr}
	\hline
	\multirow{2}[0]{*}{$n$} &\multicolumn{5}{c}{$r=0$} & \multicolumn{5}{c}{$r=0.5$}\\
	\cmidrule(lr){2-6} \cmidrule(lr){7-11}
	& \multicolumn{1}{c}{$\alpha$} & \multicolumn{1}{c}{$\beta_1$} & \multicolumn{1}{c}{$\beta_2$} & \multicolumn{1}{c}{$\beta_3$} & \multicolumn{1}{c}{$\beta_4$}  & \multicolumn{1}{c}{$\alpha$} & \multicolumn{1}{c}{$\beta_1$} & \multicolumn{1}{c}{$\beta_2$} & \multicolumn{1}{c}{$\beta_3$} & \multicolumn{1}{c}{$\beta_4$}\\
	\hline
	100 & 15.48 & 14.14 & 10.69 & 10.30 & 14.69 & 15.48 & 14.96 & 22.90 & 22.77 & 14.25 \\
	& (98\%) & (97\%) & (100\%) & (97\%) & (96\%) & (98\%) & (98\%) & (97\%) & (95\%) & (97\%)  \\
	500 & 1.96 & 2.49 & 1.56 & 1.51 & 2.27  & 1.96 & 2.46 & 3.32 & 3.74 & 2.62  \\
	& (96\%) & (92\%) & (98\%) & (97\%) & (96\%)  & (96\%) & (96\%) & (94\%) & (93\%) & (96\%)  \\
	2000 & 0.54 & 0.51 & 0.45 & 0.58 & 0.42 & 0.54 & 0.68 & 0.95 & 0.69 & 0.66 \\
	& (94\%) & (97\%) & (98\%) & (94\%) & (97\%)  & (94\%) & (95\%) & (92\%) & (97\%) & (96\%) \\
	\hline
\end{tabular}}
\end{center}\vspace{-16pt}
\end{table}%

\begin{table}
\caption{Mean Squared Errors (MSEs$\times 10^{-1}$) and empirical coverage rates (in parentheses) of the 95\% HS-PCR bootstrap confidence interval for $\gamma$ in the 2D imaging setting.
Results are shown for various sample sizes ($n = 100, 500, 2000$) and correlation levels ($r = 0, 0.5$), based on $100$ Monte Carlo replications. \label{TAB:Eg2_MSE_gamma}}
\setlength{\tabcolsep}{3pt}
\renewcommand{\arraystretch}{0.85}
\resizebox{0.95\textwidth}{!}{
\begin{tabular}{rrrrrrrrrrrrr}
	\hline 
	\multirow{2}[0]{*}{$n$} & \multicolumn{6}{c}{$r=0$}     & \multicolumn{6}{c}{$r=0.05$} \\
	\cmidrule(lr){2-7} \cmidrule(lr){8-13}
	& $\gamma_1$ & $\gamma_2$ & $\gamma_3$ & $\gamma_4$ & $\gamma_5$ & $\gamma_6$
	& $\gamma_1$ & $\gamma_2$ & $\gamma_3$ & $\gamma_4$ & $\gamma_5$ & $\gamma_6$\\
	\hline
	100 & 6.60 & 7.73 & 10.41 & 13.02 & 8.25 & 5.79 & 6.60 & 7.73 & 10.41 & 13.02 & 8.25 & 5.79 \\
	& (93\%) & (100\%) & (100\%) & (94\%) & (100\%) & (93\%)
	& (93\%) & (100\%) & (100\%) & (94\%) & (100\%) & (93\%)\\
	500 & 3.05 & 4.51 & 3.87 & 3.00 & 2.93 & 0.73 
	& 3.05 & 4.51 & 3.87 & 3.00 & 2.93 & 0.73\\
	& (92\%) & (100\%) & (98\%) & (96\%) & (98\%) & (96\%)& (92\%) & (100\%) & (98\%) & (96\%) & (98\%) & (96\%)  \\
	2000 & 0.52 & 1.32 & 1.19 & 0.84 & 0.89 & 0.12 & 0.52 & 1.32 & 1.19 & 0.84 & 0.89 & 0.12 \\
	& (96\%) & (99\%) & (93\%) & (94\%) & (93\%) & (97\%) & (96\%) & (99\%) & (93\%) & (94\%) & (93\%) & (97\%) \\
	\hline
\end{tabular}}\vspace{-6pt}
\end{table}%

\subsection{Three-dimensional imaging setting}

\begin{figure}[htbp]
\begin{tabular}{ccccccc}
	& $\phi_1$ & $\phi_2$ & $\gamma_0$ & $\widehat{\phi}_1$ & $\widehat{\phi}_2$ & $\widehat{\gamma}$ \\
	\raisebox{.2in}{\shortstack{{\fontsize{9}{11}\selectfont Transverse} \\ {\fontsize{9}{11}\selectfont (Axial View)}}}
	& \includegraphics[height=.52in]{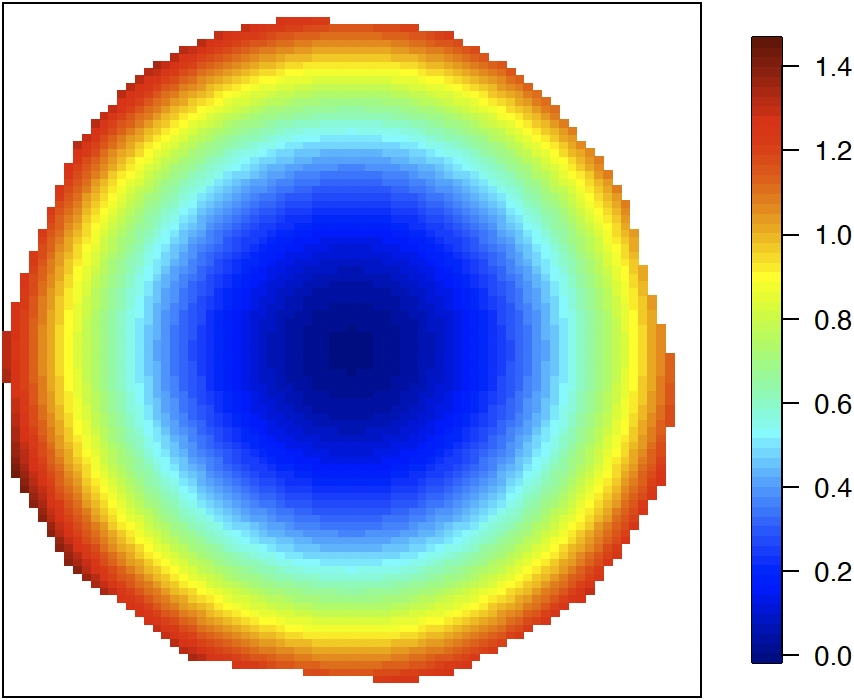} 
	& \includegraphics[height=.52in]{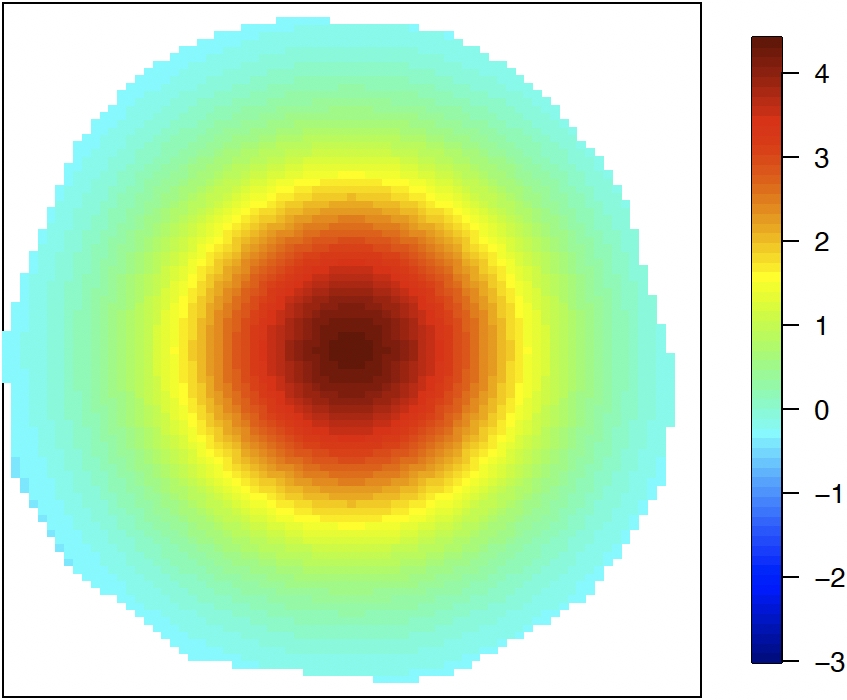} 
	& \includegraphics[height=.52in]{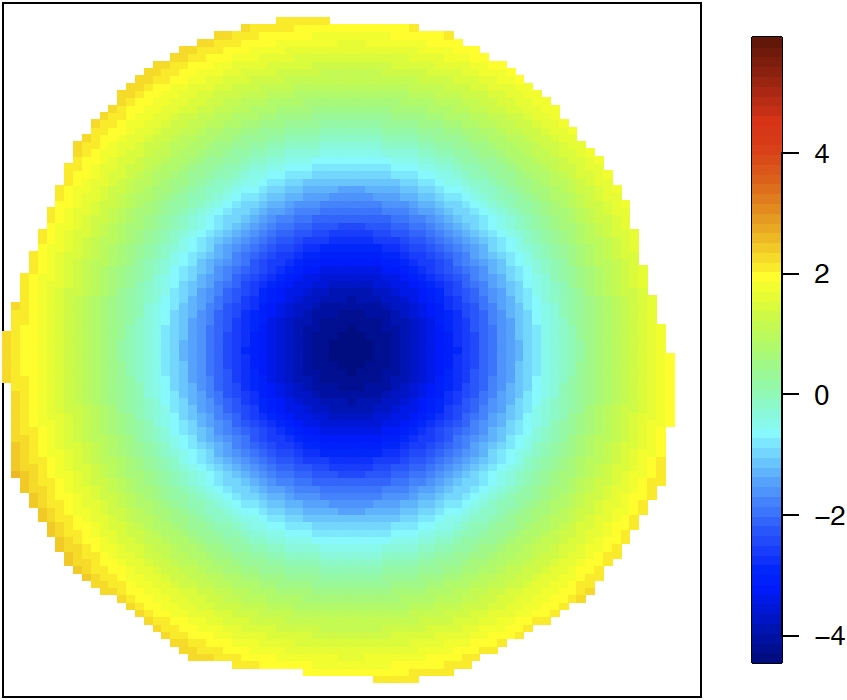}
	& \includegraphics[height=.52in]{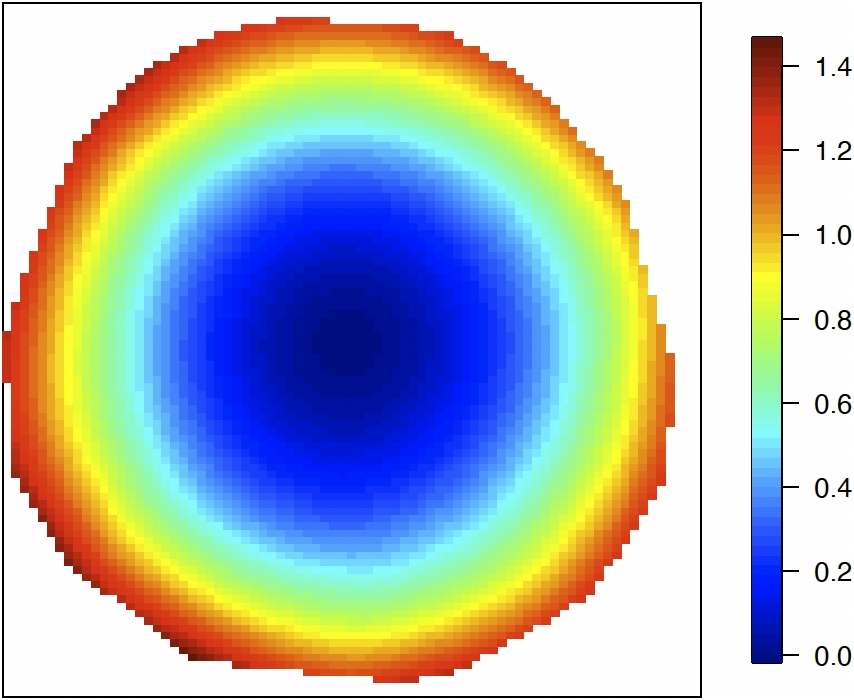} 
	& \includegraphics[height=.52in]{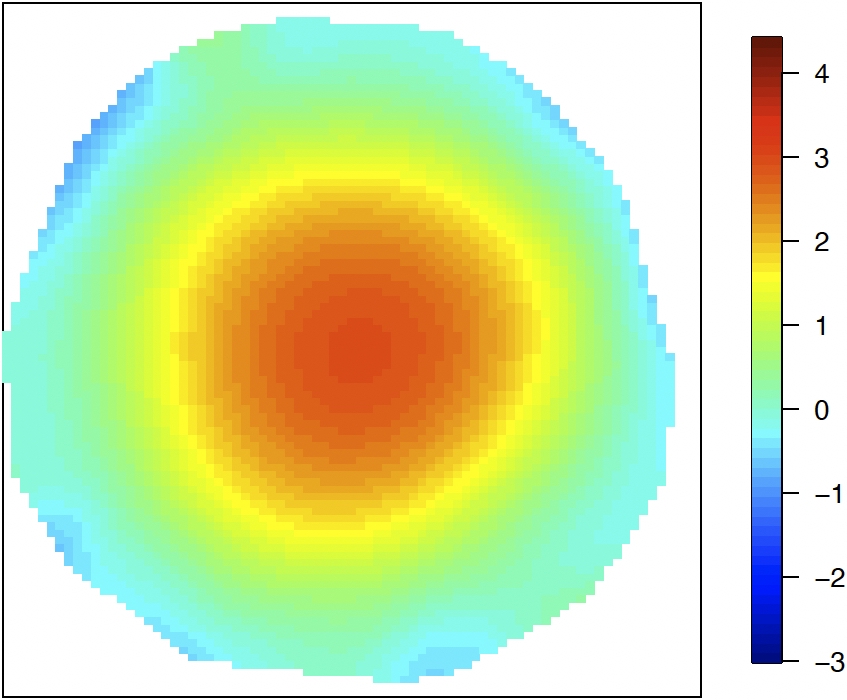} 
	& \includegraphics[height=.52in]{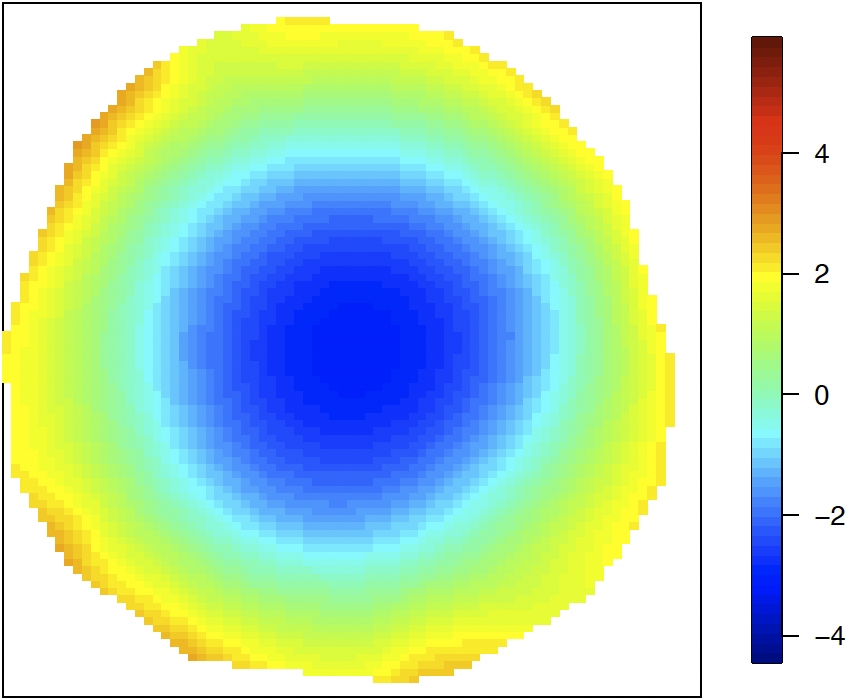}\\
	\raisebox{.2in}{\shortstack{{\fontsize{9}{11}\selectfont Coronal} \\ {\fontsize{9}{11}\selectfont (Frontal View)}}}
	& \includegraphics[height=.52in]{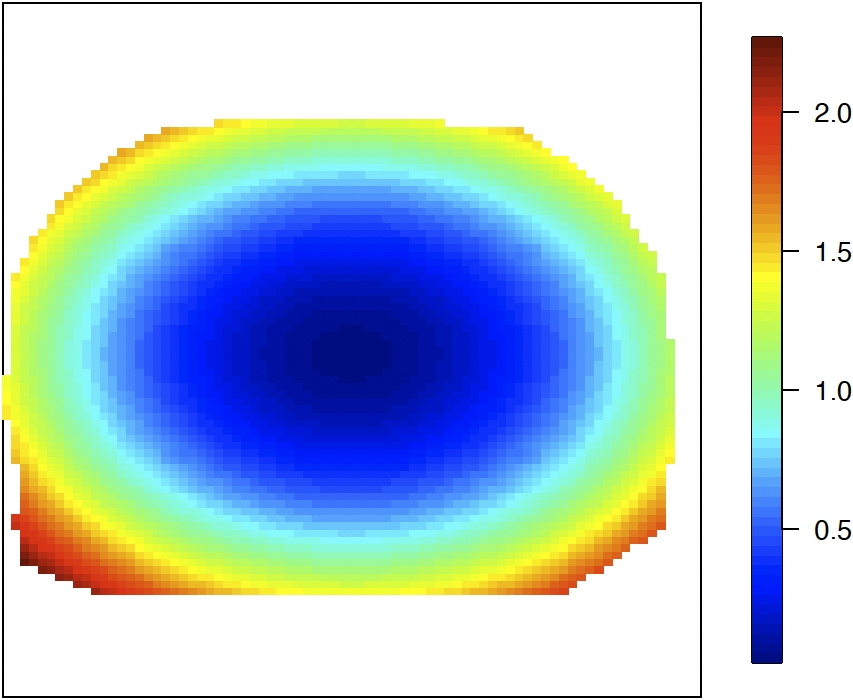}
	& \includegraphics[height=.52in]{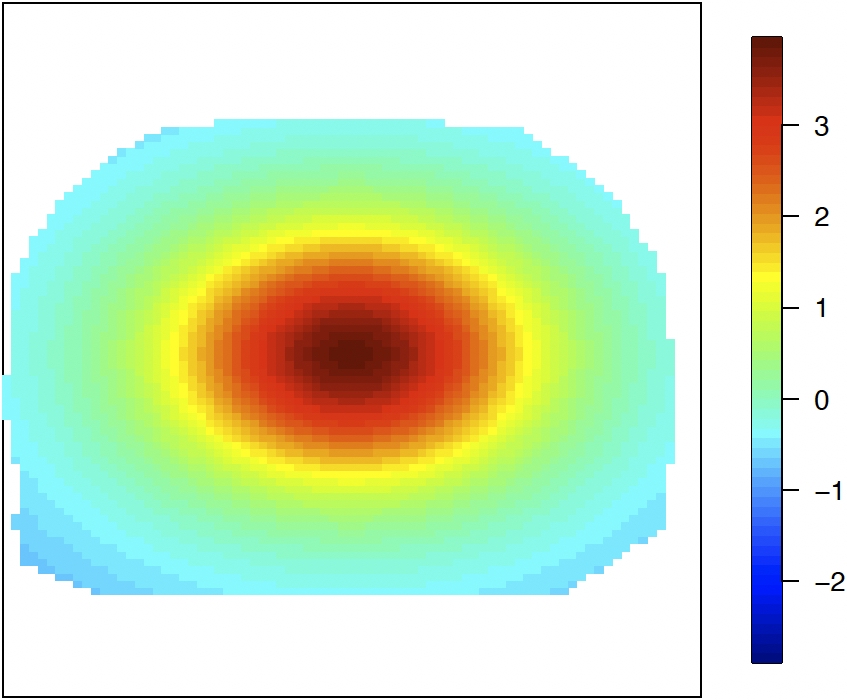} 
	& \includegraphics[height=.52in]{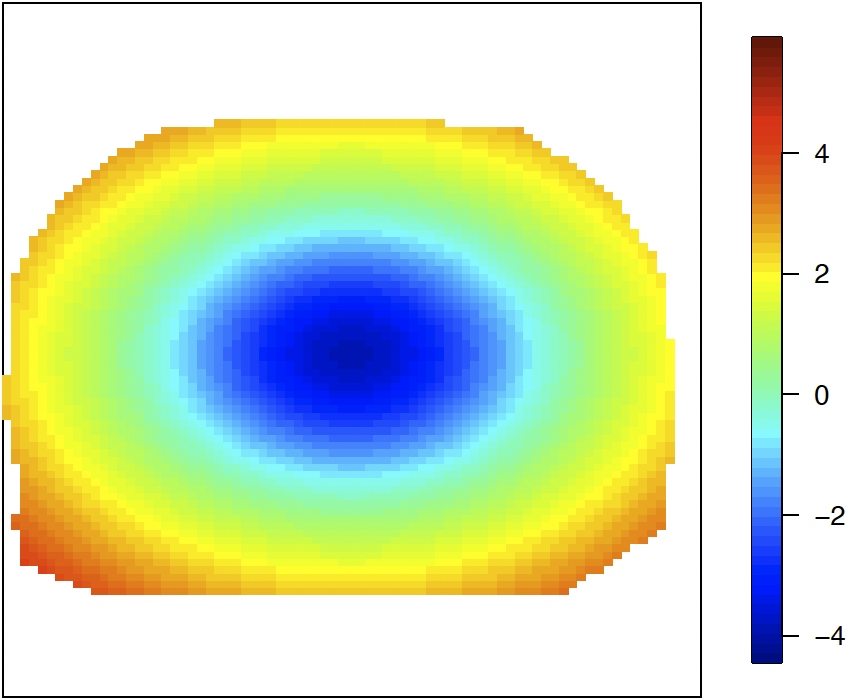}
	& \includegraphics[height=.52in]{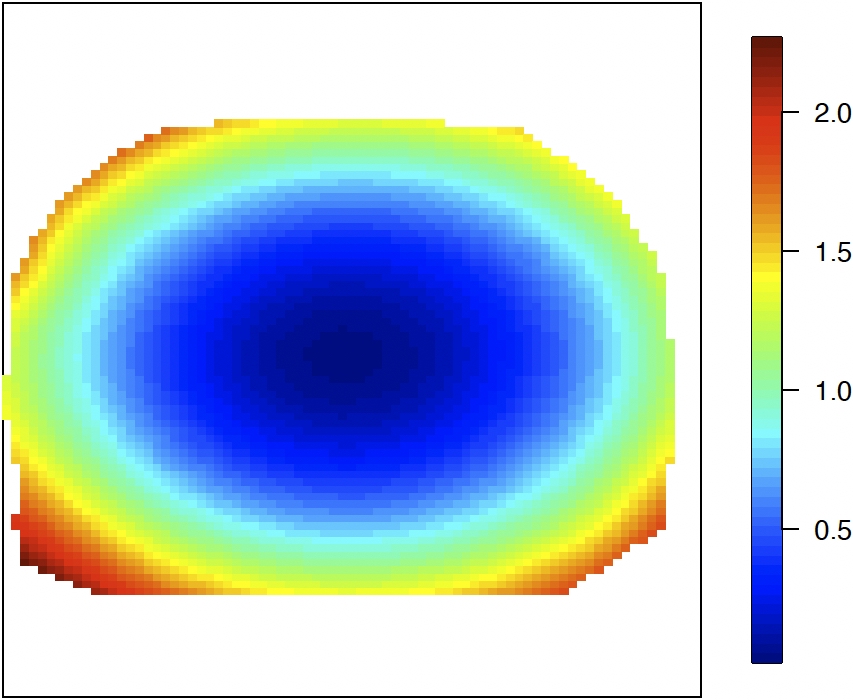} 
	& \includegraphics[height=.52in]{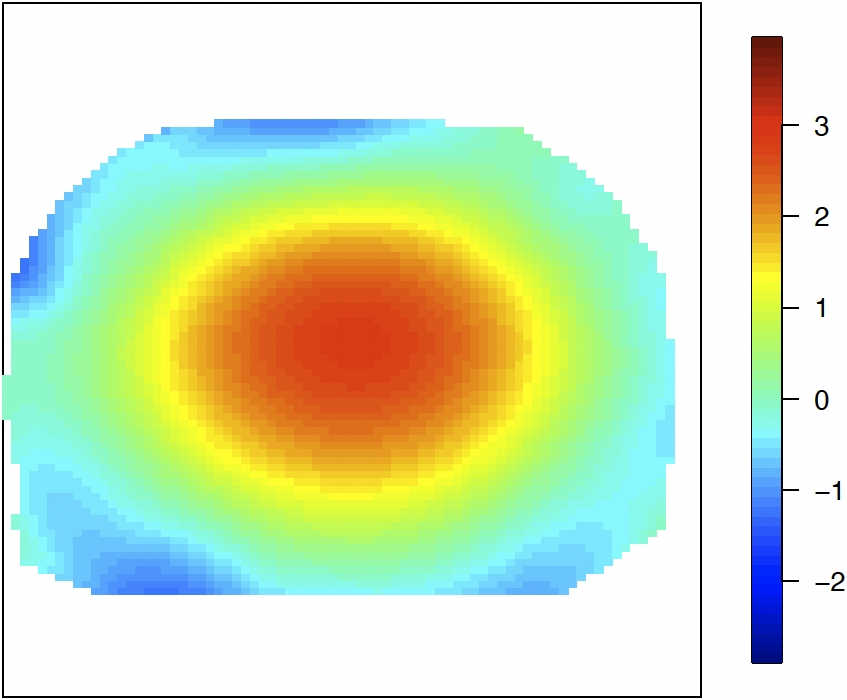} 
	& \includegraphics[height=.52in]{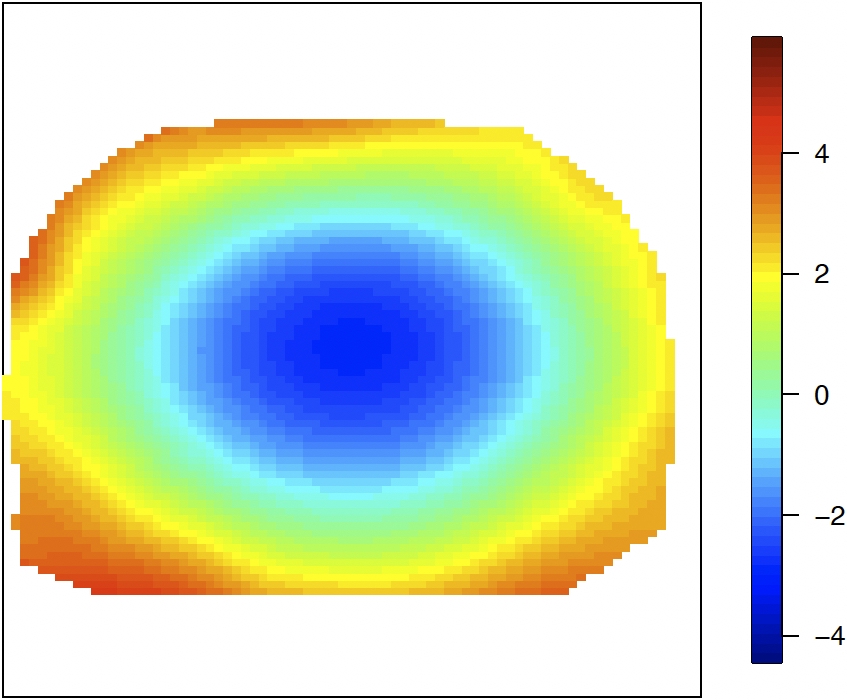} \\
	\raisebox{.2in}{\shortstack{{\fontsize{9}{11}\selectfont Sagittal} \\ {\fontsize{9}{11}\selectfont (Side View)}}}
	& \includegraphics[height=.52in]{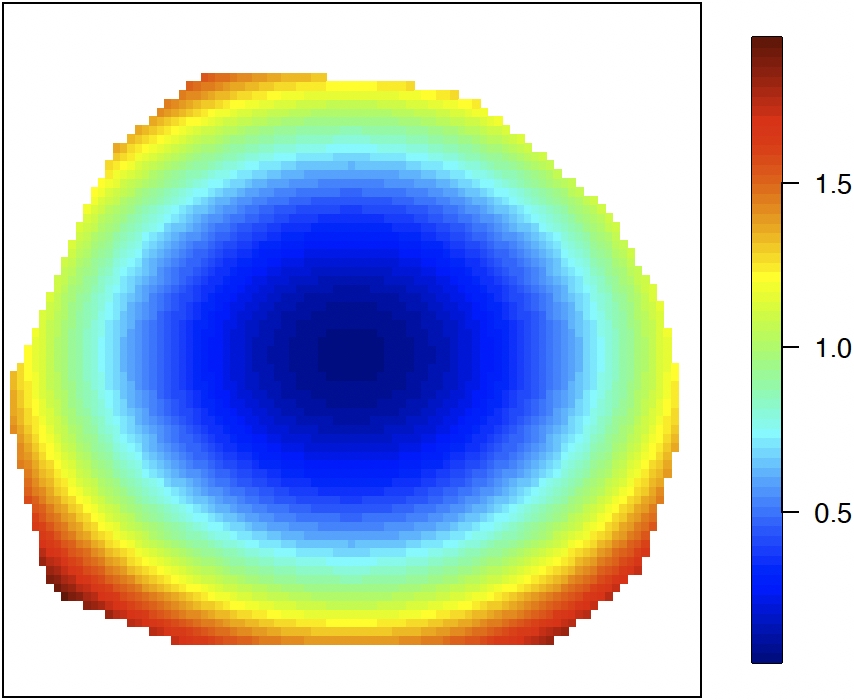}
	& \includegraphics[height=.52in]{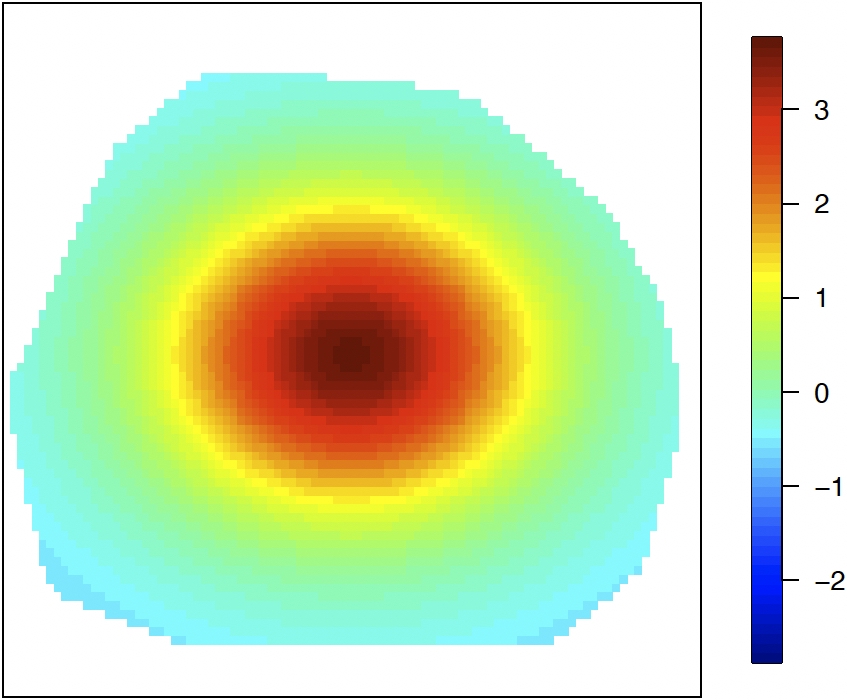} 
	& \includegraphics[height=.52in]{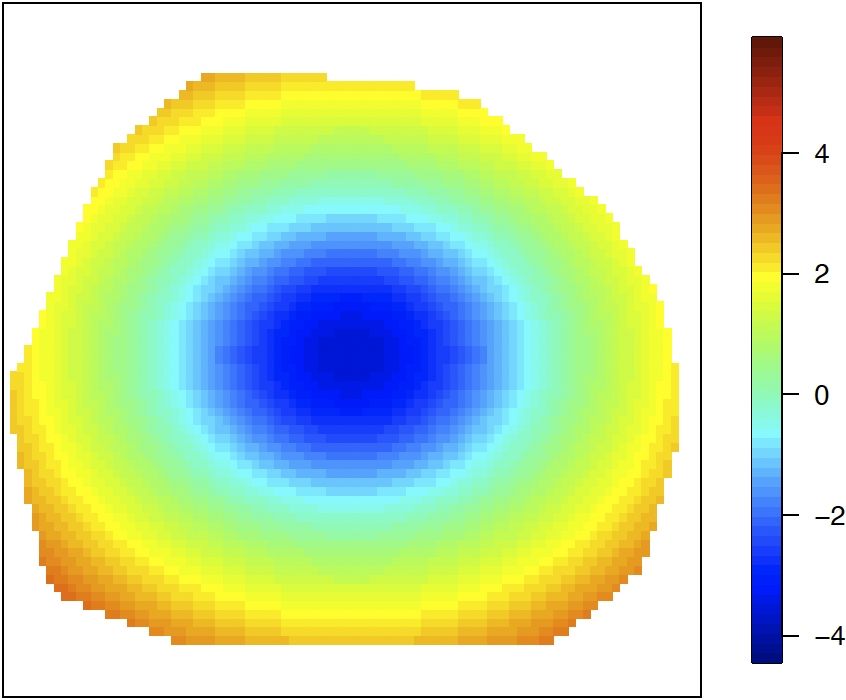}
	& \includegraphics[height=.52in]{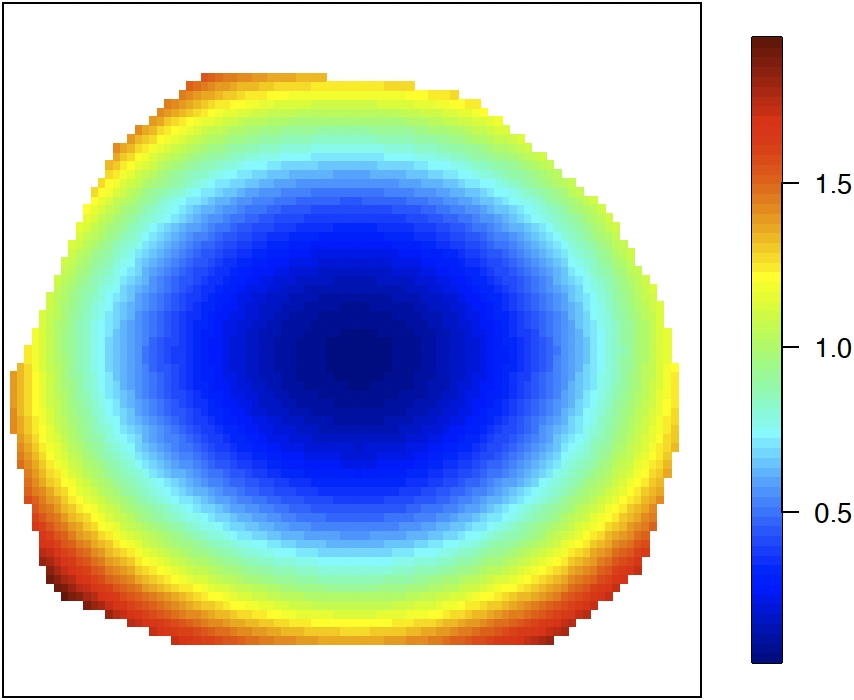} 
	& \includegraphics[height=.52in]{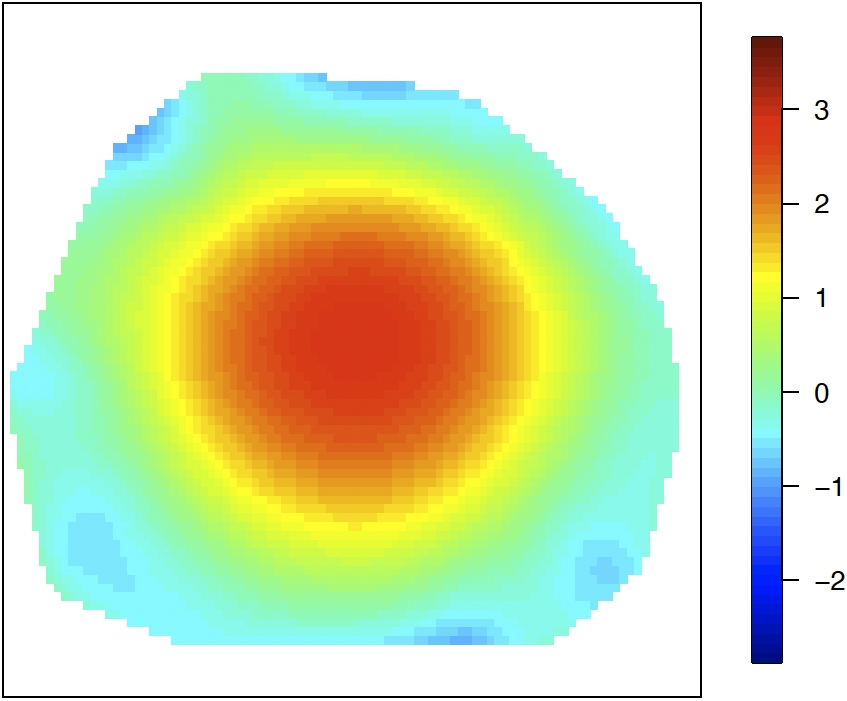}
	& \includegraphics[height=.52in]{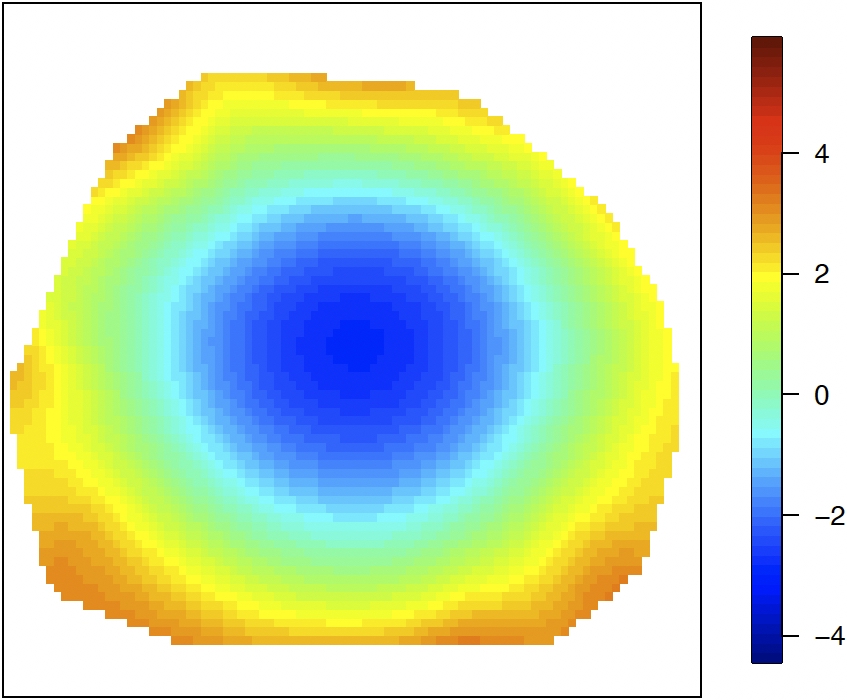} \\ 	 
\end{tabular}\vspace*{-.4cm}
\caption{Multi-planar views (transverse, coronal, sagittal) of true basis functions $\phi_1$ and $\phi_2$ and true coefficient function $\gamma_0$ and their estimates $\widehat{\phi}_1$, $\widehat{\phi}_2$ and $\widehat{\gamma}$ for a randomly selected iteration with $n=2000$ in the 3D imaging setting.}
\label{FIG:eg3-1}
\end{figure}

We next investigate the performance of the proposed method for 3D imaging data. In this scenario, for the Hilbert-valued component, we construct the imaging data as $Z_i(s_1,s_2,s_3) = \sum_{j=1}^{2} \lambda_j^{1/2}U_{ij}\phi_j(s_1,s_2,s_3)$ on a $79 \times 95 \times 66$ voxel grid, which mimics the typical resolution of well-registered and preprocessed PET images in ADNI database. The eigenvalues $(\lambda_1,\lambda_2)$ are set to be $(2,1)$, and $(U_{i1},U_{i2})^\top$ follows $\text{MVN}(0_2, I_2)$. The basis functions $\phi_1$ and $\phi_2$ are derived from $20\{(s_1-0.5)^2 + (s_2-0.5)^2 + (s_3-0.5)^2\}$ and $\exp[-15\{(s_1-0.5)^2 + (s_2-0.5)^2 + (s_3-0.5)^2\}]$ respectively, after orthonormalization. Figure \ref{FIG:eg3-1} (left two columns) displays the multi-planar representation of these basis functions. We set $\gamma_0 = 1.5\phi_1 - \phi_2$.

Similar to the analysis in Section \ref{SUBSEC:2D_ADNI}, we check the MSEs for estimation accuracy.
Figure \ref{FIG:eg3-1} demonstrates that with $n = 2000$ observations, our method achieves accuracy for a randomly selected single iteration. The estimated basis functions closely capture the main pattern and magnitude of the true basis functions, which further support the conclusions in Theorem \ref{THM:maxes}.

Table \ref{TAB:Eg3_MSE} reports the MSEs for the estimates of eigenvalues $\lambda_1$ and $\lambda_2$, regression coefficients $\alpha$, $\beta$ and $\gamma$, and the estimated number of PCA components $\widehat{m}_n$.
As the sample size increases from $n=100$ to $n=2000$, MSEs decrease significantly across all parameters, confirming estimation consistency. 
\begin{table}
\caption{Mean Squared Errors (MSEs$\times 10^{-2}$) of the estimates for AS-PCA eigenvalues $\lambda_1$ and $\lambda_2$, and the number of estimated principal components in the 3D imaging setting. Results are presented for various sample sizes ($n = 100, 500, 2000$) and correlation scenarios ($r = 0, 0.5$), based on $100$ Monte Carlo simulations. \label{TAB:Eg3_MSE}}
\begin{center}
\renewcommand{\arraystretch}{0.85}
\resizebox{0.4\textwidth}{!}{
\begin{tabular}{rrrrrrr}
	\hline
	\multirow{2}[0]{*}{$n$} & \multicolumn{3}{c}{$r=0$}     & \multicolumn{3}{c}{$r=0.05$} \\
	\cmidrule(lr){2-4} \cmidrule(lr){5-7}
	& $\lambda_1$ & $\lambda_2$ & $\widehat{m}_n$ & $\lambda_1$ & $\lambda_2$ & $\widehat{m}_n$\\
	\hline
	100 & 9.50 & 2.04 & 2 & 9.50 & 2.04 & 2\\
	500 & 1.43 & 0.83 & 2 & 1.43 & 0.83 & 2\\
	2000 & 0.38 & 0.67 & 2 & 0.38 & 0.67 & 2 \\
	\hline
\end{tabular}}%
\end{center}
\end{table}%

\begin{table}
\caption{Mean Squared Errors (MSEs$\times 10^{-2}$) and empirical coverage rates (in parentheses) of the 95\% HS-PCR bootstrap confidence interval for the estimates of coefficients $\alpha$, $\beta$, and $\gamma$ in the 3D imaging setting. Results are presented for various sample sizes ($n = 100, 500, 2000$) and correlation scenarios ($r = 0, 0.5$), based on $100$ Monte Carlo simulations. \label{TAB:Eg4_MSE}} 
\setlength{\tabcolsep}{3pt}
\renewcommand{\arraystretch}{0.85}
\resizebox{0.99\textwidth}{!}{
\begin{tabular}{rrrrrrrrrrrrrrr}
	\hline
	\multirow{2}[0]{*}{$n$} & \multicolumn{7}{c}{$r=0$}     & \multicolumn{7}{c}{$r=0.05$} \\
	\cmidrule(lr){2-8} \cmidrule(lr){9-15}
	& $\alpha$ & $\beta_1$ & $\beta_2$ & $\beta_3$ & $\beta_4$ & $\gamma_1$ & $\gamma_2$ & 
	$\alpha$ & $\beta_1$ & $\beta_2$ & $\beta_3$ & $\beta_4$ & $\gamma_1$ & $\gamma_2$\\
	\hline
	100 & 0.89 & 1.54 & 1.91 & 1.65 & 1.16 & 2.70 & 5.40 
	& 0.89 & 1.54 & 1.91 & 1.65 & 1.16 & 2.70 & 5.40 \\
	& (96\%) & (92\%) & (94\%) & (94\%) & (95\%) & (92\%) & (93\%) & (96\%) & (95\%) & (91\%) & (95\%) & (96\%) & (92\%) & (93\%)\\
	500 & 0.19 & 0.31 & 0.31 & 0.34 & 0.24 & 0.47 & 1.21 & 0.19 & 0.31 & 0.31 & 0.34 & 0.24 & 0.47 & 1.21\\
	& (95\%) & (97\%) & (95\%) & (91\%) & (95\%) & (95\%) & (94\%) & (95\%) & (94\%) & (95\%) & (98\%) & (97\%) & (95\%) & (94\%)\\
	2000 & 0.03 & 0.07 & 0.07 & 0.07 & 0.07 & 0.07 & 0.39 
	& 0.03 & 0.07 & 0.07 & 0.07 & 0.07 & 0.07 & 0.39\\\
	& (97\%) & (98\%) & (93\%) & (90\%) & (96\%) & (99\%) & (87\%) & (97\%) & (95\%) & (96\%) & (97\%) & (94\%) & (99\%) & (87\%)\\
	\hline
\end{tabular}}\vspace{-6pt}
\end{table}%

Figure \ref{FIG:eg3-1} shows the multi-planar representation of the coefficient $\gamma$ compared to its estimates.
PVE criterion produces estimates that accurately capture the true function's patterns across all three views, including the radial pattern and magnitude.

In summary, our numerical studies show that the proposed methodology delivers accurate estimation of both Euclidean and Hilbert-valued components, particularly for moderate to large sample sizes. Performance remains robust under moderate correlation among the Euclidean covariates, although estimation precision is somewhat reduced in these settings. These findings support the practical usefulness of our approach for neuroimaging applications.

\section{HS-PCR analysis of ADNI neuroimaging data}
\label{SEC:DA}

We apply HS-PCR to neuroimaging data from the Alzheimer's Disease Neuroimaging Initiative (ADNI), jointly analyzing PET imaging biomarkers, genetic risk factors, and demographic covariates to draw inferences about cognitive decline. The functional principal components of the PET imaging covariate are estimated via AS-PCA with a multivariate spline over triangulation basis (Section \ref{SEC:NeuroImple}); HS-PCR is then applied to the estimated scores. We consider two formulations: a standard linear regression setting (Section \ref{SUBSEC:LM}) and a precision-medicine extension incorporating imaging-by-treatment interactions (Section \ref{SUBSEC:PM}).

\subsection{HS-PCR: Linear Regression Formulation}
\label{SUBSEC:LM}

The data are drawn from the large neuroimaging datasets in ADNI (\url{http://adni.loni.usc.edu}).
The longitudinal cohort study in ADNI is a comprehensive neuroimaging study that collected a variety of necessary phenotypic measures, including structural, functional, and molecular neuroimaging, biomarkers, clinical and neuropsychological variables, and genomic information \citep{Petersen:etal:10,Weiner:Veitch:15}. 
These data provide unprecedented resources for statistical methods development and scientific discovery.

We now analyze the records from 441 participants through the ADNI1 and ADNI GO phases. The data contains the following variables:  \vspace{-6pt}
\begin{itemize}
	\item Mini-Mental State Examination (MMSE): response variable, integer scores from $15$ to $30$, with lower values indicating more severe AD. \vspace{-8pt}
	\item Fludeoxyglucose PET image: $79\times95$ voxel image of brain glucose metabolism,  with voxel intensities ranging from $0.013$ to $2.149$. \vspace{-8pt}
	\item Age: participant age in years ($55$--$89$). \vspace{-8pt}
	\item Education: years of education ($4$--$20$). \vspace{-8pt}
	\item Gender: indicator for female ($1$ = female, $0$ = male; $169$ females, $278$ males). \vspace{-8pt}
	\item Ethnicity: indicator for Hispanic/Latino ($1$ = Hispanic/Latino, $0$ = otherwise; $12$ Hispanic/Latino, $429$ non-Hispanic/Latino, $6$ unknown). \vspace{-8pt}
	\item Race: indicator for White ($1$ = White, $0$ = otherwise; $1$ Indian/Alaskan, $7$ Asian, $24$ Black, $413$ White, $2$ multiple). \vspace{-8pt}
	\item Marriage status: indicator for married ($1$ = married, $0$ = otherwise; $35$ divorced, $344$ married, $12$ never married, $56$ widowed). \vspace{-8pt}
	\item Apolipoprotein (APOE4) genotype: two indicators, APOE1 and APOE2, for carriers with one or two copies of the APOE4 allele, respectively (APOE4 count $0$--$2$). \vspace{-20pt}
\end{itemize} 

\begin{figure}[htbp]
\begin{center}
\begin{tabular}{ccccccc}
	\includegraphics[height=.52in]{figures/ZZeta1.jpg}
		& \includegraphics[height=.52in]{figures/ZZeta2.jpg} 
		& \includegraphics[height=.52in]{figures/ZZeta3.jpg}
		& \includegraphics[height=.52in]{figures/ZZeta4.jpg}
		& \includegraphics[height=.52in]{figures/ZZeta5.jpg} 
		& \includegraphics[height=.52in]{figures/ZZeta6.jpg}
		& \includegraphics[height=.52in]{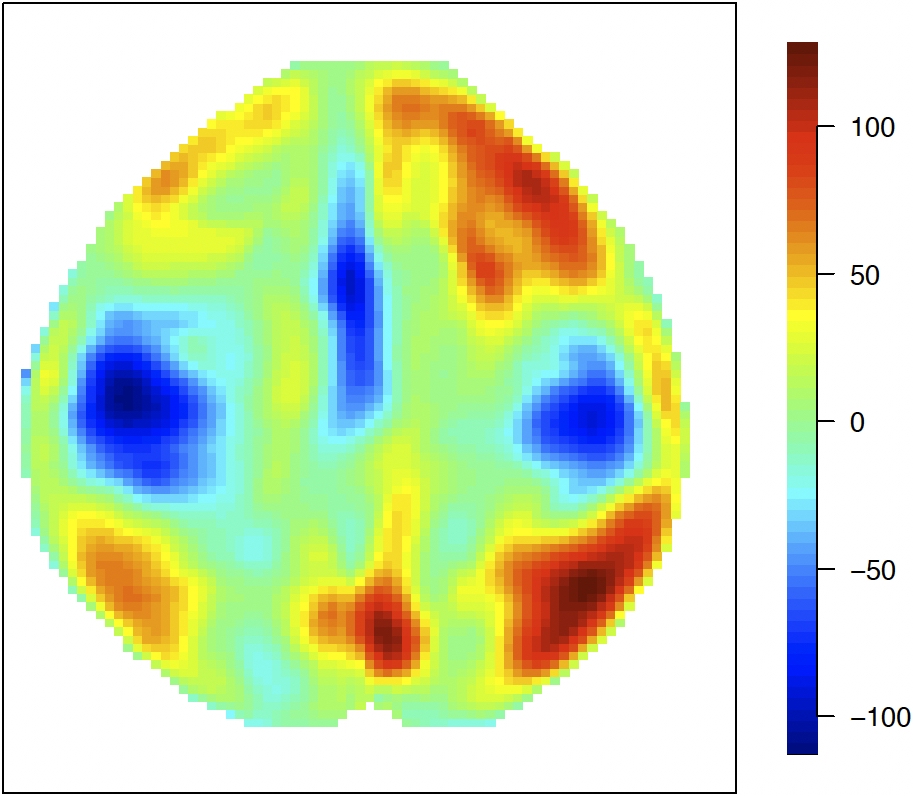} \vspace{-2pt}  \\ 
	$\widehat{\phi}_{1}$  & $\widehat{\phi}_{2}$  & $\widehat{\phi}_{3}$ & $\widehat{\phi}_{4}$  & $\widehat{\phi}_{5}$  & $\widehat{\phi}_{6}$ & $\widehat{\gamma}$ \\ 
\end{tabular}
\end{center}\vspace*{-.6cm}
\caption{Left six panels: selected leading principal component maps from the ADNI PET data. Right panel: estimated coefficient map $\widehat{\gamma}$ based on these PCs. }
\label{FIG:beta2hat}
\end{figure}

We apply the proposed method as follows. We first select the six leading principal component bases using the PVE criterion. The selected principal component maps are shown in the left three columns of Figure \ref{FIG:beta2hat}. We then estimate the regression model using these six PCs, with the estimated coefficient map for $\gamma$ shown in the right panel of Figure \ref{FIG:beta2hat}. Table \ref{TAB:Est_BCI} reports the estimated coefficients for the nonfunctional covariates and the corresponding 95\% bootstrap confidence intervals.

With increasing age, the mean MMSE scores decrease significantly, and higher education is associated with higher MMSE scores, both in line with previous findings \citep{Fratiglioni:Winblad:vonStrauss:07,Stern:12}. For the APOE4 gene, the MMSE scores decrease with the number of copies, reflecting the higher risk for earlier onset of AD \citep{Schneider:11}. Gender and ethnicity effects are not statistically significant, whereas race and marital status show moderate effects.

\begin{table}
\caption{Estimated coefficients (Esti) and corresponding 95\% bootstrap confidence intervals (BCI) in linear regression setting. \label{TAB:Est_BCI}}
\setlength{\tabcolsep}{3pt}
\renewcommand{\arraystretch}{0.85}
\begin{center}
\resizebox{0.8\textwidth}{!}{
\begin{tabular}{lrc | lrc | lrc}
	\hline
	Term & Esti & 95\% BCI &
	Term & Esti & 95\% BCI &
	Term & Esti & 95\% BCI \\
	\hline
	Intercept & 17.814 & (11.79, 24.59) 
		& Gender & -0.162 & (-0.65, 0.35) 
		& Ethnicity & -0.151 & (-1.12, 0.80) \\
	Age & -0.071 & (-0.11, -0.04) 
		& APOE1 & -0.572 & (-1.01, -0.13)
		&Race & 1.048 & (0.12, 1.95)\\
	Education & 0.219 & (0.15, 0.29) 
		& APOE2 & -1.529 & (-2.46, -0.65) 
		& Marriage & -0.398 & (-0.94, 0.17) \\
	\hline
\end{tabular}}
\end{center}
\end{table}%

\subsection{HS-PCR: Precision Medicine Formulation}
\label{SUBSEC:PM}

In addition to the linear regression setting of Section \ref{SUBSEC:LM}, we now incorporate treatment as a binary covariate and consider a precision-medicine formulation. 
To be specific, for the model $f(\alpha,\beta,\gamma)=\alpha+\beta^{\top} X + \langle \gamma,Z\rangle$,
we write
\begin{align*}
	\mathbb{E}(Y|X,Z,A) 
	&= f(\alpha_1,\beta_1,\gamma_1) + A f(\alpha_2,\beta_2,\gamma_2) \\
	&= \left\{\alpha_1+\beta_1^{\top} X + \langle \gamma_1,Z\rangle\right\} 
		+ A\left\{\alpha_2+\beta_2^{\top} X + \langle \gamma_2,Z\rangle\right\},
\end{align*}
where $A$ is the treatment. 
The first term models the baseline mean response, and the second term models the treatment effect and its interaction with covariates.

During the ADNI1 and ADNI GO study periods, the US FDA-approved therapies for AD symptoms included cholinesterase inhibitors and the NMDA-partial receptor antagonist memantine. Cholinesterase inhibitors, including donepezil, galantamine, and rivastigmine, are prescribed for mild-to-moderate-stage AD. Memantine is prescribed for the treatment of AD either as monotherapy or in combination with one of the cholinesterase inhibitors for moderate-to-severe stage AD \citep{Schneider:11}. 
We define $A=1$ for participants who received at least one FDA-approved cholinesterase inhibitor (donepezil, galantamine, rivastigmine) or memantine during the ADNI1/GO study period, and $A=0$  otherwise.

\begin{figure}[htbp]
\begin{center}
\begin{tabular}{cccc}
	\includegraphics[height=.62in]{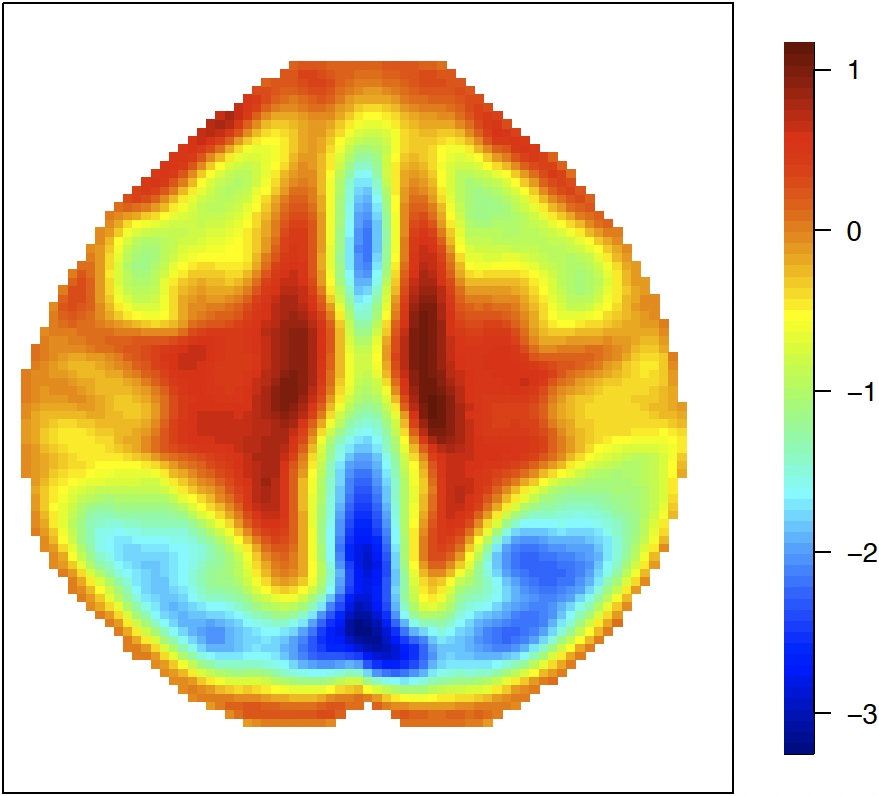}
		& \includegraphics[height=.62in]{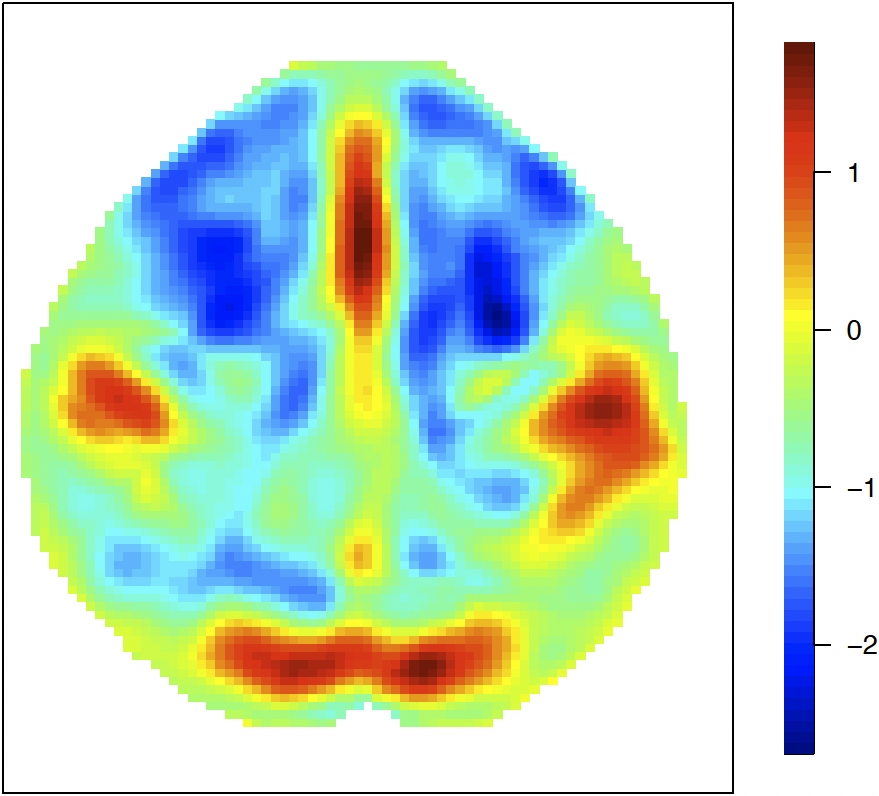} 
		& \includegraphics[height=.62in]{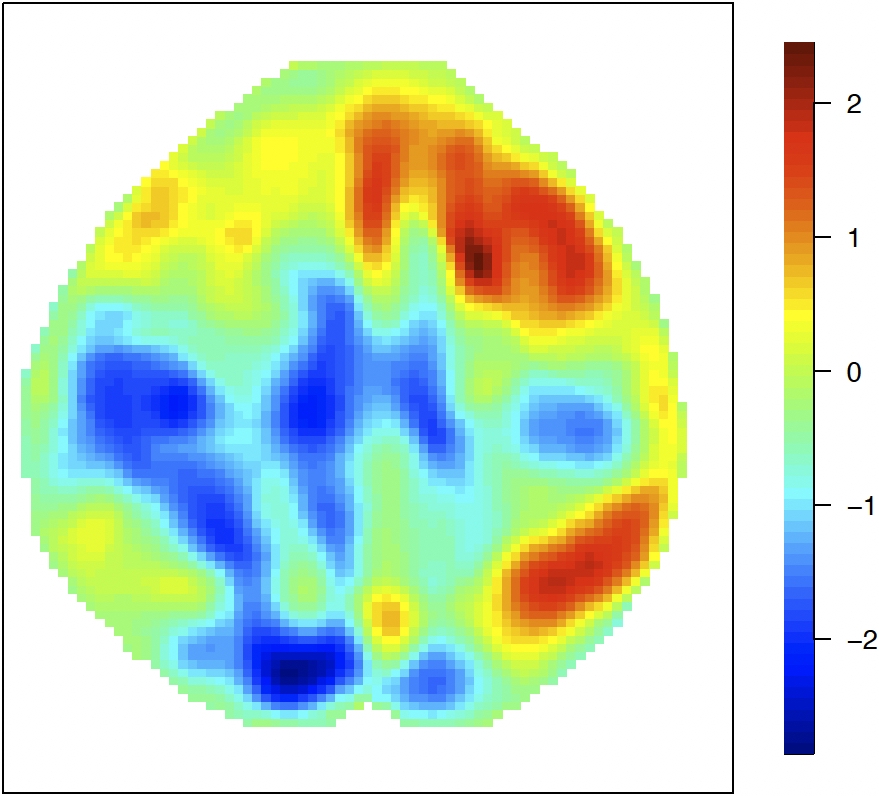}
		& \includegraphics[height=.62in]{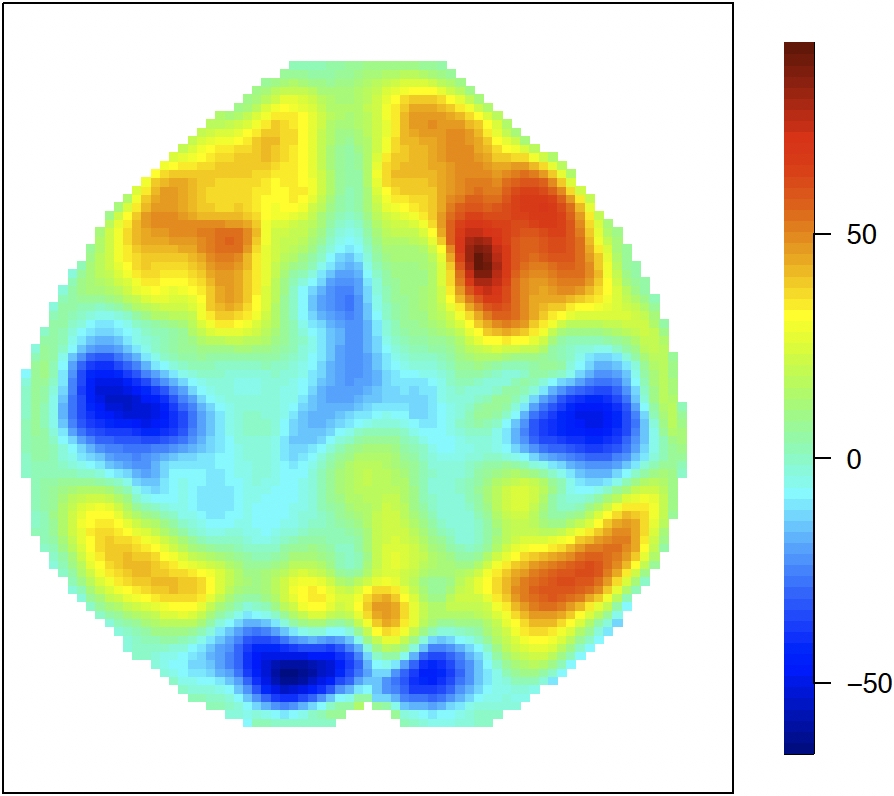} \vspace{-2pt} \\ 
	$\widehat{\phi}_{1}$ for $Z~~$ 
		& $\widehat{\phi}_{2}$ for $Z~~$ 
		& $\widehat{\phi}_{3}$ for $Z~~$ 
		& $\widehat{\gamma}_1~~$ \\
	\includegraphics[height=.62in]{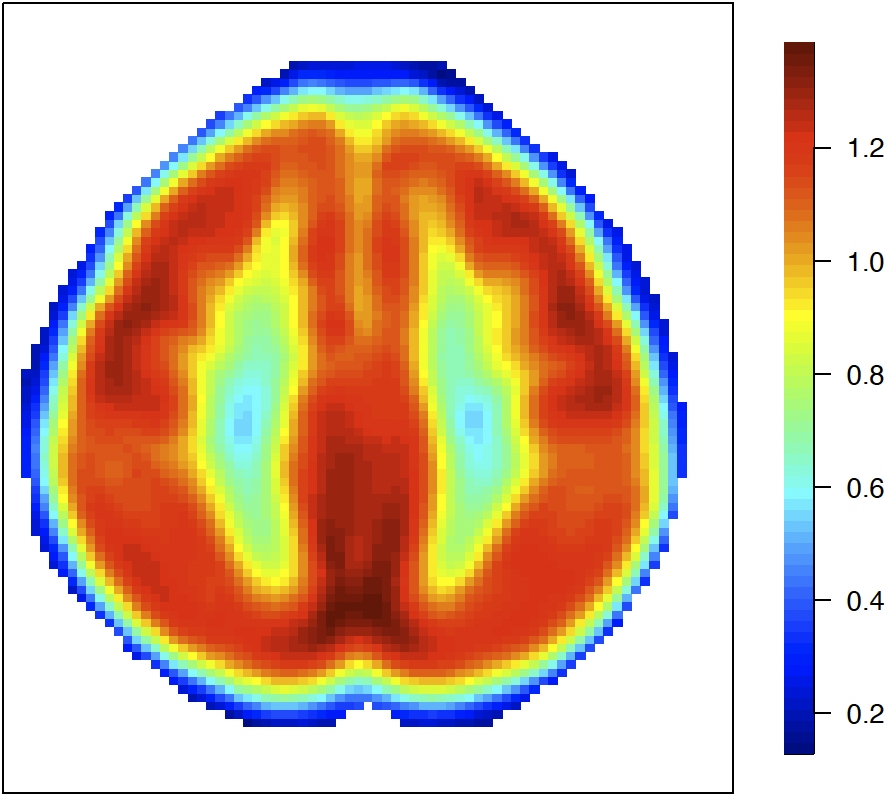}
		& \includegraphics[height=.62in]{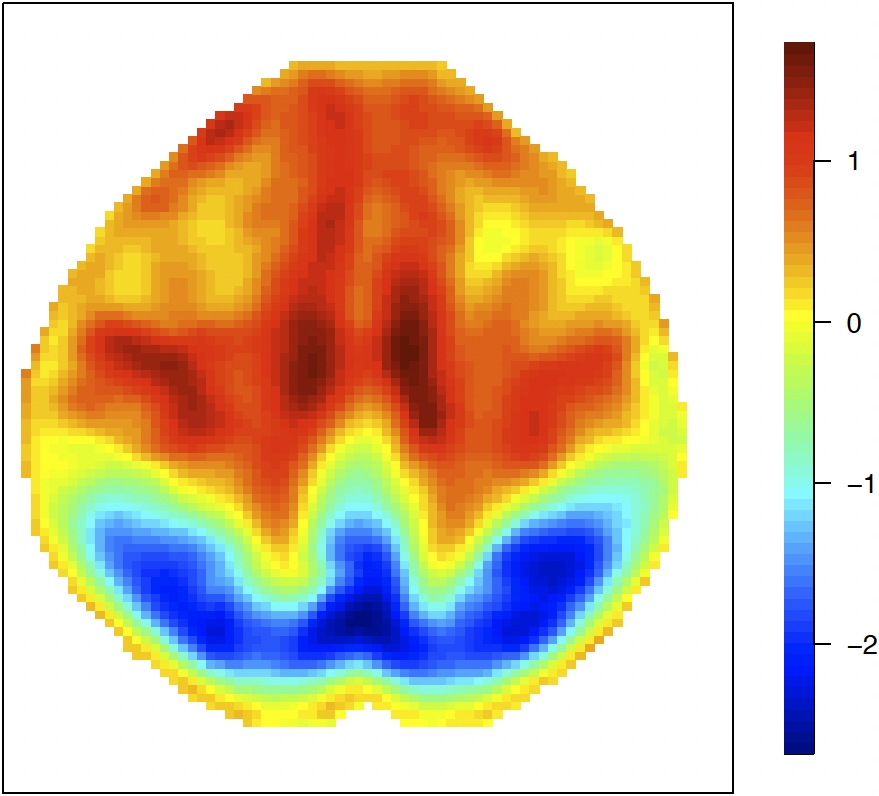} 
		& \includegraphics[height=.62in]{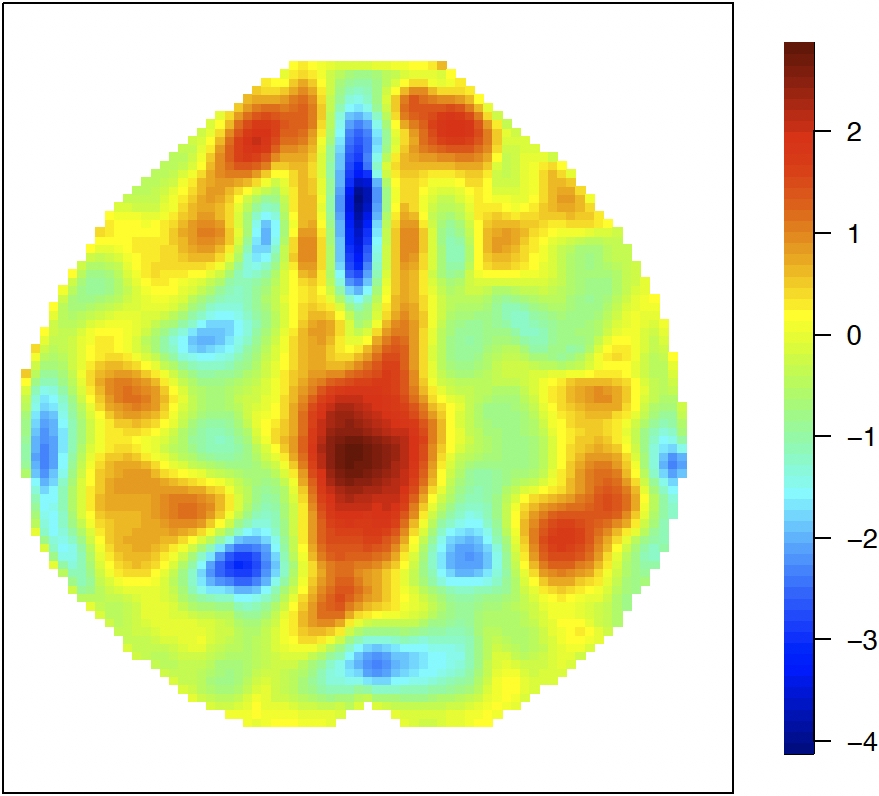} 
		& \includegraphics[height=.62in]{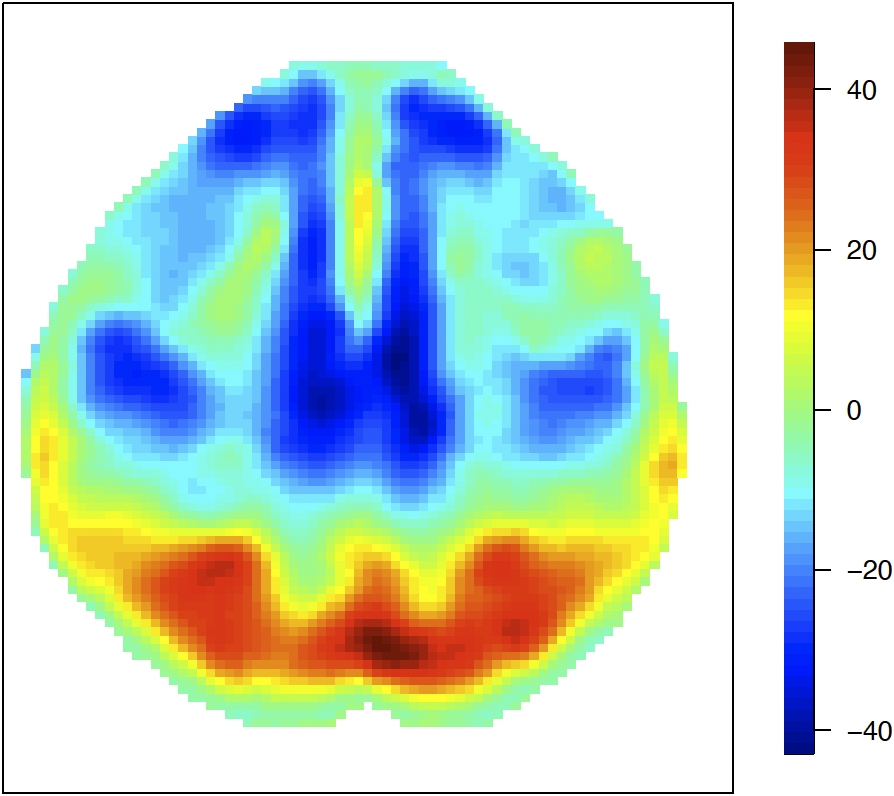} \vspace{-2pt}\\ 
	$\widehat{\phi}_{1}$ for $AZ~~$ 
		& $\widehat{\phi}_{2}$ for $AZ~~$ 
		& $\widehat{\phi}_{3}$ for $AZ~~$
	& $\widehat{\gamma}_2~~$\\
\end{tabular}
\end{center} \vspace*{-.6cm}
\caption{Left Three Columns: selected leading principal component basis maps for $Z$ (top) and $AZ$ (bottom). 
Right Column: Estimated coefficient maps of $\gamma_1$ (top) and $\gamma_2$ (bottom).}
\label{FIG:PM_phiHat_gammaHat}
\end{figure}

We again use the PVE criterion to select the leading principal components. The selected basis maps for $Z$ and the treatment–imaging interaction $AZ$ are shown in the left three columns of Figure \ref{FIG:PM_phiHat_gammaHat} (top and bottom rows, respectively). The estimated coefficient maps for $\gamma_1$ and $\gamma_2$ based on these PCs are displayed in the right column of Figure \ref{FIG:PM_phiHat_gammaHat}, and both reveal interpretable brain structures.
Table \ref{TAB:PM_Est_BCI} summarizes the estimated coefficients for the nonfunctional predictors and their 95\% bootstrap confidence intervals. The estimated main effect of treatment corresponds to an average improvement of approximately $9.85$ MMSE points, with conclusions broadly consistent with those from Section \ref{SUBSEC:LM}.

\begin{table}
\caption{Estimated coefficients and 95\% bootstrap confidence intervals (BCI) for the linear covariates $\alpha$ and $\beta$ for ADNI data.}
\label{TAB:PM_Est_BCI}%
\setlength{\tabcolsep}{3pt}
\renewcommand{\arraystretch}{0.85}
\begin{center}
\resizebox{0.8\textwidth}{!}{
\begin{tabular}{lrc | lrc}
	\hline
	Term & Estimate & 95\% BCI & Term & Estimate & 95\% BCI \\
	\hline
	Intercept & 12.488 & (4.59, 20.68) & Treatment & 9.850 & (-4.03, 23.34) \\
	Age & 0.005 & (-0.03, 0.04) & Treatment$\times$Age & -0.120 & (-0.18, -0.07)  \\
	Education & 0.228 & (0.13, 0.31) & Treatment$\times$Education & -0.032 & (-0.16, 0.12) \\
	Gender & 0.014 & (-0.62, 0.65) & Treatment$\times$Gender & -0.090 & (-1.02, 0.77) \\
	APOE1 & -0.785 & (-1.43, -0.13) & Treatment$\times$APOE1 & 0.529 & (-0.44, 1.40) \\
	APOE2 & -0.362 & (-1.77, 0.74) & Treatment$\times$APOE2 & -1.141 & (-2.71, 0.54) \\
	Ethnicity & -0.396 & (-1.77, 0.95) & Treatment$\times$Ethnicity & 0.657 & (-2.26, 3.92) \\
	Race & 0.120 & (-0.80, 1.12) & Treatment$\times$Race & 2.594 & (1.01, 4.25) \\
	Marriage & 0.070 & (-0.61, 0.81) & Treatment$\times$Marriage & -0.214 & (-1.34, 0.93) \\
	\hline
\end{tabular}
}
\end{center}
\vspace{-10pt}
\end{table}%

\section{Conclusions and discussions}
\label{SEC:Conclusions}

We have developed Adaptive Subspace PCA (AS-PCA), a framework for principal component analysis of random elements in a general separable Hilbert space, and applied it to construct Hilbert-Space Principal Component Regression (HS-PCR). The theoretical contributions of AS-PCA, including a Donsker theorem under second-moment conditions only, joint Gaussian limits for leading eigenpairs, a data-driven projection accuracy diagnostic, and bootstrap validity, together provide the first complete inferential pipeline for PCA and regression with abstract Hilbert-space-valued covariates, without requiring kernel specification, Gaussianity, or effective rank constraints.

HS-PCR translates the theoretical guarantees of AS-PCA into a practical regression tool: root-$n$ consistent and asymptotically normal estimation of all regression parameters, with bootstrap confidence regions valid for scalar, Euclidean, and functional coefficients simultaneously. As demonstrated in Section \ref{SEC:DA}, HS-PCR accommodates complex neuroimaging domains through the multivariate spline over triangulation implementation of AS-PCA, and extends naturally to precision-medicine formulations involving imaging-by-treatment interactions.

An important direction for future work is to extend AS-PCA to settings with incomplete or noisy observations of $Z$, where the projection scores $\langle \psi_\ell, Z_i \rangle$ are not directly observed. Extending the AS-PCA diagnostic and inferential theory to this setting would substantially broaden the applicability of HS-PCR to sparse functional data problems.

\section*{Acknowledgements}

Research reported in this publication was supported in part by the National Institute of General Medical Sciences of the National Institutes of Health under Award Number P20GM139769 (Xinyi Li), National Science Foundation awards DMS-2210658 (Xinyi Li) and DMS-2210659 (Michael R.\ Kosorok). The content is solely the responsibility of the authors and does not necessarily represent the official views of the National Institutes of Health. The investigators within the ADNI contributed to the design and implementation of ADNI and/or provided data but did not participate in analysis or writing of this report. A complete listing of ADNI investigators can be found at \url{http://adni.loni.usc.edu/wp-content/uploads/how_to_apply/ADNI_Acknowledgement_List.pdf}.

\bibliography{references}

@book{Lai:Schumaker:07,
	Author = {Lai, Ming-Jun and Schumaker, L. L.},
	Publisher = {Cambridge University Press},
	Title = {Spline functions on triangulations.},
	Year = {2007}}

@article{Lai:Wang:13,
	Author = {Lai, Ming-Jun and Wang, Li},
	Journal = {Statistica Sinica},
	Pages = {1399--1417},
	Publisher = {JSTOR},
	Title = {Bivariate penalized splines for regression},
	Volume = {23},
	Number = {3},
	Year = {2013}}

@article{Zhu:Fan:Kong:14,
	Author = {Zhu, Hongtu and Fan, Jianqing and Kong, Linglong},
	Journal = {Journal of the American Statistical Association},
	Number = {507},
	Pages = {1084--1098},
	Publisher = {Taylor \& Francis},
	Title = {Spatially varying coefficient model for neuroimaging data with jump discontinuities},
	Volume = {109},
	Year = {2014}}

@article{Yuan:Cai:10,
    Author={Yuan, Ming and Cai, T. Tony},
    Journal={Annals of Statistics},
    Volume={38},
    Number={6},
    Pages={3412-3444},
    Title={A reproducing kernel Hilbert space approach to functional linear regression},
    Year={2010}}

@book{Ramsay:Silverman:05,
	Author = {Ramsay, James O and Silverman, Bernard W},
	Publisher = {New York, Springer.},
	Title = {Functional Data Analysis},
	Year = {2005}}

@article{Morris:15,
	Author = {Morris, Jeffrey S.},
	Journal = {Annual Reviews of Statistics and its Application},
	Pages = {321--359},
	Title = {Functional regression},
	Volume = {2},
	Year = {2015}}

@book{Kosorok:08,
	title={Introduction to Empirical Processes and Semiparametric Inference},
	author={Kosorok, Michael R},
	year={2008},
	publisher={Springer: New York}
}

@article{Bateman:etal:12,
  title={Clinical and biomarker changes in dominantly inherited Alzheimer's disease},
  author={Bateman, Randall J and Xiong, Chengjie and Benzinger, Tammie LS and Fagan, Anne M and Goate, Alison and Fox, Nick C and Marcus, Daniel S and Cairns, Nigel J and Xie, Xianyun and Blazey, Tyler M and Holtzman, David M. and Santacruz,  Anna and Buckles, Virginia and Oliver, Angela and Moulder, Krista and Aisen, Paul S. and Ghetti, Bernardino and Klunk, William E. and McDade, Eric and Martins, Ralph N. and Masters, Colin L. and Mayeux, Richard and Ringman, John M. and Rossor, Martin N. and Schofield, Peter R. and Sperling, Reisa A. and Salloway, Stephen and Morris, John C.},
  journal={New England Journal of Medicine},
  volume={367},
  number = {9},
  pages={795--804},
  year={2012},
  publisher={Mass Medical Soc}
}

@article{Petersen:etal:10,
  title={Alzheimer's disease neuroimaging initiative ({ADNI}): clinical characterization},
  author={Petersen, Ronald Carl and Aisen, PS and Beckett, Laurel A and Donohue, MC and Gamst, AC and Harvey, Danielle J and Jack, CR and Jagust, WJ and Shaw, LM and Toga, AW and Trojanowski, J.Q.},
  journal={Neurology},
  volume={74},
  number={3},
  pages={201--209},
  year={2010},
  publisher={AAN Enterprises}
}

@article{Weiner:Veitch:15,
  title={Introduction to special issue: overview of {A}lzheimer's {D}isease {N}euroimaging {I}nitiative},
  author={Weiner, Michael W and Veitch, Dallas P},
  journal={Alzheimer's \& Dementia},
  volume={11},
  number={7},
  pages={730--733},
  year={2015},
  publisher={Elsevier}
}

@article{Schneider:11,
  title={Treatment with cholinesterase inhibitors and memantine of patients in the Alzheimer's Disease Neuroimaging Initiative},
  author={Schneider, Lon S and Insel, Philip S and Weiner, Michael W and Alzheimer's Disease Neuroimaging Initiative and others},
  journal={Archives of neurology},
  volume={68},
  number={1},
  pages={58--66},
  year={2011},
  publisher={American Medical Association}
}

@article{Fratiglioni:Winblad:vonStrauss:07,
  title={Prevention of Alzheimer's disease and dementia. Major findings from the Kungsholmen Project},
  author={Fratiglioni, Laura and Winblad, Bengt and von Strauss, Eva},
  journal={Physiology \& behavior},
  volume={92},
  number={1-2},
  pages={98--104},
  year={2007},
  publisher={Elsevier}
}

@article{Stern:12,
    title = {Cognitive reserve in ageing and Alzheimer's disease},
    author = {Yaakov Stern},
    journal = {The Lancet Neurology},
    volume = {11},
    number = {11},
    pages = {1006-1012},
    year = {2012},
    issn = {1474-4422},
    doi = {https://doi.org/10.1016/S1474-4422(12)70191-6},
    url = {https://www.sciencedirect.com/science/article/pii/S1474442212701916},
}

@article{Kong:Staicu:Maity:16,
  title={Classical testing in functional linear models},
  author={Kong, Dehan and Staicu, Ana-Maria and Maity, Arnab},
  journal={Journal of nonparametric statistics},
  volume={28},
  number={4},
  pages={813--838},
  year={2016},
  publisher={Taylor \& Francis}
}

@article{Wang:Chiou:Muller:16,
  title={Functional data analysis},
  author={Wang, Jane-Ling and Chiou, Jeng-Min and M{\"u}ller, Hans-Georg},
  journal={Annual Review of Statistics and its application},
  volume={3},
  pages={257--295},
  year={2016},
  publisher={Annual Reviews}
}

@article{Li:etal:24,
  title={Nonparametric Regression for 3D Point Cloud Learning},
  author={Li, Xinyi and Yu, Shan and Wang, Yueying and Wang, Guannan and Wang, Li and Lai, Ming-Jun},
  journal={Journal of Machine Learning Research},
  volume={25},
  number={102},
  pages={1--56},
  year={2024}
}

@article{Dai:Muller:18,
  title={Principal component analysis for functional data on Riemannian manifolds and spheres},
  author={Dai, Xiongtao and M{\"u}ller, Hans-Georg},
  journal={Annals of Statistics},
  volume={46},
  number={6B},
  pages={3334--3361},
  year={2018}
}

@article{Lin:Yao:19,
  title={Intrinsic Riemannian functional data analysis},
  author={Lin, Zhenhua and Yao, Fang},
  journal={Annals of Statistics},
  volume={47},
  number={6},
  pages={3533--3577},
  year={2019}
}

@article{Kim:etal:20,
  title={Principal component analysis for Hilbertian functional data},
  author={Kim, Dongwoo and Lee, Young Kyung and Park, Byeong U},
  journal={CSAM (Communications for Statistical Applications and Methods)},
  volume={27},
  number={1},
  pages={149--161},
  year={2020}
}

@article{Perry:etal:25,
	title = {Inference on the Proportion of Variance Explained in Principal Component Analysis},	
	author = {Perry, Ronan and Panigrahi, Snigdha and Bien, Jacob and Witten, Daniela},
	journal = {Journal of the American Statistical Association},
	volume = {0},
	number = {0},
	pages = {1--11},
	year = {2025},
	publisher = {Taylor \& Francis},
	doi = {10.1080/01621459.2025.2538895}
}

@book{Hsing:Eubank:15,
	Author = {Hsing, Tailen and Eubank, Randall},
	Publisher = {John Wiley \& Sons},
	Title = {Theoretical foundations of functional data analysis, with an introduction to linear operators},
	Year = {2015}}

@article{Li:Wang:Wang:21,
	Author = {Li, Xinyi and Wang, Lily and Wang, Huixia Judy},
	Title = {Sparse Learning and Structure Identification for Ultrahigh-Dimensional Image-on-Scalar Regression},
	Journal = {Journal of the American Statistical Association},
	Volume = {116},
	Number = {536},
	Pages = {1994--2008},
	Year = {2021},
	Publisher = {Taylor \& Francis},
	doi = {10.1080/01621459.2020.1753523}
}

@article{Zhang:etal:23,
  title={Image response regression via deep neural networks},
  author={Zhang, Daiwei and Li, Lexin and Sripada, Chandra and Kang, Jian},
  journal={Journal of the Royal Statistical Society Series B: Statistical Methodology},
  volume={85},
  number={5},
  pages={1589--1614},
  year={2023},
  publisher={Oxford University Press US}
}

@article{Shi:etal:22,
  title={Two-dimensional functional principal component analysis for image feature extraction},
  author={Shi, Haolun and Yang, Yuping and Wang, Liangliang and Ma, Da and Beg, Mirza Faisal and Pei, Jian and Cao, Jiguo},
  journal={Journal of Computational and Graphical Statistics},
  volume={31},
  number={4},
  pages={1127--1140},
  year={2022},
  publisher={Taylor \& Francis}
}

@article{Happ:Greven:18,
  title={Multivariate functional principal component analysis for data observed on different (dimensional) domains},
  author={Happ, Clara and Greven, Sonja},
  journal={Journal of the American Statistical Association},
  volume={113},
  number={522},
  pages={649--659},
  year={2018},
  publisher={Taylor \& Francis}
}

@article{Lila:Aston:Sangalli:16,
  title={Smooth principal component analysis over two-dimensional manifolds with an application to neuroimaging},
  author={Lila, Eardi and Aston, John AD and Sangalli, Laura M},
  journal={The Annals of Applied Statistics},
  volume={10},
  number={4},
  pages={1854--1879},
  year={2016}
}

@article{Chen:Goldsmith:Ogden:19,
  title={Functional data analysis of dynamic PET data},
  author={Chen, Yakuan and Goldsmith, Jeff and Ogden, R Todd},
  journal={Journal of the American Statistical Association},
  volume={114},
  number={526},
  pages={595--609},
  year={2019},
  publisher={Taylor \& Francis}
}

@article{Dauxois:etal:82,
    title = {Asymptotic theory for the principal component analysis of a vector random function: Some applications to statistical inference},
    journal = {Journal of Multivariate Analysis},
    volume = {12},
    number = {1},
    pages = {136-154},
    year = {1982},
    issn = {0047-259X},
    doi = {https://doi.org/10.1016/0047-259X(82)90088-4},
    url = {https://www.sciencedirect.com/science/article/pii/0047259X82900884},
    author = {J. Dauxois and A. Pousse and Y. Romain}
}

@book{Bosq:00,
  title={Linear processes in function spaces: theory and applications},
  author={Bosq, Denis},
  volume={149},
  year={2000},
  publisher={Springer Science \& Business Media}
}

@article{Cardot:etal:99,
	title={Functional linear model},
	author={Cardot, Herv{\'e} and Ferraty, Fr{\'e}d{\'e}ric and Sarda, Pascal},
	journal={Statistics \& Probability Letters},
	volume={45},
	number={1},
	pages={11--22},
	year={1999},
	publisher={Elsevier}
}

@article{Cardot:etal:03,
	title={Spline estimators for the functional linear model},
	author={Cardot, Herv{\'e} and Ferraty, Fr{\'e}d{\'e}ric and Sarda, Pascal},
	journal={Statistica Sinica},
	pages={571--591},
	year={2003},
	volume = {13},
	number = {3},
	publisher={JSTOR}
}

@article{Hall:HosseiniNasab:06,
	title={On properties of functional principal components analysis},
	author={Hall, Peter and Hosseini-Nasab, Mohammad},
	journal={Journal of the Royal Statistical Society Series B: Statistical Methodology},
	volume={68},
	number={1},
	pages={109--126},
	year={2006},
	publisher={Oxford University Press}
}

@article{Koltchinskii:Lounici:17:Bernoulli,
	title={Concentration inequalities and moment bounds for sample covariance operators},
	author={Koltchinskii, Vladimir and Lounici, Karim},
	journal={Bernoulli},
	volume = {23},
	number = {1},
	pages={110--133},
	year={2017},
	publisher={JSTOR}
}

@article{Koltchinskii:Lounici:17:AOS,
author = {Vladimir Koltchinskii and Karim Lounici},
title = {{Normal approximation and concentration of spectral projectors of sample covariance}},
volume = {45},
journal = {The Annals of Statistics},
number = {1},
publisher = {Institute of Mathematical Statistics},
pages = {121 -- 157},
year = {2017},
doi = {10.1214/16-AOS1437},
URL = {https://doi.org/10.1214/16-AOS1437}
}

@article{Koltchinskii:Lounici:17:Sankhya,
  title={New asymptotic results in principal component analysis},
  author={Koltchinskii, Vladimir and Lounici, Karim},
  journal={Sankhya A},
  volume={79},
  number={2},
  pages={254--297},
  year={2017},
  publisher={Springer}
}

@article{Koltchinskii:Lounici:16,
author = {Vladimir Koltchinskii and Karim Lounici},
title = {{Asymptotics and concentration bounds for bilinear forms of spectral projectors of sample covariance}},
volume = {52},
journal = {Annales de l'Institut Henri Poincaré, Probabilités et Statistiques},
number = {4},
publisher = {Institut Henri Poincaré},
pages = {1976 -- 2013},
year = {2016},
doi = {10.1214/15-AIHP705},
URL = {https://doi.org/10.1214/15-AIHP705}
}

@article{Koltchinskii:Loffler:Nickl:20,
author = {Vladimir Koltchinskii and Matthias L{\"o}ffler and Richard Nickl},
title = {{Efficient estimation of linear functionals of principal components}},
volume = {48},
journal = {The Annals of Statistics},
number = {1},
publisher = {Institute of Mathematical Statistics},
pages = {464 -- 490},
year = {2020},
doi = {10.1214/19-AOS1816},
URL = {https://doi.org/10.1214/19-AOS1816}
}

@book{vanderVaart:Wellner:23,
	Author = {van der Vaart, A W and Wellner, Jon A},
	Publisher = {Springer International Publishing},
	Title = {Weak Convergence and Empirical Processes: With Applications to Statistics},
	Year = {2023}
}


\clearpage

\setcounter{section}{0}
\setcounter{subsection}{0}
\setcounter{equation}{0}
\setcounter{figure}{0}
\setcounter{table}{0}
\setcounter{thmcounter}{0}

\renewcommand{\theequation}{S.\arabic{equation}}
\renewcommand{\thesection}{S.\arabic{section}}
\renewcommand{\thesubsection}{S.\arabic{section}.\arabic{subsection}}
\renewcommand{\thefigure}{S.\arabic{figure}}
\renewcommand{\thetable}{S.\arabic{table}}
\renewcommand{\thethmcounter}{S.\arabic{thmcounter}}

\makeatletter
\renewenvironment{theorem}[1][]%
  {\refstepcounter{thmcounter}\par\medskip\noindent%
   \textbf{{\sc Theorem} \thethmcounter%
   \ifx\\#1\\\else\ (#1)\fi.}\itshape}%
  {\par\medskip}
\renewenvironment{corollary}[1][]%
  {\refstepcounter{thmcounter}\par\medskip\noindent%
   \textbf{{\sc Corollary} \thethmcounter%
   \ifx\\#1\\\else\ (#1)\fi.}\itshape}%
  {\par\medskip}
\renewenvironment{remark}[1][]%
  {\refstepcounter{thmcounter}\par\medskip\noindent%
   \textbf{{\sc Remark} \thethmcounter%
   \ifx\\#1\\\else\ (#1)\fi.}\itshape}%
  {\par\medskip}
\makeatother

\begin{center}
{\large Supplemental Materials for ``Linear Regression Using Principal Components\\[2pt]
from General Hilbert-Space-Valued Covariates''}
\end{center}

\bigskip
\normalsize In the supplement, we provide the technical details of the conclusions presented in the main paper. Specifically, we give detailed proofs of Theorems~\ref{THM:expansion}--\ref{THM:bootstrap} in the paper, and present a number of technical lemmas and other supporting results used in the proofs.

\section{Technical Details}

\subsection{Proof of Theorem \ref{THM:expansion}}

The expansion and properties of $\{\lambda_j\}$, $\{\phi_j\}$ and $\{U_j\}$ follow from the Karhunen-Lo\`{e}ve expansion theorem for separable random variables in Hilbert spaces. Existence of the mean follows from $\mathbb{E}\|Z\|^2<\infty$ and so does the finiteness of $\sum_{j=1}^\infty \lambda_j$.

\subsection{Proof of Theorem \ref{THM:product}}

In the following, we use the subscript $0$ in the space to denote the separable space.

\begin{proof}[Proof of Theorem~\ref{THM:product}.]
By Theorem~\ref{THM:expansion}, each $Z_k$ can be written as $Z_k = \sum_{j=1}^\infty \lambda_{kj}^{1/2}U_{kj}\phi_{kj}$, where $\infty > \lambda_{k1} \geq \lambda_{k2} \geq \cdots \geq 0$, $\{\phi_{k1}, \phi_{k2}, \cdots\}$ is an orthonormal basis on $\mathcal{H}_k$, $\mathbb{E} U_{kj} = 0$, $\mathbb{E} U_{kj}^2 = 1$, and $\mathbb{E}U_{kj}U_{kj'} = 0$, for all $1 \leq j \ne j' <\infty$ and all $1 \leq k \leq K$. Denote $\mathcal{H}_{k0} \subset \mathcal{H}_k$ as the separable subspace spanned by $\{\phi_{k1}, \phi_{k2}, \cdots\}$, and let $\mathcal{B}_{k0} = \mathcal{B}_k \cap \mathcal{H}_{k0}$, $1 \leq k \leq K$.
Since $\langle c, Z_k \rangle_k = 0$ almost surely for all $c \in \mathcal{H}_k \backslash \mathcal{H}_{k0}$, define
\[
    \mathcal{F}_0 = \left\{f(Z)=\prod_{k=1}^K \langle a_k, Z_k\rangle_k : a_k \in \mathcal{B}_{k0}, 1 \leq k \leq K\right\},
\]
we have $\mathcal{F} = \mathcal{F}_0$ almost surely. Since $F = \prod_{k=1}^K \|Z_k\|_k$ is an envelope for $\mathcal{F}$, and $\mathbb{E} F^2 <\infty$ by assumption, we have all finite subsets $\mathcal{G} \subset \mathcal{F}$ are Donsker; and thus, we have that all finite-dimensional distributions of $\{\mathbb{G}_n g : g \in \mathcal{G}\}$ converge. Thus, by Theorem~2.1 and Lemma~7.20 of \citet{Kosorok:08}, our proof will become complete if we can show that for every $\epsilon, \eta > 0$, there exists a finite partition $\mathcal{F}_0 = \bigcup_{m = 1}^M \mathcal{F}_m$, so that $M < \infty$ and
\begin{equation}
\label{EQN:zerostar}
    \limsup_{n \rightarrow \infty} P^\ast\!\left(\sup_{1 \leq m \leq M} \sup_{f, g \in \mathcal{F}_m} \left|\mathbb{G}_n(f-g) \right| > \epsilon\right) < \eta,
\end{equation}
where $P^\ast$ is the outer probability. Thus $\mathcal{F}_0$, and hence also $\mathcal{F}$, is Donsker.

We now prove~(\ref{EQN:zerostar}). First, let $a = (a_1, \cdots, a_K)$ and $b = (b_1, \cdots, b_K)$ satisfy $a_k, b_k \in \mathcal{B}_{k0}$, $1 \leq k \leq K$. Then $a_k=\sum_{j=1}^{\infty} a_{kj}\phi_{kj}$ and $b_k=\sum_{j=1}^{\infty} b_{kj}\phi_{kj}$, where $a_{kj} = \langle \phi_{kj}, a_k \rangle_k$ and $b_{kj} = \langle \phi_{kj}, b_k \rangle_k$, respectively. Consequently,
\begin{align*}
&\mathbb{G}_n \left(\prod_{k=1}^K\langle a_k, Z_k\rangle_k - \prod_{k=1}^K \langle b_k, Z_k \rangle_k\right) \\
=\;& \sum_{k=1}^K \mathbb{G}_n\left\{\left(\prod_{k'=1}^{k-1}\langle b_{k'}, Z_{k'} \rangle_{k'}\right)\langle a_k - b_k, Z_k\rangle_k\left(\prod_{k' = k+1}^K \langle a_{k'}, Z_{k'} \rangle_{k'}\right)\right\}\\
=\;& \sum_{k=1}^K \sum_{(j_1,\ldots,j_K)=(1, \cdots, 1)}^{(\infty, \cdots, \infty)}
    \!\left\{\!\left(\prod_{k'=1}^{k-1}b_{{k'} j_{k'}}\!\right)
        \!\left(a_{k j_k}\! -\! b_{k j_k}\right)\!
        \left(\!\prod_{k'=k+1}^K a_{{k'}j_{k'}}\right)\!\right\}\!
    \mathbb{G}_n\!\left(\prod_{k'=1}^K U_{k' j_{k'}}^\ast\!\right) \!\equiv\! D_1,
\end{align*}
where $U_{k' j_{k'}}^{\ast} = \langle \phi_{k' j_{k'}}, Z_{k'}\rangle_{k'}=\lambda_{k' j_{k'}}^{1/2}U_{k' j_{k'}}$ and products over empty sets have the value 1. Continuing,
\begin{align*}
    |D_1|
    \leq \sum_{k=1}^K
    \left[\sum_{(1,\cdots, 1)}^{(\infty, \cdots, \infty)}
        \left\{\left(\prod_{k'=1}^{k-1} b_{k' j_{k'}}^2\right)
        \left(a_{k j_k} - b_{k j_k}\right)^2
        \left(\prod_{k'=k+1}^K a_{k' j_{k'}}^2\right)\right\}\right]^{1/2}\!\! \times
    \left[\sum_{(1,\cdots, 1)}^{(\infty, \cdots, \infty)}
        \left\{\mathbb{G}_n\left(\prod_{k'=1}^KU_{k' j_{k'}}^{\ast}\right)\right\}^2\right]^{1/2}\!\!,
\end{align*}
which implies
\begin{align}
\label{EQN:ed1}
   \mathbb{E}|D_1|
   &\leq \left(\sum_{k=1}^K \| a_k - b_k\|_k\right) \left[\sum_{(1, \cdots, 1)}^{(\infty, \cdots,\infty)} \mathbb{E} \left\{\prod_{k=1} ^K\left(U_{k j_{k}}^{\ast}\right)^2\right\}\right]^{1/2}
   = \left(\sum_{k=1}^K \|a_k - b_k \|_k \right)\left\{\mathbb{E}\left(\prod_{k=1}^K \|Z_{k}\|_{k}^2\right)\right\}^{1/2}.
\end{align}

Now fix $\delta >0$. Let $r = (r_1, \cdots, r_K)$ be $K$ positive integers. For $a = (a_1, \cdots, a_K) \in \mathcal{B}_{10} \times \cdots \times \mathcal{B}_{K0} \equiv \mathcal{B}_0^\ast$, let $a^{(r)} = (a_1^{(r_1)}, \cdots, a_K^{(r_K)})$, where $a_k^{(r_k)} = \sum_{j=1}^{r_k} a_{kj}\phi_{kj}$. For any $f \in \mathcal{F}_0$, we can write $f(Z) = \prod_{k=1}^K\langle a_k, Z_k \rangle_k$, for some $a \in \mathcal{B}_0^\ast$. Let $f^{(r)}(Z) = \prod_{k=1}^K \langle a_k^{(r_k)}, Z_k\rangle_k$. Then
\begin{align*}
    \sup_{f \in \mathcal{F}_0} \mathbb{G}_n \left(f - f^{(r)}\right)
    &= \sup_{a \in \mathcal{B}_0^\ast} \mathbb{G}_n\left(\prod_{k=1}^K \langle a_k, Z_k \rangle_k - \prod_{k=1}^K \langle a_k^{(r_k)}, Z_k \rangle_k\right) \\
    &= \sup_{a \in \mathcal{B}_0^\ast} \sum_{k=1}^K \sum_{L_k(r)}^{U_k(r)}
        \left\{\left(\prod_{k' = 1}^K a_{k' j_{k'}}\right)\mathbb{G}_n\left(\prod_{k' = 1}^K U_{k' j_{k'}}^\ast\right)\right\}
    \equiv D_2,
\end{align*}
where $L_k(r) = (j_1 = 1, \cdots, j_{k-1} = 1, j_k = r_k + 1, j_{k+1} = 1, \cdots, j_K = 1)$ and $U_k(r) = (j_1 = r_1, \cdots, j_{k-1}=r_{k-1}, j_k = \infty, \cdots, j_K = \infty)$.

Notice that for a single $k$ we have
\begin{align*}
   & \sup_{a \in \mathcal{B}_0^\ast} \sum_{L_k(r)}^{U_k(r)} \left(\prod_{k' = 1}^K a_{k' j_{k'}}\right)
        \mathbb{G}_n \left(\prod_{k' = 1}^K U_{k' j_{k'}}^{\ast}\right) \\
    \leq\; & \sup_{a \in \mathcal{B}_0^\ast} \left\{\sum_{L_k(r)}^{U_k(r)}\left(\prod_{k' = 1}^K a_{k' j_{k'}}^2\right)\right\}^{1/2}
        \left[\sum_{L_k(r)}^{U_k(r)}
        \left\{\mathbb{G}_n\left(\prod_{k' = 1}^K U_{k' j_{k'}}^{\ast}\right)\right\}^2\right]^{1/2} \\
    \leq\; & \left[\sum_{L_k(r)}^{U_k(r)} \left\{\mathbb{G}_n\left(\prod_{k' = 1}^K U_{k' j_{k'}}^{\ast}\right)\right\}^2\right]^{1/2}
    \equiv D_{3k}.
\end{align*}
Now,
\begin{align*}
    \mathbb{E}|D_{3k}|
    \leq \left\{\mathbb{E} \sum_{L_k(r)}^{U_k(r)}\prod_{k' = 1}^K \left(U_{k' j_{k'}}^\ast\right)^2\right\}^{1/2}
    &\leq \left\{\mathbb{E}\sum_{L_k(r)}^{(\infty, \cdots, \infty)}\prod_{k' = 1}^K \left(U_{k' j_{k'}}^\ast\right)^2\right\}^{1/2}\\
    &= \left(\mathbb{E}\left[\left\{\sum_{j=r_k+1}^\infty \left(U_{kj}^\ast\right)^2\right\}\prod_{1 \leq k' \leq K : k'\ne k}\!\!\!\!\!\!\!\!\|Z_{k'}\|_{k'}^2\right]\right)^{1/2}
    \!\!\!\!\!\!\!\!\equiv \left(\mathbb{E}C_k\right)^{1/2}.
\end{align*}
Note that for $1 \leq k \leq K$, $\mathbb{E}C_k \leq \mathbb{E} (\prod_{k' = 1}^K \|Z_{k'} \|_{k'}^2) < \infty$, and since $\sum_{j = r_k+1}^\infty( U_{kj}^\ast)^2 \xrightarrow{P} 0$ as $r_k \rightarrow \infty$, we have by the bounded convergence theorem $\mathbb{E} C_k \rightarrow 0$ as $r_k \rightarrow \infty$, which implies that $\exists\, r_k < \infty$ such that $\mathbb{E} C_k \leq \delta^2$.

Now do this for all $k$, $1 \leq k \leq K$. Thus, there exists $r = (r_1, \cdots, r_K)$ such that $r_k < \infty$ for all $1 \leq k \leq K$, and moreover, $\mathbb{E}|D_2| \leq K\delta$. Let $\mathcal{B}_0^\ast(r) = \mathcal{B}_{10}^\ast(r_1) \times \cdots \times \mathcal{B}_{K0}^\ast(r_K)$, where $\mathcal{B}_{k0}^\ast(r_k)=\{a_k^{(r_k)}=\sum_{j=1}^{r_k} a_{kj}\phi_{kj}: a_k=\sum_{j=1}^{\infty} a_{kj}\phi_{kj}\in\mathcal{B}_{k0}\}$ for $1 \leq k \leq K$. Since $\mathcal{B}_0^\ast(r)$ is compact, we have that there exists a finite subset $\mathcal{T}_\delta \subset \mathcal{B}_0^\ast(r)$ such that $\sup_{a \in \mathcal{B}_0^\ast(r)} \inf_{b \in \mathcal{T}_\delta} \sum_{k=1}^K \|a_k - b_k\|_{k} \leq \delta$, where $a = (a_1, \cdots, a_K)$ and $b = (b_1, \cdots, b_K)$. For the given set $\mathcal{T}_{\delta}$ of $K$-tuples of finite sequences, define $\mathcal{T}_{\delta}^\ast$ as:
\[
\mathcal{T}_{\delta}^\ast
\!=\!\left\{a^\ast \!=\! (a_1^\ast, \ldots, a_K^\ast) : \exists\, a
\in \mathcal{T}_{\delta} \text{ such that } a_k^\ast
\!=\! (a_{k1}, \ldots, a_{kr_k}, 0, 0, \ldots) \text{ for } k \!=\! 1, \ldots, K\right\},
\]
where $a_{kj}$ denotes the $j$th element of the finite sequence $a_k$ of length $r_k$, and $a_k^\ast$ is the corresponding infinite sequence with all elements beyond position $r_k$ set to zero.

Let $M = |\mathcal{T}_\delta^\ast| = | \mathcal{T}_\delta|$, the cardinality of $\mathcal{T}_\delta$. Let $c_1, \cdots, c_M$ be an enumeration of $\mathcal{T}_\delta^\ast$. For each $m = 1, \cdots, M$, we define $\mathcal{F}_m = \{f(Z)=\prod_{k=1}^K\langle a_k, Z_k\rangle \in \mathcal{F}_0 : \sum_{k=1}^K \| a_k^{(r_k)} - c_{m k} \| \leq \delta\}$.
Note that
\begin{align*}
    &\sup_{f, g \in \mathcal{F}_m} |\mathbb{G}_n(f-g)|
    \leq \sup_{f, g \in \mathcal{F}_m} \bigg\| \mathbb{G}_n\left(f - f^{(r)}\right) + \mathbb{G}_n\Big(f^{(r)} - \prod_{k=1}^K \langle c_{m k}, Z_k\rangle_k\Big)\\
 &\qquad\qquad\qquad\qquad\qquad\qquad\qquad\quad - \mathbb{G}_n\Big(g^{(r)}-\prod_{k=1}^K \langle c_{m k}, Z_k\rangle_k\Big) - \mathbb{G}_n\left(g-g^{(r)}\right)\bigg\|\\
    &\leq  2\left\{\sup_{f \in \mathcal{F}_0} \left\|\mathbb{G}_n\left(f-f^{(r)}\right)\right\| + \sup_{f, g \in \mathcal{F}_0 : \sum_{k=1}^K \|a_k-b_k\|_k \leq \delta} \|\mathbb{G}_n(f-g)\|\right\} \equiv D_4,
\end{align*}
where $f(Z) = \prod_{k=1}^K \langle a_k, Z_k \rangle_k$ and $g(Z) = \prod_{k=1}^K \langle b_k, Z_k \rangle_k$. Then
\[
    \mathbb{E} |D_4|
    \leq 2\left[K\delta + \delta\left\{\mathbb{E}\left(\prod_{1 \leq k \leq K}\|Z_k\|_k^2\right)\right\}^{1/2}\right],
\]
where the first term follows from $\mathbb{E}|D_2| \leq K\delta$; the second term follows from the fact that
\[
    \sup_{f \in \mathcal{F}_m} \left\|\mathbb{G}_n
    \left(f^{(r)}-\prod_{k=1}^K \langle c_{m k}, Z_k\rangle_k\right)\right\|
    \leq \sup_{f, g \in \mathcal{F}_m : \sum_{k=1}^K \|a_k-b_k\|_k \leq \delta} \left\|\mathbb{G}_n(f-g)\right\| \equiv D_6,
\]
since for all $f \in \mathcal{F}_m$, $\sum_{k=1}^K \| a_k^{(r_k)} - \sum_{j=1}^{r_k}c_{m kj}\phi_{kj}\| \leq \delta$ by construction, where $f(Z) = \prod_{k=1}^K \langle a_k, Z_k \rangle_k$ as above. Now $D_6$ is less than or equal to the second term of $D_4$ since $\mathcal{F}_m \subset \mathcal{F}_0$ for all $m = 1, \cdots, M$. We then use~(\ref{EQN:ed1}) to obtain the final form $\delta\{\mathbb{E}(\prod_{1 \leq k \leq K} \|Z_k\|_k^2)\}^{1/2}$. We make a note that $\bigcup_{m = 1}^M \mathcal{F}_m \supset \mathcal{F}$, since for each $f \in \mathcal{F}$, $f^{(r)}(Z) = \prod_{k=1}^K \langle a_k^{(r_k)}, Z_k \rangle_k$ and all such $a^{(r)} = (a_1^{(r_1)}, \cdots, a_K^{(r_k)})$ are contained in the ball $\mathcal{B}_0^\ast(r)$ by construction. Now since $D_4$ does not depend on $m$, we have that
\[
\sup_{1 \leq m \leq M} \sup_{f, g \in \mathcal{F}_m} \|\mathbb{G}_n(f-g)\| \leq D_4.
\]
Recall that $U_{kj}^\ast =\langle \phi_{kj}, Z_k \rangle_k$. Let
\[
\mathcal{G}_0 = \left\{\prod_{k=1}^K \left(\sum_{j=1}^{r_k} a_{kj} U_{kj}^\ast\right): a_k \in \mathcal{B}^{r_k} \cap \mathbb{Q}^{r_k}, 1 \leq r_k < \infty, 1 \leq k \leq K\right\},
\]
where $\mathcal{B}^{r_k} = \{x \in \mathbb{R}^{r_k} : \|x\| \leq 1\}$ and $\mathbb{Q}$ are the rationals. It is easy to verify that $\mathcal{G}_0$ is countable, and moreover, that for every $f \in \mathcal{F}_0$, there exists a sequence $\{g_n\} \in \mathcal{G}_0$ such that $g_n(Z) \rightarrow f(Z)$ pointwise on $Z$. Thus $\mathcal{F}_0$ is pointwise measurable; as it is easy to verify that the two suprema in $D_4$ are also measurable, we have
\begin{align*}
   &\limsup_{n \rightarrow \infty} P^\ast\left\{\sup_{1 \leq m \leq M} \sup_{f, g \in \mathcal{F}_m}\|\mathbb{G}_n(f-g)\| > \epsilon \right\} \\
   \leq\; &\limsup_{n \rightarrow \infty} P\left[2\left\{\sup_{f \in \mathcal{F}_0} \|\mathbb{G}_n(f-f^{(r)})\|
   + \sup_{f, g \in \mathcal{F}_0 : \sum_{k=1}^K \|a_k-b_k\|_k \leq \delta} \|\mathbb{G}_n(f-g)\|\right\} > \epsilon\right]\\
   \leq\;& 2\epsilon^{-1}\left[K\delta + \delta\left\{\mathbb{E}\left(\prod_{k=1}^K\|Z_k\|_k^2\right)\right\}^{1/2}\right].
\end{align*}
For the above to be less than $\eta$, it is sufficient for $\delta \leq \epsilon \eta(3[K + \{\mathbb{E}(\prod_{k =1}^K \|Z_k\|_k^2)\}]^{1/2})^{-1}$, since $\delta$, as well as $\epsilon$ and $\eta$ are arbitrary. We have thus satisfied~(\ref{EQN:zerostar}) as needed, and our proof is complete.
\end{proof}

\subsection{Proof of Corollary \ref{COR1:product}}

\begin{proof}[Proof of Corollary~\ref{COR1:product}.]
For any $f \in \mathcal{F}$, there exists some $h \in \mathcal{H}$ with $\|h\| \leq 1$, such that $\mathbb{G}_n f = \mathbb{G}_n \langle h, Z \rangle = \mathbb{G}_n \langle h, Z - \mu \rangle + \mathbb{G}_n \langle h, \mu \rangle$. Since $\langle h, \mu \rangle$ is not random, $\mathbb{G}_n \langle h, \mu \rangle = 0$ for all $h \in \mathcal{H}$ with $\|h\| \leq 1$. The desired conclusion now follows from Theorem~\ref{THM:product} by setting $K = 1$.
\end{proof}

\subsection{Proof of Corollary \ref{COR2:product}}

\begin{proof}[Proof of Corollary~\ref{COR2:product}.]
We have
\begin{align*}
    \prod_{k=1}^K \langle h_k, Z_k\rangle_k
    &= \prod_{k=1}^K \left(\langle h_k, Z_k - \mu_k\rangle_k + \langle h_k , \mu_k\rangle_k\right)= \sum_{k'=0}^K \sum_{r \in \mathcal{N}_{k'}} \prod_{k=1}^K \langle h_{k}, Z_{k}- \mu_{k}\rangle_{k}^{r_{k}}\langle h_{k}, \mu_{k} \rangle_{k}^{(1-r_{k})}.
\end{align*}
This implies that $\mathcal{F} \subset \sum_{k'=0}^K \sum_{r \in \mathcal{N}_{k'}} \mathcal{F}_{(r)}$, where $\mathcal{F}_{(r)} = \{c \mathcal{G}_{(r)} : |c| \leq c_\ast\}$, with $c_\ast = \prod_{k=1}^K (1 + \| \mu_k \|_k) < \infty$, and where $\mathcal{G}_{(r)} = \{\prod_{k =1}^K \langle h_{k}, Z_{k} - \mu_{k}\rangle_{k}^{r_{k}} : h_{k} \in \mathcal{B}_{k}, 1 \leq k \leq K\}$; here, the addition between classes is the Minkowski addition. By Theorem~\ref{THM:product}, each $\mathcal{G}_{(r)}$ is Donsker, and it is easy to see that $\mathcal{F}_{(r)}$ is also Donsker. Since finite sums of Donsker classes are also Donsker, we have that $\mathcal{F}$ is Donsker.
\end{proof}

\subsection{Proof of Theorem \ref{THM:op}}

For random element $Z\in\mathcal{H}$, define the centered variable $Z^c = Z- \mu$, and the sample-centered variable $Z_i^{sc} = Z_i- \bar{Z}_n$.

\begin{proof}[Proof of Theorem~\ref{THM:op}.]
From~(A1), we know that $\mathcal{F}_1\equiv\{\langle a,Z\rangle:\; a\in\mathcal{B}\}$ is Donsker. Thus $\mathcal{F}_1$ is also Glivenko-Cantelli, and hence $\|\bar{Z}_n-\mu\|\xrightarrow{as\ast}0$. Since also
\[
    \widehat{V}_n=G_N\mathbb{P}_n\left(Z^c{Z^c}^{\top}\right)G_N-G_N\left\{\left(\bar{Z}_n-\mu\right)\left(\bar{Z}_n-\mu\right)^{\top}\right\}G_N,
\]
we only need to show that $\sup_{a_1,a_2\in\mathcal{B}}|a_1^{\top}(\widehat{V}_n^{\ast}-V_0)a_2|\xrightarrow{as\ast}0$, where $\widehat{V}_n^{\ast}=\mathbb{P}_n G_N(Z^c{Z^c}^{\top})G_N$. Since $\|Z^c\|$ is an envelope for $\mathcal{F}_1$, and since also $\mathbb{E}\|Z^c\|^2<\infty$, we have that $\mathcal{F}_1\cdot\mathcal{F}_1$ is Glivenko-Cantelli. Since also $G_N a\subset \mathcal{B}$, for any $a\in \mathcal{B}$ by the properties of projections, we therefore conclude that
$C_n\equiv\sup_{a_1,a_2\in\mathcal{B}}|a_1^{\top}\{G_N(\mathbb{P}_n-P) (Z^c{Z^c}^{\top}) G_N\}a_2|\xrightarrow{as\ast}0.$
Note that this implies that $C_n^{\ast}\xrightarrow{as\ast}0$ for some measurable majorant $C_n^{\ast}$ of $C_n$. Now fix $\epsilon>0$. We then have that
\begin{align*}
    &P\left(\limsup_{n\rightarrow\infty} \sup_{a_1, a_2 \in \mathcal{B}} \left| a_1^\top (\widehat{V}_n^{\ast} -V_0)a_2 \right|>\epsilon\right)\\
    &\leq P\left(\limsup_{n\rightarrow\infty}C_n^{\ast}>\epsilon/2\right)+\mathbbm{1}\!\left[\limsup_{n\rightarrow\infty}\sup_{a_1, a_2 \in \mathcal{B}} \left| a_1^\top \left\{G_N P\left(Z^cZ^{c\top}\right)G_N-V_0\right\}a_2 \right|>\epsilon/2\right]\\
    &\leq 0+\mathbbm{1}\!\left\{\limsup_{n\rightarrow\infty}\mathbb{E}\left(\|G_NZ^c-Z^c\|\|G_NZ^c\|+\|Z^c\|\|G_NZ^c-Z^c\|\right)>\epsilon/2\right\}\\
    &\leq \mathbbm{1}\!\left\{2\left(\mathbb{E}\|Z\|^2\right)^{1/2}\limsup_{n\rightarrow\infty}\left(\mathbb{E}\|G_NZ^c-Z^c\|^2\right)^{1/2}>\epsilon/2\right\}
    =0,
\end{align*}
by Assumption~(A2). Thus Part~(a) follows.

To prove Part~(b), let $a_1, a_2 \in \mathcal{B}$. Now
\begin{align*}
    D_1\equiv\;& a_1^\top \left\{\sqrt{n} \left(\widehat{V}_n - V_0\right) - \mathbb{G}_n\left(Z^c {Z^c}^\top\right)\right\}a_2 \\
    =\;& a_1^\top \left[\sqrt{n}\mathrm{G}_N \mathbb{P}_n\left\{(Z- \bar{Z}_n)(Z-\bar{Z}_n)^\top \right\}\mathrm{G}_N  - \sqrt{n}P(Z^c {Z^c}^\top) - \mathbb{G}_n(Z^c {Z^c}^\top)\right]a_2 \\
    =\;& a_1^\top\!\!\left[\!\mathrm{G}_N\!\!\left\{\!\mathbb{G}_n\!(Z^c {Z^c}^\top\!\!) \!-\! \sqrt{n}(\bar{Z}_n \!\!-\!\! \mu)(\bar{Z}_n \!\!-\!\! \mu )^\top\!\!\right\}\!\mathrm{G}_N
    \!+\! \sqrt{n}P(\mathrm{G}_N\! Z^c {Z^c}^\top \! \mathrm{G}_N \!\!-\!\! Z^c {Z^c}^\top \!) \!-\! \mathbb{G}_n\!(Z^c {Z^c}^\top\!)\!\right]\! a_2.
\end{align*}
Now $\sup_{a_1, a_2 \in \mathcal{B}} |a_1^\top \mathrm{G}_N\sqrt{n}(\bar{Z}_n - \mu)(\bar{Z}_n - \mu)^\top \mathrm{G}_N a_2| = o_P(1)$, since $\mathcal{F}_1$ is Donsker and also Glivenko-Cantelli as verified above. Also, we have
\begin{align*}
    &\sqrt{n}\! P \!\left| a_1^\top\!\left(\!\mathrm{G}_N Z^c {Z^c}^\top \mathrm{G}_N \!-\! Z^c {Z^c}^\top\!\right)\!a_2 \right|
    \leq \sqrt{n}P\left(2 \|\mathrm{G}_N Z^c - Z^c\|\| Z^c \|\right)
    \leq 2\sqrt{n}\left(P\|\mathrm{G}_N Z^c - Z^c\|^2\right)^{1/2}\!\!\left(P\|Z^c \|^2\right)^{1/2}
    = o(1).
\end{align*}
Thus, $D_1= a_1^\top\{\mathrm{G}_N \mathbb{G}_n(Z^c {Z^c}^\top)\mathrm{G}_N - \mathbb{G}_n(Z^c {Z^c}^\top)\}a_2 + o_P(1)$, where the $o_P(1)$ is uniform over $a_1, a_2 \in \mathcal{B}$. Now
\[
    a_1^\top \left\{\mathrm{G}_N \mathbb{G}_n \! \left(Z^c {Z^c}^\top\right)\mathrm{G}_N \!-\! \mathbb{G}_n\left(Z^c {Z^c}^\top\right)\right\}a_2
    \!=\! \mathbb{G}_n\left(\!\langle a_1, \mathrm{G}_N Z^c \rangle \langle \mathrm{G}_N Z^c, a_2 \rangle \!-\! \langle a_1, Z^c\rangle \langle Z^c, a_2 \rangle\!\right),
\]
and
\begin{align*}
    D_2 \equiv\;&
    \mathbb{E} \left|\mathbb{G}_n\left(\langle a_1, \mathrm{G}_N Z^c \rangle \langle \mathrm{G}_N Z^c, a_2 \rangle - \langle a_1, Z^c \rangle\langle Z^c, a_2 \rangle\right)\right| \\
    =\;& \mathbb{E}\left|\mathbb{G}_n\left(\langle a_1, \mathrm{G}_N Z^c - Z^c \rangle \langle \mathrm{G}_N Z^c, a_2 \rangle + \langle a_1, Z^c \rangle \langle\mathrm{G}_N Z^c - Z^c, a_2 \rangle\right)\right|\\
    \leq\;& \sqrt{2} \left(\mathbb{E}\left|\langle a_1, \mathrm{G}_N Z^c - Z^c\rangle^2 \langle \mathrm{G}_N Z^c, a_2\rangle^2 + \langle a_1, Z^c \rangle ^2\langle \mathrm{G}_N Z^c - Z^c, a_2 \rangle^2 \right|\right)^{1/2}\\
    \leq\;& 2 \left\{\mathbb{E}\left(\|\mathrm{G}_N Z^c - Z^c\|^2\|Z^c\|^2\right)\right\}^{1/2},
\end{align*}
which does not depend on $a_1, a_2 \in \mathcal{B}$. We have $D_2 \leq 4(\mathbb{E}\|Z^c \|^4)^{1/2} < \infty$ by assumption. Moreover, $\mathbb{E}\|\mathrm{G}_N Z^c - Z^c\|^2 = o(n^{-1})$ by Assumption~(A2). This implies that $\|\mathrm{G}_N Z^c - Z^c\|^2 \xrightarrow{P} 0$, and thus, by the dominated convergence theorem, $D_2 \rightarrow 0$, and part~(i) is proved. Part~(ii) follows from Theorem~\ref{THM:product}, and part~(iii) is a direct result of parts~(i) and~(ii).
\end{proof}

\subsection{Proof of Theorem \ref{THM:maxes}}

By the spectral theorem, $V_0= \sum_{j=1}^\infty \lambda_{j}\phi_j\phi_j^\top$. As a consequence of Part~(a) of Theorem~\ref{THM:op}, $E_n \equiv \sup_{a_1, a_2 \in \mathcal{B}} |a_1^\top (\widehat{V}_n - V_0) a_2| \xrightarrow{as\ast} 0$. From the definitions of $\widehat{\lambda}_1$ and $\lambda_1$, $\widehat{\phi}_1$ and $\phi_1$, and $\widehat{V}_n$ and $V_0$, we have
\[
    \widehat{\lambda}_1 \!=\! \widehat{\phi}_1^\top \widehat{V}_n \widehat{\phi}_1 \!\geq\! \phi_1^\top \widehat{V}_n \phi_1 \!\geq\! \phi_1^\top V_0\phi_1 \!-\!E_n\!=\! \lambda_1 \!-\! E_n \!\geq\! \widehat{\phi}_1^\top V_0\widehat{\phi}_1 \!-\! E_n \!\geq\! \widehat{\phi}_1^\top\widehat{V}_n\widehat{\phi}_1 \!-\! 2E_n\!=\! \widehat{\lambda}_1 \!-\! 2E_n,
\]
which implies that $|\widehat{\lambda}_1 - \lambda_1| \leq E_n$, and thus $|\widehat{\lambda}_1-\lambda_1| \xrightarrow{as\ast}0$. Now, $0 \geq \widehat{\phi}_1^\top V_0 \widehat{\phi}_1 - \phi_1^\top V_0 \phi_1 \geq \widehat{\phi}_1^\top \widehat{V}_n\widehat{\phi}_1 - \phi_1^\top \widehat{V}_n\phi_1 - 2E_n \geq -2E_n$, and it follows that
\begin{equation}\label{eq2}
     \sum_{j=1}^\infty \lambda_j \langle \widehat{\phi}_1, \phi_j\rangle^2 - \lambda_1 \xrightarrow{as\ast} 0.
\end{equation}
Let $\rho_n = 1 - \sum_{j=1}^\infty \langle\widehat{\phi}_1,\phi_j\rangle^2$, and note that $\rho_n \geq 0$, almost surely. Then~(\ref{eq2}) implies that $\sum_{j=1}^\infty \lambda_j\langle \widehat{\phi}_1, \phi_j\rangle^2 - \lambda_1(\rho_n + \sum_{j=1}^\infty\langle \widehat{\phi}_1, \phi_j \rangle^2) \xrightarrow{as\ast} 0$; equivalently, $- \sum_{j=1}^\infty (\lambda_1 - \lambda_j)\langle \widehat{\phi}_1, \phi_j \rangle^2 - \rho_n\lambda_1 \xrightarrow{as\ast} 0$. Therefore, both $\sum_{j=2}^{\infty}(\lambda_1-\lambda_j)\langle \widehat{\phi}_1, \phi_j\rangle^2\xrightarrow{as\ast} 0$ and $\rho_n \xrightarrow{as\ast} 0$. As a result, we have $\sum_{j\geq 2}\langle \widehat{\phi}_1, \phi_j\rangle^2\xrightarrow{as\ast} 0$, implying that $\langle \widehat{\phi}_1, \phi_{1} \rangle^2 \xrightarrow{as\ast} 1$. Thus either $\widehat{\phi}_1 \xrightarrow{as\ast} \phi_1$ or $\widehat{\phi}_1 \xrightarrow{as\ast} - \phi_1$. Without loss of generality, we can change the sign of $\phi_1$, since either sign is correct; hence $\widehat{\phi}_1 \xrightarrow{as\ast} \phi_1$.

Now, suppose for some $1 \leq j < m$, $\max_{1 \leq j' \leq j} \{\|\widehat{\phi}_{j'} - \phi_{j'} \| \vee|\widehat{\lambda}_{j'} - \lambda_{j'}|\} \xrightarrow{as\ast} 0$. Let $\widehat{V}_n(j) = \widehat{V}_n - \sum_{{j'} = 1}^j \widehat{\lambda}_{j'} \widehat{\phi}_{j'} \widehat{\phi}_{j'}^\top$ and $V_0(j) = V_0 - \sum_{j' =1}^j \lambda_{j'} \phi_{j'} \phi_{j'}^\top$. Now $E_n(j) = \sup_{a_1, a_2 \in \mathcal{B}} | a_1^\top \{\widehat{V}_n(j) - V_0(j)\} a_2 | \xrightarrow{as\ast} 0$, by the supposition and previous results. Moreover,
\begin{align*}
    \widehat{\lambda}_{j+1} &= \widehat{\phi}_{j+1}^\top \widehat{V}_n(j) \widehat{\phi}_{j+1} \geq \phi_{j+1}^\top \widehat{V}_n(j) \phi_{j+1} \geq \phi_{j+1}^\top V_0(j)\phi_{j+1} - E_n(j) = \lambda_{j+1} - E_n(j)\\
    &\geq \widehat{\phi}_{j+1}^\top V_0(j)\widehat{\phi}_{j+1} - E_n(j)\geq \widehat{\phi}_{j+1}^\top \widehat{V}_n(j)\widehat{\phi}_{j+1} - 2E_n(j)= \widehat{\lambda}_{j+1}-2E_n(j),
\end{align*}
which implies that $| \widehat{\lambda}_{j+1} - \lambda_{j+1} | \xrightarrow{as\ast} 0$. Now,
\begin{align*}
    0 \!\geq\! \widehat{\phi}_{j+1}^\top \!V_0(j)\!\widehat{\phi}_{j+1} \!-\! \phi_{j+1}^\top\! V_0(j)\! \phi_{j+1}\!\geq\! \widehat{\phi}_{j+1}^\top\! \widehat{V}_n(j)\!\widehat{\phi}_{j+1} \!-\! \phi_{j+1}^\top\! \widehat{V}_n(j) \!\phi_{j+1} \!-\! 2E_n(j) \!\geq\! -\!2E_n(j),
\end{align*}
which implies that $\sum_{j' = j+1}^\infty \lambda_{j'}\langle\widehat{\phi}_{j+1}, \phi_{j'}\rangle^2 - \lambda_{j+1} \xrightarrow{as\ast} 0$. Note that
\[\rho_n(j) \equiv 1 - \sum_{j' = j+1}^\infty \langle\widehat{\phi}_{j+1}, \phi_{j'}\rangle^2 \geq 0.\]
Consequently, $\sum_{j' = j+1}^\infty \lambda_{j'} \langle\widehat{\phi}_{j+1}, \phi_{j'}\rangle^2 - \lambda_{j+1}\{\sum_{j' = j+1}^\infty \langle\widehat{\phi}_{j+1}, \phi_{j'}\rangle^2 + \rho_n(j)\} \xrightarrow{as\ast} 0$; it follows that $- \sum_{j' = j+1}^\infty (\lambda_{j+1}- \lambda_{j'})\langle\widehat{\phi}_{j+1}, \phi_{j'}\rangle^2 \xrightarrow{as\ast} 0$ and $\rho_n(j)\xrightarrow{as\ast} 0$. This implies that $\sum_{j' = j+2}^\infty \langle\widehat{\phi}_{j+1}, \phi_{j'}\rangle^2 \xrightarrow{as\ast} 0$, and as a result we have $\{\langle\widehat{\phi}_{j+1}, \phi_{j+1}\rangle^2 -1\}\xrightarrow{as\ast} 0$, which implies that $\langle\widehat{\phi}_{j+1}, \phi_{j+1}\rangle^2 \xrightarrow{as\ast}1$. Provided we switch the sign of $\phi_{j+1}$ if needed, we obtain that $\| \widehat{\phi}_{j+1} - \phi_{j+1} \| \xrightarrow{as\ast} 0$. Repeat this process until $j+1 = m$, and we have $\widehat{V}_n(m)$ and $V_0(m)$ so that $E_n(m) \xrightarrow{as\ast} 0$. Then
\begin{align*}
    \widehat{\lambda}_{m+1}
    &= \sup_{\phi \in \mathcal{B}}\phi^\top \widehat{V}_n(m) \phi
    \geq \sup_{\phi \in \mathcal{B}} \phi^\top V_0(m) \phi - E_n(m)
    = \lambda_{m+1} - E_n(m)\\
    &\geq \sup_{\phi \in \mathcal{B}} \phi^\top \widehat{V}_n(m) \phi - 2E_n(m)
    = \widehat{\lambda}_{m+1}-2E_n(m),
\end{align*}
which implies $\widehat{\lambda}_{m+1} \xrightarrow{as\ast} \lambda_{m+1}$, even if $\lambda_{m+1} = 0$.

\subsection{Proof of Theorem \ref{THM:asymp}}

By Theorem~\ref{THM:maxes}, we have $\max_{1 \leq j \leq m} |\widehat{\lambda}_j - \lambda_j |\vee \|\widehat{\phi}_j - \phi_j\| \xrightarrow{as\ast} 0$. By the eigenvalue structure, for any $1 \leq j \leq m$, it follows that $\widehat{V}_n\widehat{\phi}_j = \widehat{\phi}_j\widehat{\lambda}_j$ and $V_0\phi_j = \phi_j \lambda_j$. Then
    \begin{align}\label{EQN:Vphi}
        (\widehat{V}_n-V_0)\widehat{\phi}_j + V_0(\widehat{\phi}_j - \phi_j)
        = (\widehat{\phi}_j-\phi_j)\widehat{\lambda}_j + \phi_j(\widehat{\lambda}_j-\lambda_j),
    \end{align}
    which indicates that $\phi_j^\top (\widehat{V}_n - V_0)\widehat{\phi}_j + \phi_j^\top V_0(\widehat{\phi}_j - \phi_j) = \phi_j^\top(\widehat{\phi}_j - \phi_j)\widehat{\lambda}_j + \phi_j^\top\phi_j(\widehat{\lambda}_j - \lambda_j)$, and consequently, $\phi_j^\top (\widehat{V}_n - V_0)\phi_j + \phi_j^\top(\widehat{V}_n-V_0)(\widehat\phi_j - \phi_j) = \phi_j^\top(\widehat\phi_j-\phi_j)(\widehat{\lambda}_j -\lambda_j)+(\widehat\lambda_j -\lambda_j)$.
    Therefore, $\widehat\lambda_j - \lambda_j = \phi_j^\top(\widehat{V}_n - V_0)\phi_j + o_P(n^{-1/2}) + o_P(|\widehat{\lambda}_j - \lambda_j|)$, and we have
    \[
        \sqrt{n}(\widehat\lambda_j - \lambda_j) = \sqrt{n}\phi_j^\top(\widehat{V}_n - V_0)\phi_j + o_P(1) = \phi_j^\top\mathbb{G}_n(Z^c {Z^c}^\top)\phi_j + o_P(1),
    \]
where the last equality follows from Theorem~\ref{THM:op}. It is also the case that~(\ref{EQN:Vphi}) implies
\begin{align*}
   &V_0(\widehat\phi_j - \phi_j) = (\widehat\phi_j-\phi_j)\lambda_j - (\widehat{V}_n - V_0)\phi_j -(\widehat{V}_n - V_0)(\widehat\phi_j-\phi_j)+ (\widehat\phi_j - \phi_j)(\widehat\lambda_j - \lambda_j) + \phi_j(\widehat\lambda_j - \lambda_j)\\
    &\Rightarrow (V_0 - \lambda_j \mathrm{I})(\widehat\phi_j - \phi_j) = -(\widehat{V}_n-V_0)\phi_j + o_P(n^{-1/2}) + \phi_j\phi_j^\top (\widehat{V}_n - V_0)\phi_j+o_P(n^{-1/2})\\
    &\Rightarrow (V_0 - \lambda_j \mathrm{I} + \phi_j\phi_j^\top)(
    \widehat\phi_j - \phi_j) = -(\widehat{V}_n - V_0)\phi_j +\phi_j\phi_j^\top(\widehat{V}_n-V_0)\phi_j + \phi_j\phi_j^\top(\widehat{\phi_j}-\phi_j)+o_P(n^{-1/2}).
\end{align*}
Since $\phi_j\phi_j^\top(\widehat{\phi}_j - \phi_j) = O_P(\|\widehat{\phi}_j - \phi_j\|^2 ) = o_P(\|\widehat{\phi}_j - \phi_j\|)$, this implies that
\begin{align*}
    &\widehat\phi_j - \phi_j = -\left(V_0 - \lambda_j \mathrm{I} + \phi_j \phi_j^\top\right)^{-1}\!\Big\{\left(\widehat{V}_n -V_0\right)\phi_j - \phi_j\phi_j^\top\left(\widehat{V}_n-V_0\right)\phi_j +o_P\left(n^{-1/2}+\|\widehat\phi_j -\phi_j\|\right)\!\Big\}\\
    &= \left\{\sum_{j' \neq j, j' \geq 1}\left(\lambda_j - \lambda_{j'}\right)^{-1}\phi_{j'}\phi_{j'}^\top\right\}\left\{\left(\widehat{V}_n-V_0\right)\phi_j + o_P\left(n^{-1/2}+\|\widehat\phi_j-\phi_j\|\right)\right\}.
\end{align*}
We then have $\sqrt{n}\|\widehat\phi_j - \phi_j\| = O_P(1) + o_P(\sqrt{n}\|\widehat{\phi}_j-\phi_j\|)$, and $\sqrt{n}\|\widehat\phi_j - \phi_j \| = O_P(1)$. Thus
\[
\sqrt{n}(\widehat\phi_j - \phi_j) =\left(\sum_{j' \neq j, j' \geq 1} (\lambda_j - \lambda_{j'})^{-1}\phi_{j'}\phi_{j'}^\top\right)\sqrt{n}(\widehat{V}_n - V_0)\phi_j + o_P(1),
\]
and the desired result follows from Theorem~\ref{THM:op}.

\subsection{Proof of Theorem \ref{THM:delta}}

We now describe an approach to evaluating whether $\delta_n=o(n^{-1})$ from Assumption~(A2). Let $\widehat{\delta}_n = n^{-1} \sum_{i=1}^n \| Z_i - \bar{Z}_n \|^2 - \sum_{\ell = 1}^N n^{-1}\sum_{i=1}^n \langle \psi_\ell, Z_i - \bar{Z}_n\rangle^2$, let $Z_i^c = Z_i - \mu$, let $\bar{Z}_n^c = n^{-1}\sum_{i=1}^n Z_i^c$, and let $\rho_{Nj}^2 = \sum_{\ell = 1}^N \langle \psi_\ell, \phi_j \rangle^2$. Recall that we can write $Z_i^c = \sum_{j=1}^\infty \lambda_{j}^{1/2}U_{ji}\phi_j$; we will denote $\bar{U}_{jn} = n^{-1}\sum_{i=1}^n U_{ji}$. Note that $\widehat{\delta}_n$ can be computed from the data and, moreover,
\begin{align*}
    &\widehat{\delta}_n = n^{-1}\sum_{i=1}^n\|Z_i^c - \bar{Z}_n^c\|^2 - \sum_{\ell = 1}^N n^{-1}\sum_{i=1}^n \langle \psi_\ell, Z_i^c - \bar{Z}_n^c\rangle^2\\
    &= \sum_{j=1}^\infty \lambda_jn^{-1}\sum_{i=1}^n\left(U_{ji} -  \bar{U}_{jn}\right)^2 - \sum_{\ell = 1}^N \psi_\ell^\top \sum_{j=1}^\infty \lambda_j n^{-1}\sum_{i=1}^n\left(U_{ji} - \bar{U}_{jn}\right)^2\phi_j\phi_j^\top\psi_j\\
    &=\sum_{j=1}^\infty \lambda_j n^{-1} \sum_{i=1}^n \left(U_{ji} - \bar{U}_{jn}\right)^2 - \sum_{j=1}^\infty \lambda_j n^{-1}\sum_{i=1}^n \left(U_{ji}-\bar{U}_{jn}\right)^2 \rho_{Nj}^2\\
    &= \sum_{j=1}^{\infty} \lambda_j n^{-1} \sum_{i=1}^n \left(U_{ji}-\bar{U}_{jn}\right)^2\left(1-\rho_{Nj}^2\right)= n^{-1}\sum_{i=1}^n \|Z_i^c - \mathrm{G}_NZ_i^c - \bar{Z}_n^c + \mathrm{G}_N\bar{Z}_n^c\|^2 \\
    &= \delta_n^\ast - \|\bar{Z}_n^c - \mathrm{G}_N\bar{Z}_n^c\|^2,
\end{align*}
where $\delta_n^\ast = n^{-1}\sum_{i=1}^n \|Z_i^c - \mathrm{G}_NZ_i^c\|^2$, and thus $\widehat{\delta}_n \leq \delta_n^\ast$. Let $\delta_n = \mathbb{E}\delta_n^\ast$, let $\widehat{S}_n^2 = n^{-1}\sum_{i=1}^n (\|Z_i^c - \mathrm{G}_NZ_i^c - \bar{Z}_n^c + \mathrm{G}_N\bar{Z}_n^c\|^2 - \widehat{\delta}_n)^2$ (which also can be computed from the data), and let $\widehat{T}_n = \sqrt{n}\widehat{\delta}_n/\sqrt{\widehat{S}_n^2+1/n}$. We will denote $\sigma_0^2 \equiv \text{Var}(\|Z^c-\mathrm{G}_NZ^c\|^2)$. In the following theorem, the given moment condition implies that $\mathbb{E}\|Z^c\|^4<\infty$. The proof requires this slightly stronger condition.

\begin{proof}[Proof of Theorem~\ref{THM:delta}.]
    For part~(i), we have
    \[
        \widehat{T}_n = \frac{\sqrt{n}\widehat{\delta}_n}{\sqrt{\widehat{S}_n^2+ \frac{1}{n}}} \leq \frac{\sqrt{n}
        \delta_n^\ast}{\sqrt{\widehat{S}_n^2+ \frac{1}{n}}} = \frac{n\delta_n^\ast}{\sqrt{1 + n\widehat{S}_n^2}} \leq n\delta_n^\ast,
    \]
    and since $\mathbb{E}(n\delta_n^\ast) = n\delta_n \rightarrow 0$, it follows that $\widehat{T}_n \xrightarrow{P} 0$.

    For part~(ii), we have
    \begin{align*}
        \widehat{T}_n &\leq \frac{\sqrt{n}\delta_n^\ast}{\sqrt{\widehat{S}_n^2+ \frac{1}{n}}} = \frac{\sqrt{n}\delta_n^\ast}{\sqrt{\sigma_0^2 + \frac{1}{n}}}\left(\frac{\widehat{S}_n^2+ \frac{1}{n}}{\sigma_0^2+ \frac{1}{n}}\right)^{-1/2}
        = \left(\frac{\widehat{S}_n^2 + \frac{1}{n}}{\sigma_0^2+\frac{1}{n}}\right)^{-1/2}\left(\frac{\sqrt{n}(\delta_n^\ast - \delta_n)}{\sqrt{\sigma_0^2+ \frac{1}{n}}} + \frac{\sqrt{n}\delta_n}{\sqrt{\sigma_0^2+ \frac{1}{n}}}\right)\\
        &\leq \left(\frac{\widehat{S}_n^2+ \frac{1}{n}}{\sigma_0^2+ \frac{1}{n}}\right)^{-1/2}\left(\frac{-\sqrt{n}|\delta_n^\ast - \delta_n|}{\sigma_0}+\frac{\sqrt{n}\delta_n}{\sqrt{\sigma_0^2+ \frac{1}{n}}}\right).
    \end{align*}
    Fix $\eta > 0$, then
    \begin{align}\label{eqn:prob}
        \Pr\!\left\{\!\left(\frac{\widehat{S}_n^2\!+\!\frac{1}{n}}{\sigma_0^2\!+\!\frac{1}{n}}\right)^{-1/2}\!\! < \eta\!\right\} &= \Pr\!\left\{\!\left(\frac{\widehat{S}_n^2 + \frac{1}{n}}{\sigma_0^2+\frac{1}{n}}\right)\! >\! \eta^{-2}\right\}\!\leq\! \eta^2\mathbb{E}\!\left(\frac{\widehat{S}_n^2+\frac{1}{n}}{\sigma_0^2 + \frac{1}{n}}\right) \\
        &\leq \eta^2 \left\{\frac{\sigma_0^2\left(1-\frac{1}{n}\right) + \frac{8\mathbb{E}\|Z^c\|^4}{n^2}+\frac{1}{n}}{\sigma_0^2 + \frac{1}{n}}\right\}
        \leq \eta^2(1+O(n^{-1})), \nonumber
    \end{align}
    for all $n$ large enough. We will verify that $\mathbb{E}\widehat{S}_n^2\leq \sigma_0^2(1-1/n)+8\mathbb{E}\|Z^c\|^4/n^2$ soon. For now, fix $c_1 < \infty$; then $\Pr(-\sqrt{n}|\delta_n^\ast - \delta_n|/\sigma_0 < -c_1) \leq 1/c_1^2$. Now,
    \begin{align*}
        \sigma_0^2 &= \text{Var}\left(\|Z^c- \mathrm{G}_NZ^c\|^2\right) = \text{Var}\left\{\sum_{j=1}^\infty \lambda_j U_j^2\left(1 - \rho_{Nj}^2\right)\right\}\\
        &\leq \left\{\sum_{j=1}^\infty \lambda_j \left(\mathbb{E}U_j^4\right)^{1/2}\left(1-\rho_{Nj}^2\right)\right\}^2 \leq \delta_n^2 \sup_{j \geq 1} \mathbb{E}U_j^4 \leq \delta_n^2c_2,
    \end{align*}
    where $c_2 = \sup_{j \geq 1}\mathbb{E}U_j^4 < \infty$. We then have $\sqrt{n} \delta_n/\sqrt{\sigma_0^2+1/n} = n\delta_n/\sqrt{1 + n\sigma_0^2} \geq n\delta_n/\sqrt{1+ c_2n\delta_n^2} = \{1/(n^2\delta_n^2) + c_2/n\}^{-1/2}$, where $\liminf_{n \rightarrow \infty}\{1/(n^2\delta_n^2)+c_2/n\}^{-1/2} \geq k_1$, provided $\liminf_{n \rightarrow \infty} n\delta_n \geq k_1$. Now, fix $\epsilon > 0$. Let $1/c_1^2 = \epsilon/2$, or alternatively, $c_1 = \sqrt{2/\epsilon}$. This implies that $\Pr(-\sqrt{n}|\delta_n^\ast - \delta_n|/\sigma_0 > -c_1) \geq 1- \epsilon/2$. Let $\eta^2 = \epsilon/2$; then $\eta = \sqrt{\epsilon/2}$ and this implies that $\liminf_{n\rightarrow\infty}\Pr\{(\widehat{S}_n^2 + 1/n)^{-1/2}(\sigma_0^2 + 1/n)^{1/2})> \eta\} \geq 1- \epsilon/2$. We can now pick $k_1$ such that $\sqrt{\epsilon/2}\left(-\sqrt{2/\epsilon}+k_1\right) \geq k$ and thus $k_1 \geq k\sqrt{2/\epsilon}+\sqrt{2/\epsilon}$. As long as~(\ref{eqn:prob}) holds, this implies that
    \begin{align*}
        \liminf_{n \rightarrow \infty} \Pr\!\left(\!\widehat{T}_n \!\geq\! k\!\right) \!\geq\! \liminf_{n \rightarrow \infty} \Pr\!\left\{\!\frac{-\sqrt{n}|\delta_n^\ast\!-\!\delta_n|}{\sigma_0} \!>\! -c_1,\! \left(\!\frac{\widehat{S}_n^2\!+\!\frac{1}{n}}{\sigma_0^2\!+\!\frac{1}{n}}\!\right)^{-1/2} \!>\! \eta \!\right\}
        \!\geq\! 1 \!-\!\frac{2\epsilon}{2} \!=\! 1\!-\! \epsilon.
    \end{align*}
    To complete the proof, we need to verify that $\mathbb{E}\widehat{S}_n^2 \leq \sigma_0^2(1-1/n)+ 8\mathbb{E}\|Z^c\|^4/n^2$. To that end, we have
    \begin{align*}
        &\widehat{S}_n^2 = \frac{1}{n}\sum_{i=1}^n \left(\|Z_i^c - \mathrm{G}_NZ_i^c -\bar{Z}_n^c + \mathrm{G}_N\bar{Z}_n^c\|^2 - \widehat{\delta}_n\right)^2\\
        &= \frac{1}{n} \sum_{i=1}^n \Big\{\left(\|Z_i^c -\mathrm{G}_N Z_i^c\|^2 - \delta_n^\ast \right)^2 + 4\langle Z_i^c - \mathrm{G}_N Z_i^c, \bar{Z}_n^c -\mathrm{G}_N\bar{Z}_n^c \rangle^2\\
        &\qquad+ 4\|\bar{Z}_n^c - \mathrm{G}_N \bar{Z}_n^c\|^4 - 4\|Z_i^c -\mathrm{G}_N Z_i^c\|^2\langle Z_i^c - \mathrm{G}_NZ_i^c, \bar{Z}_n^c - \mathrm{G}_N\bar{Z}_n^c\rangle \\
        &\qquad+ 4\delta_n^\ast \|\bar{Z}_n^c - \mathrm{G}_N\bar{Z}_n^c\|^2 - 8 \|\bar{Z}_n^c - \mathrm{G}_N \bar{Z}_n^c\|^4\Big\}\\
        &= \frac{1}{n}\sum_{i=1}^n \left(\|Z_i^c - \mathrm{G}_NZ_i^c\|^2 - \delta_n^\ast\right)^2 + \frac{4}{n}\sum_{i=1}^n \langle Z_i^c - \mathrm{G}_NZ_i^c, \bar{Z}_n^c - \mathrm{G}_N\bar{Z}_n^c\rangle^2 \\
        &\qquad-\frac{4}{n}\sum_{i=1}^n \|Z_i^c -\mathrm{G}_NZ_i^c\|^2\langle Z_i^c - \mathrm{G}_NZ_i^c, \bar{Z}_n^c - \mathrm{G}_N\bar{Z}_n^c\rangle + 4\delta_n^\ast \|\bar{Z}_n^c - \mathrm{G}_N\bar{Z}_n^c \|^2 - 4\|\bar{Z}_n^c - \mathrm{G}_N \bar{Z}_n^c\|^4.
    \end{align*}
    We have $\mathbb{E}(\|Z_i^c - \mathrm{G}_NZ_i^c\|^2 - \delta_n^\ast)^2 = \sigma_0^2(1-1/n)$ and
    \begin{align*}
        \mathbb{E}\langle Z_i^c - \mathrm{G}_NZ_i^c, \bar{Z}_n^c - \mathrm{G}_N\bar{Z}_n^c\rangle^2 &= \mathbb{E}\left\{\sum_{j=1}^\infty \lambda_j U_{ji}(\bar{U}_{jn})(1-\rho_{Nj}^2)\right\}^2\\
        &= \sum_{j=1}^\infty\sum_{j'=1}^\infty \lambda_j\lambda_{j'}\mathbb{E}\left(U_{ji}\bar{U}_{jn}U_{j'i}\bar{U}_{j'n}\right)(1-\rho_{Nj}^2)(1-\rho_{Nj'}^2),
    \end{align*}
    where
    \begin{align*}
        \mathbb{E}\left(U_{ji}\bar{U}_{jn}U_{j'i}\bar{U}_{j'n}\right)&=n^{-2}\mathbb{E}\left\{U_{j1}U_{j'1}\!\left(\sum_{i=1}^n U_{ji}\right)\!\left(\sum_{i=1}^n U_{j'i}\right)\right\}
        = n^{-2}\mathbb{E}\!\left(U_{j1}^2U_{j'1}^2\right),
    \end{align*}
    which implies that $\mathbb{E}\langle Z_i^c - \mathrm{G}_NZ_i^c, \bar{Z}_n^c - \mathrm{G}_N\bar{Z}_n^c\rangle^2 = n^{-2}\mathbb{E}\|Z^c\|^4$. Further,
    \begin{align*}
        \mathbb{E}\left(\|Z_i - \mathrm{G}_NZ_i \|^2\langle Z_i - \mathrm{G}_NZ_i, \bar{Z}_n^c - \mathrm{G}_N\bar{Z}_n^c\rangle\right)
        = \frac{\mathbb{E}\|Z^c\|^4}{n}.
    \end{align*}
    Additionally, we have $\mathbb{E}\left(\delta_n^\ast\|\bar{Z}_n^c - \mathrm{G}_N\bar{Z}_n^c\|^2\right) \leq\mathbb{E}\|Z^c\|^4/n^2 + \delta_n^2/n \leq \mathbb{E}\|Z^c\|^4/n^2 + \mathbb{E}\|Z^c\|^4/n$. We can ignore the $-4\mathbb{E}\|\bar{Z}_n^c-\mathrm{G}_N\bar{Z}_n^c\|^4$ term since it is $\leq 0$ almost surely. Thus we have
    \[
        \mathbb{E}\widehat{S}_n^2\!\leq\!\sigma_0^2\!\left(\!1\!-\!\frac{1}{n}\!\right)\! +\! \frac{4\mathbb{E}\|Z^c\|^4}{n^2} \!-\! \frac{4\mathbb{E}\|Z^c\|^4}{n}\!+\!\frac{4\mathbb{E}\|Z^c\|^4}{n^2}\!+\!\frac{4\mathbb{E}\|Z^c\|^4}{n}
        =\sigma_0^2\!\left(\!1\!-\!\frac{1}{n}\!\right)\! +\! \frac{8\mathbb{E}\|Z^c\|^4}{n^2},
    \]
    and the proof is complete.
\end{proof}

\subsection{Proof of Corollary \ref{COR:PVE_Consistency}}

We know from previous arguments that $\mathbb{P}_n \|Z- \bar{Z}_n \|^2 \xrightarrow{as\ast} \mathbb{E}\|Z^c\|^2 > 0$. We also know by Theorem~\ref{THM:maxes} that $\widehat{\lambda}_j \rightarrow \lambda_j$ for all $1 \leq j \leq m+1$. If $m_0(\alpha)=1$, then $\lambda_1 > (1-\alpha)\mathbb{E}\|Z^c\|^2$ and $\mathbbm{1}\{\widehat{\lambda}_1 > (1-\alpha)\mathbb{P}_n\|Z-\bar{Z}_n\|^2\} \xrightarrow{as\ast} 1$ as $n \rightarrow \infty$. If $m_0(\alpha) > 1$, then, for $m = m_0(\alpha)$, we have $\mathbbm{1}\{\sum_{j=1}^{m-1} \widehat{\lambda}_j < (1-\alpha)\mathbb{P}_n\|Z - \bar{Z}_n\|^2 < \sum_{j=1}^m \widehat{\lambda}_j\} \xrightarrow{as\ast} 1$, as $n \rightarrow \infty$, and the desired result is proved.

\subsection{Proof of Theorem \ref{THM:LM_Consistency}}

Let $\widehat{\theta}_n^\ast = \begin{pmatrix} \widehat\alpha_n^\ast & \widehat\beta_n^\ast & \widehat\gamma_n^\ast \end{pmatrix}^\top = \{P(\widehat{U}\widehat{U}^\top)\}^{-1}P(\widehat{U}Y)$. Then
    \begin{align*}
    \sqrt{n}(\widehat{\theta}_n - \widehat{\theta}_n^\ast)
    &= \sqrt{n} \left\{\mathbb{P}_n\left(\widehat{U}\widehat{U}^\top\right)\right\}^{-1}\mathbb{P}_n\left[\widehat{U}\left\{Y-(\widehat{\theta}_n^{\ast})^{\top}\widehat{U}\right\}\right]\\
    &= \left\{\mathbb{P}_n\left(\widehat{U}\widehat{U}^\top\right)\right\}^{-1}\sqrt{n}\left(\mathbb{P}_n-P\right)\left[\widehat{U}\left\{Y-(\widehat{\theta}_n^\ast)^\top\widehat{U}\right\}\right]\\
    &= \left\{P\left(UU^\top\right)\right\}^{-1}\sqrt{n}\left(\mathbb{P}_n-P\right)\left\{U(Y-\theta_0^\top U)\right\} +o_P(1),
    \end{align*}
    where, via recycling previous arguments, the first equality is almost surely true for all $n$ large enough since the first term is almost surely consistent for a full rank matrix. The next two equalities follow from algebra, and the last equality follows from recycled empirical process arguments, including reuse of Corollary~\ref{COR2:product}. Next, we have
    \begin{align*}
        &\sqrt{n}(\widehat{\theta}_n^\ast - \theta_0)
        = \left\{P\left(\widehat{U}\widehat{U}^\top\right)\right\}^{-1}\sqrt{n}P\left\{\widehat{U}(Y-\theta_0^\top\widehat{U})\right\}\\
        &=\left\{P\left(\widehat{U}\widehat{U}^\top\right)\right\}^{-1}\left[P\left\{\begin{pmatrix}
        0 \\ 0 \\ \sqrt{n}(\widehat{Z}_\ast-Z_\ast)
    \end{pmatrix}(Y - \theta_0^\top\widehat{U})\right\}-P\{U\sqrt{n}(\widehat{Z}_\ast-Z_\ast)^\top\gamma_0\}\right]\\
    &=\left\{P\left(UU^\top\right)\right\}^{-1}\left[P\left\{\begin{pmatrix}
        0 \\ 0 \\ \sqrt{n}(\widehat{Z}_\ast-Z_\ast)
    \end{pmatrix}(Y - \theta_0^\top\widehat{U})\right\}-P\{U\sqrt{n}(\widehat{Z}_\ast-Z_\ast)^\top\gamma_0\}\right] + o_P(1)
    \equiv C_1 + o_P(1).
    \end{align*}
    The remainder of the proof follows by first recalling that
    \[
    \sqrt{n}(\widehat{Z}_{\ast}-Z_{\ast})=\begin{pmatrix}\langle \sqrt{n}(\widehat{\phi}_1-\phi_1),Z\rangle\\ \vdots \\ \langle \sqrt{n}(\widehat{\phi}_m-\phi_m),Z\rangle\end{pmatrix},
    \]
    applying some expectations, and then applying some algebraic derivations. This completes the proof.

\subsection{Proof of Theorem \ref{THM:bootstrap_v}}

For any $a_1, a_2 \in \mathcal{B}$, we have
    \begin{align*}
    &a_1^\top\left[\sqrt{n}(\widetilde{V}_n^{\ast} - \widehat{V}_n) - \widetilde{\mathbb{G}}_n\{Z^c{Z^c}^\top\}\right]a_2\\
    &= a_1^\top n^{-1/2}\sum_{i=1}^n\left(W_{ni} -1\right) \Big[\mathrm{G}_N^\top Z_i^c{Z_i^c}^\top\mathrm{G}_N - Z_i^c {Z_i^c}^\top - \mathrm{G}_N^\top\left\{\left(\widetilde{Z}_n^{\ast}-\mu\right)\left(\widetilde{\mathbb{G}}_n Z\right)+\left(\widetilde{\mathbb{G}}_n Z\right)(\widetilde{Z}_n^{\ast} - \mu)^{\top}\right\}\mathrm{G}_N\Big]a_2 \\
    &= a_1^\top(A_n - B_n)a_2,
    \end{align*}
    where $A_n = n^{1/2}\sum_{i=1}^n(W_{ni} -1)(\mathrm{G}_N^\top Z_i^c{Z_i^c}^\top\mathrm{G}_N - Z_i^c {Z_i^c}^\top)$ and $B_n = -n^{1/2}\sum_{i=1}^n (W_{ni}-1)\mathrm{G}_N^\top\{(\widetilde{Z}_n^\ast - \mu)(\widetilde{\mathbb{G}}_nZ)+(\widetilde{\mathbb{G}}_nZ)(\widetilde{Z}_n^\ast - \mu)^\top\}\mathrm{G}_N$.
    Note that $\widetilde{Z}_n^{\ast} - \mu = \widetilde{Z}_n^{\ast} - \bar{Z}_n + \bar{Z}_n - \mu$; this implies that
    \[\Pr(| \widetilde{Z}_n^{\ast} - \mu | > \epsilon| Z_1, \cdots, Z_n) \leq \Pr(|\widetilde{Z}_n^{\ast} - \bar{Z}_n| > \epsilon/2| Z_1, \cdots, Z_n) + \mathbbm{1}(|\bar{Z}_n - \mu| > \epsilon/2).\]
    As $\Pr(|\widetilde{Z}_n^{\ast} - \bar{Z}_n| > \epsilon/2| Z_1, \cdots, Z_n) \xrsquigarrow{P}{Z_1, \cdots, Z_n} 0$ and $\mathbbm{1}(|\bar{Z}_n - \mu| > \epsilon/2) \xrightarrow{P} 0$, we have $\Pr(| \widetilde{Z}_n^{\ast} - \mu | > \epsilon| Z_1, \cdots, Z_n)\xrsquigarrow{P}{Z_1, \cdots, Z_n} 0$. Thus, we can proceed with $A_n$.

    Now, $|a_1^\top A_n a_2| = |a_1^\top\{n^{1/2}\sum_{i=1}^n(W_{ni} -1)(\mathrm{G}_N^\top Z_i^c{Z_i^c}^\top\mathrm{G}_N - Z_i^c{Z_i^c}^\top)\}a_2|$. Since $\mathbb{E}(|W_{ni}-1| | Z_1, \cdots, Z_n) \leq \mathbb{E}(W_{ni}+1 | Z_1, \cdots, Z_n) = 2$, this implies that
    \begin{align*}
    \mathbb{E}(a_1^\top A_n a_2 | Z_1, \cdots, Z_n) &\leq 4\sqrt{n}\left\{n^{-1}\sum_{i=1}^n \left(\|Z_i\|\|\mathrm{G}_N^\top Z_i - Z_i \|\right)\right\} \equiv C_1,
    \end{align*}
    since $\|\mathrm{G}_N Z_i \| \leq \|Z_i \|$. Thus
    \begin{align*}
    C_1 &\leq 4\left(n^{-1}\sum_{i=1}^n \|Z_i \|^2\right)^{1/2}\sqrt{n}\left(n^{-1}\sum_{i=1}^n \|\mathrm{G}_N Z_i - Z_i \|^2\right)^{1/2}
    = 4\left\{\mathbb{E}\|Z\|^2 + o_P(1)\right\}^{1/2}\left(\sum_{i=1}^n \|\mathrm{G}_N Z_i - Z_i \|^2\right)^{1/2} =o_P(1),
    \end{align*}
    since $\mathbb{E}\|\mathrm{G}_N Z_i - Z_i \|^2 = o(n^{-1})$. This implies that $\sup_{a_1, a_2 \in \mathcal{B}} |a_1^\top A_n a_2 | \xrsquigarrow{P}{W} 0$, and thus the proof is complete.

\subsection{Proof of Theorem \ref{THM:bootstrap_v2}}

The bootstrap and measurability conditions of Theorem~12.1 of \citet{Kosorok:08} are satisfied by $\widetilde{V}_n^{\ast}$, since $\sqrt{n}(\widetilde{V}_n^{\ast} - \widehat{V}_n) \xrsquigarrow{P}{W} \mathbb{G}_n(Z^c{Z^c}^\top) \rightsquigarrow \mathbb{X}$ which is tight (since the Donsker class has bootstrap validity by Theorem~2.6 of \citet{Kosorok:08}). Hence, the results given in the proof of Theorem~12.1 of \citet{Kosorok:08} apply. Namely, expression~(12.1) of \citet{Kosorok:08} yields the desired result.

\subsection{Proof of Theorem \ref{THM:bootstrap}}

Let $\widetilde{\theta}_n^\ast = \begin{pmatrix} \widetilde\alpha_n^\ast & \widetilde\beta_n^{\ast} & \widetilde\gamma_n^\ast \end{pmatrix}^\top = \left[\mathbb{P}_n\left\{\begin{pmatrix} 1 & \widetilde{X}^\top & \widetilde{Z}_\ast^\top \end{pmatrix}^{\top}\right\}^{\otimes 2}\right]^{-1}\mathbb{P}_n\left\{\begin{pmatrix} 1 & \widetilde{X}^\top & \widetilde{Z}_\ast^\top \end{pmatrix}^\top Y\right\}$.
Then
    \begin{align*}
    &\sqrt{n}(\widetilde{\theta}_n - \widetilde{\theta}_n^\ast) = \left\{\mathbb{P}_n\begin{pmatrix} 1 \\ \widetilde{X} \\ \widetilde{Z}_\ast \end{pmatrix}^{\otimes 2}\right\}^{-1}\sqrt{n}(\widetilde{\mathbb{P}}_n - \mathbb{P}_n)\left[\begin{pmatrix} 1 \\ \widetilde{X} \\ \widetilde{Z}_\ast \end{pmatrix}\left\{Y - \widetilde{\alpha}_n^\ast -(\widetilde{\beta}_n^\ast)^\top \widetilde{X} -(\widetilde{\gamma}_n^\ast)^\top \widetilde{Z}_\ast\right\}\right]\\
    &= \left\{P\begin{pmatrix} 1 \\ X \\ Z_\ast \end{pmatrix}^{\otimes 2}\right\}^{-1}\sqrt{n}(\widetilde{\mathbb{P}}_n - \mathbb{P}_n)\left\{\begin{pmatrix} 1 \\ X \\ Z_\ast \end{pmatrix}\left(Y - \alpha_0 -\beta_0^\top X -\gamma_0^\top Z_\ast\right)\right\} +E_n,
    \end{align*}
    where $E_n \xrsquigarrow{P}{W} 0$, recycling previous arguments. Next,
    \begin{align*}
        &\sqrt{n}(\widetilde{\theta}_n^\ast - \widehat{\theta}_n)
        = \left\{P\begin{pmatrix} 1 \\ X \\ Z_\ast \end{pmatrix}^{\otimes 2}\right\}^{-1}\!\Bigg[P\left\{\begin{pmatrix} 0 \\ 0 \\ \sqrt{n}(\widetilde{Z}_\ast - \widehat{Z}_\ast) \end{pmatrix}\!\left(Y - \alpha_0 - \beta_0^\top X - \gamma_0 Z_{\ast}\right)\right\} \\
        &\qquad\qquad\qquad\qquad\qquad\qquad\qquad\qquad - P\left\{\begin{pmatrix} 1 \\ X \\ Z_{\ast} \end{pmatrix}\sqrt{n}(\widetilde{Z}_\ast-\widehat{Z}_\ast)^\top\gamma_0\right\}\Bigg] + E_{n1} \equiv C_2,
    \end{align*}
    where $E_{n1} \xrsquigarrow{P}{W} 0$, recycling previous arguments, and thus $C_2 = L_0\{\widetilde{\mathbb{G}}_n({Z^c}{Z^c}^\top)\} + E_{n2}$, where $E_{n2} \xrsquigarrow{P}{W} 0$. Thus the proof is complete.

\end{document}